\documentclass[12pt, reqno]{amsart}
\usepackage[T1]{fontenc}
\usepackage[utf8]{inputenc}
\usepackage{lmodern,amsmath,amssymb,amsthm,mathtools,mathrsfs}
\usepackage{enumitem}
\usepackage[expansion=false]{microtype}
\usepackage{graphicx,tikz}
\usepackage[hidelinks,bookmarksnumbered]{hyperref}
\hypersetup{colorlinks=true,linkcolor=blue, anchorcolor=blue,
    citecolor=cyan, urlcolor=blue, filecolor=magenta, pdftoolbar=true}

\def\XXint#1#2#3{{\setbox0=\hbox{$#1{#2#3}{\int}$ }
        \vcenter{\hbox{$#2#3$ }}\kern-.6\wd0}}

\newtheorem{theorem}{Theorem}[section]
\newtheorem{lemma}[theorem]{Lemma}
\newtheorem{proposition}[theorem]{Proposition}
\newtheorem{corollary}[theorem]{Corollary}
\theoremstyle{definition}
\newtheorem{definition}[theorem]{Definition}

\theoremstyle{remark}
\newtheorem{remark}[theorem]{Remark}
\numberwithin{equation}{section}
\newcommand{\R}{\mathbb R}
\newcommand{\N}{\mathbb N}
\newcommand{\Z}{\mathbb Z}

\newcommand{\supp}{\operatorname{supp}}
\newcommand{\BMO}{\operatorname{BMO}}
\newcommand{\WBP}{\operatorname{WBP}}

\begin{document}

    \title[Heat kernels and a $T1$ theorem]{Kernel Estimates and   Calder\'{o}n-Zygmund Theory in the $(k,a)$-generalized Laguerre setting}

    \author{Mingwei Shi}
    \address{School of Mathematics and Physics, University of Science and Technology Beijing, Beijing 100083, China}
    \email{shimingweimath@163.com}

\author{Yongsheng Han}
\address{Department of Mathematics, Auburn University, AL 36849-5310, USA}
\email{hanyong@auburn.edu}

    \author{Yu Liu}
    \address{School of Mathematics and Physics, University of Science and Technology Beijing, Beijing 100083, China}
    \email{liuyu75@pku.org.cn}

    \thanks {The corresponding author
     is Yu Liu, E-mail: liuyu75@pku.org.cn}

    \maketitle

    {\Small {\bf Abstract:}
In this paper, we establish heat and Poisson kernel estimates for
$\mathfrak L=\|x\|^{2-a}\Delta_k$ in the $(k,a)$-generalized
Laguerre setting, where $\Delta_k$
   is the Dunkl Laplacian associated with a finite reflection group.
   The heat kernel is nonnegative, symmetric, and conservative, and satisfies a Gaussian bound in the orbit distance
    together with polynomial decay in the distance between the transformed coordinates. The Poisson kernel, obtained by subordination, satisfies corresponding
    size and time-derivative estimates, while $-t\partial_t p_t$ has integral
    zero in each variable.
Motivated by these estimates, we introduce a class of
$(k,a)$-Calder\'{o}n-Zygmund singular integral operators and show
that $L^2_{k,a}$-boundedness implies $L^p_{k,a}$-boundedness for
$1<p<\infty$, together with weak type $(1,1)$ and
$L^\infty$-to-$\BMO$ estimates. Finally,
      we establish a $T1$  theorem characterizing the
       $L^2$-boundedness of such operators by the conditions $T1\in BMO$, $T^*1\in BMO$, and the weak boundedness property.

    \vspace{0.2cm} {\bf Keywords:} { heat kernel, Poisson kernel,
      Calder\'{o}n-Zygmund singular integral operator, $T1$ theorem}

    }

\tableofcontents
\section{Introduction}\label{sec:intro}
Dunkl theory extends classical Euclidean Fourier analysis by replacing
Lebesgue measure with a weighted measure invariant under a finite
reflection group and partial derivatives with differential-difference
operators. The theory originated in the work of Dunkl \cite{Dunkl} and
has led to developments in orthogonal polynomials, integral transforms,
heat equations, and singular integrals; see
\cite{ADH,A1,A4,DX,Rosler,RV}.

The \((k,a)\)-generalized Fourier transform is a two-parameter family
introduced by Ben Sa\"id, Kobayashi, and {\O}rsted \cite{BKOshort,BKO}.
It contains the Euclidean Fourier transform when \(k=0\) and \(a=2\),
the Dunkl transform when \(a=2\) \cite{A25}, and a Hankel-type transform when
\(a=1\). The associated holomorphic Laguerre semigroup acts on a
weighted space \(L^2_{k,a}(\mathbb{R}^n)\) and unifies the corresponding Hermite and
Laguerre constructions. We recall the differential-difference
operators
\begin{equation}\label{eq:operators}
    \Delta_{k,a}=\mathfrak L-\|x\|^a \qquad \mathrm{ and}
    \qquad \mathfrak L=\|x\|^{2-a}\Delta_k,
\end{equation} where $\Delta_k$
   is the Dunkl Laplacian associated with a finite reflection group.
The heat semigroup studied in this paper is generated by \(\mathfrak
L\), not by the oscillator \(\Delta_{k,a}\).

Weighted inequalities for the deformed transform were studied in
\cite{A9}. Kernel bounds and parameter-dependent behavior were
investigated in \cite{A27,A11}. Positive radial product formulas and
convolution-based maximal estimates in Dunkl analysis were treated
in \cite{A16,A20}, while the \((k,1)\)-translation and the maximal
operator were studied in \cite{A26}.

When \(k=0\) and \(a=2\), the operator \(\mathfrak L\) is the
Euclidean Laplacian. The Dunkl heat equation and its positive heat
kernel were studied by R\"osler \cite{Rosler}. Detailed heat and
Poisson estimates in the rational Dunkl setting were subsequently
established in \cite{ADH,DH}. For \(k=0\) and general \(a\), Aoyama
\cite{Aoyama} constructed a deformed heat kernel and studied its
heat flow, integral identities, and the associated Markov process.
Gaussian-semigroup Hardy-space methods are developed in \cite{A28},
and Riesz-transform boundedness in the Dunkl setting is studied
in \cite{A122}. %The specific relationship mentioned above is shown in the following figure.

Motivated by these results, we introduce the heat kernel associated with the \((k,a)\)-generalized Laplace operator \(\mathfrak L\)
 (see Definition~\ref{def:heat}). A natural next question is to establish an upper bound that retains Gaussian decay
 in the orbit distance while also reflecting the distance between the transformed coordinates. In the Dunkl case
 where \(a=2\), the estimates in \cite{DH} were obtained from first- and second-order identities for the
  heat kernel. For general \(a\), we instead estimate the second-order displacement moment. The constant
   and first-order spherical harmonic components reduce the calculation to radial Bessel kernels.
   Combining the moment bound with the semigroup identity and an initial orbit-Gaussian estimate
   yields additional polynomial decay in the coordinate distance, without using either a \((k,a)\)-translation formula
    or spatial derivative formulas for \(B_{k,a}\).

Throughout this paper, we assume that the parameters satisfy \eqref{eq:parameters-minimal}.
 Under this assumption, we obtain the following heat kernel estimates.

\begin{theorem}\label{thm:heat}
    For every $t>0$, the heat kernel is nonnegative and symmetric:
    \[
    h_t(x,y)\geq0,\qquad h_t(x,y)=h_t(y,x)
    \]
    for $dw_{k,a}(x)\,dw_{k,a}(y)$-almost every $(x,y)\in\R^n\times\R^n$.

    It is conservative, in the sense that
    \[
    \int_{\R^n}h_t(x,y)\,dw_{k,a}(y)=1
    \]
    for $dw_{k,a}$-almost every $x\in\R^n$.
Moreover, there exists a constant \(C>0\) such that, for every \(t>0\),
\begin{equation}\label{eq:moment}
    \int_{\mathbb{R}^n}\|\widehat{x}-\widehat{y}\|^2h_t(x,y)\,dw_{k,a}(y)\leq Ct
\end{equation}
for \(dw_{k,a}\)-almost every \(x\in\mathbb{R}^n\), where
\(\widehat{x}\) and \(\widehat{y}\) are given as in \eqref{eq:phi}.

    For every $0<\sigma<2$, there exist constants
    $C_\sigma,c_\sigma>0$ such that, for every $t>0$,
    \begin{equation}\label{eq:heat-improved}
        0\leq h_t(x,y)\leq\frac{C_\sigma}{V(x,y,\sqrt t)}
        \left(1+\frac{\|\widehat x-\widehat y\|}{\sqrt t}\right)^{-\sigma}
        e^{-c_\sigma d(x,y)^2/t}
    \end{equation}
    for $dw_{k,a}(x)\,dw_{k,a}(y)$-almost every $(x,y)\in\R^n\times\R^n$,  where \(\widehat{x}\) and \(\widehat{y}\) are given as  in \eqref{eq:phi}.
\end{theorem}

Starting from the heat kernel, we introduce the Poisson semigroup
$P_t=e^{-t\sqrt{-\mathfrak L}}$ and its kernel $p_t(x,y)$
(see Definition~\ref{def:poisson}). By subordination, we obtain
the following size and time-derivative estimates, together with
cancellation identities.

\begin{theorem}\label{thm:poisson}
    Let \(m\in\N_0\) and \(0<\sigma<\min\{2,1+\min(n,n_{k,a})\}\). There exists a constant \(C_{m,\sigma}>0\) such that, for every \(t>0\),
    \begin{equation}\label{eq:poisson-improved}
        |t^m\partial_t^m p_t(x,y)|\leq
        \frac{C_{m,\sigma}}{V(x,y,t+d(x,y))}\frac{t}{t+d(x,y)}
        \left(\frac{t+d(x,y)}{t+\|\widehat x-\widehat y\|}\right)^\sigma
    \end{equation}
    for \(dw_{k,a}(x)\,dw_{k,a}(y)\)-almost every \((x,y)\in\R^n\times\R^n\),  where \(\widehat{x}\) and \(\widehat{y}\) are given as in \eqref{eq:phi}.
  In particular, there exists a constant \(C_1>0\) such that
    \begin{equation}\label{eq:poisson-one}
        p_t(x,y)\leq\frac{1}{V(x,y,t+d(x,y))}
        \frac{C_1t}{t+\|\widehat x-\widehat y\|}
    \end{equation}
    for every \(t>0\) and \(dw_{k,a}(x)\,dw_{k,a}(y)\)-almost every \((x,y)\in\R^n\times\R^n\), where \(\widehat{x}\) and \(\widehat{y}\) are given
    as in \eqref{eq:phi}.

    For every \(t>0\), the kernels
    \(q_t=-t\partial_t p_t\) have integral zero in each variable for
    \(dw_{k,a}\)-almost every value of the other variable. Their associated
    operators are bounded on \(L^p_{k,a}(\mathbb{R}^n)\),
    \(1\leq p\leq\infty\), with bounds independent of \(t>0\).
\end{theorem}

Motivated by the Poisson kernel estimates, we introduce and study
a class of $(k,a)$-Calder\'{o}n-Zygmund singular integral operators
in the $(k,a)$-generalized Laguerre setting.
The Poisson estimates provide the size bounds motivating this
class; the spatial difference conditions are imposed explicitly
in Definition~\ref{def:CZ}.

\begin{definition}\label{def:CZ}
   Let \(n\geq1\), \(a>0\), \(k\geq0\), and \(n+2\gamma+a-2>0\),
   and let \(\widehat{x}\) and \(\widehat{y}\) be as in \eqref{eq:phi}. A
   continuous linear operator
   \(T:\mathscr D_a\longrightarrow\mathscr D_a'\), where
   $\mathscr D_a$ is defined in Section~\ref{sec:analytic-symbols},
   is called a
   \((k,a)\)-Calder\'{o}n-Zygmund singular integral operator if its
   kernel \(K(x,y)\) satisfies, for some
   \(0<\varepsilon\leq1\) and \(0 < C_K<\infty\), the following estimates
   for pairs satisfying \(d(x,y)>0\):
    \begin{equation}\label{eq:CZ-size}
        |K(x,y)|\leq\frac{C_K}{w_{k,a}(B_a(x,d(x,y)))}
        \left(\frac{d(x,y)}{\|\widehat x-\widehat y\|}\right)^\varepsilon,\qquad d(x,y)>0,
    \end{equation}
    and
    \begin{equation}\label{eq:CZ-y}
        |K(x,y)-K(x,y')|
        \leq C_K\left(\frac{\|\widehat y-\widehat{y'}\|}{\|\widehat x-\widehat y\|}\right)^\varepsilon
        \frac1{w_{k,a}(B_a(x,d(x,y)))},
    \end{equation}
    whenever $\|\widehat y-\widehat{y'}\|\leq d(x,y)/2$, and
    \begin{equation}\label{eq:CZ-x}
        |K(x',y)-K(x,y)|
        \leq C_K\left(\frac{\|\widehat x-\widehat{x'}\|}{\|\widehat x-\widehat y\|}\right)^\varepsilon
        \frac1{w_{k,a}(B_a(x,d(x,y)))},
    \end{equation}
    whenever $\|\widehat x-\widehat{x'}\|\leq d(x,y)/2$.

    Moreover, $K$ represents $T$ off the ordinary diagonal:
    \begin{equation}\label{eq:CZ-representation}
        Tf(x)=\int_{\mathbb{R}^n}K(x,y)f(y)\,dw_{k,a}(y),
        \qquad x\notin\supp f,
    \end{equation}
    where the equality is understood almost everywhere and
    distributionally off $\supp f$. Equivalently, for disjointly
    supported tests $f$ and $g$,
    \begin{equation}\label{eq:offdiagonal-pairing}
        \langle Tf,g\rangle
        =\iint_{\mathbb R^n\times\mathbb R^n} K(x,y)f(y)\overline{g(x)}\,dw_{k,a}(y)\,dw_{k,a}(x).
    \end{equation}
    A $(k,a)$-Calder\'{o}n-Zygmund singular integral operator is called
    a $(k,a)$-Calder\'{o}n-Zygmund operator if it extends to a bounded
    operator on $L^2_{k,a}(\mathbb{R}^n)$.
\end{definition}

For a $(k,a)$-Calder\'{o}n-Zygmund operator, we obtain the
following boundedness properties.

\begin{theorem}\label{thm:Lp}
   Let $T$ be a $(k,a)$-Calder\'{o}n-Zygmund operator as in
   Definition~\ref{def:CZ}. Then
    \[
    T:L^p_{k,a}(\mathbb{R}^n)\longrightarrow L^p_{k,a}(\mathbb{R}^n)
    \quad(1<p<\infty),
    \qquad
    T:L^1_{k,a}(\mathbb{R}^n)\longrightarrow L^{1,\infty}_{k,a}(\mathbb{R}^n),
    \]
    and
    \[
    T:L^\infty_{k,a}(\mathbb{R}^n)\longrightarrow
    \BMO(\mathbb{R}^n,dw_{k,a}).
    \]
    Moreover, the operator norms depend only on the  parameters,
    $\|T\|_{2\to2}$, and the kernel constant $C_K$; the
    $L^p_{k,a}(\mathbb R^n)$ operator norm may also depend on $p$.
\end{theorem}

The $T1$ theorem in the rational Dunkl setting was established
in \cite{THHLL}. We formulate the corresponding criterion for
the coordinate map $\phi$ and the measure $dw_{k,a}$ considered
here. Its proof combines the mixed-distance molecular estimates
with reproducing formulae from the general RD-space theory.

\begin{theorem}\label{thm:T1-intro}
    Let $T$ be a $(k,a)$-Calder\'{o}n-Zygmund singular integral
    operator as in Definition~\ref{def:CZ}. Then $T$ extends to a bounded
    operator on $L^2_{k,a}(\mathbb{R}^n)$ if and only if the following
    three conditions hold:
    \begin{enumerate}[label=\textup{(\roman*)},ref=\thetheorem(\roman*),leftmargin=2.2em,labelsep=0.55em,itemsep=6pt,topsep=6pt,parsep=0pt]
        \item\label{eq:T1-iff} $T1\in\BMO(\mathbb R^n,dw_{k,a})$;
        \item $T^*1\in\BMO(\mathbb R^n,dw_{k,a})$;
        \item $T\in\WBP$.
    \end{enumerate}

Moreover,
    \begin{equation}\label{eq:T1-bound}
        \|T\|_{2\to2}\leq C\bigl(C_K+\|T\|_{\WBP}
        +\|T1\|_{\BMO}+\|T^*1\|_{\BMO}\bigr),
    \end{equation}
    where $C$ depends only on the parameters.
\end{theorem}

The kernel conditions in Definition~\ref{def:CZ} involve both
the coordinate distance $\|\widehat x-\widehat y\|$ and the orbit
distance $d(x,y)$. These distances can differ substantially near
reflected diagonals, so the class is not covered directly by the
classical coordinate-metric Calder\'on--Zygmund estimates.
The coordinate reproducing formulae provide auxiliary decompositions;
the mixed-distance estimates are handled by the molecular and
almost-orthogonality arguments in Section~\ref{sec:T1}.

The paper is organized as follows. Section~\ref{sec:prelim} recalls
the Dunkl operators, weighted function spaces, spherical harmonics,
and generalized Fourier transform. Section~\ref{sec:kernels}
establishes the heat and Poisson kernel estimates.
Section~\ref{sec:CZ} presents the kernel integration calculation and
proves Theorem~\ref{thm:Lp}. Section~\ref{sec:T1} proves the $T1$
criterion via Meyer's commutation argument, smooth molecules, almost
orthogonality, reproducing formula, and paraproducts.

\subsection{Notation}\label{sec:notation}
We use the following notation throughout the paper.
\begin{itemize}[leftmargin=*,itemsep=4pt,topsep=6pt]
    \item We write $\mathbb{N}=\{1,2,\ldots\}$ and
    $\mathbb{N}_0=\{0,1,2,\ldots\}$.
    \item $\langle x,y\rangle$ and $\|x\|$ denote the Euclidean inner
    product and norm. We write $B(z,r)=\{u:\|u-z\|<r\}$ and
    $S^{n-1}=\{u:\|u\|=1\}$.
    \item $\mathbf1_E$ denotes the indicator of a measurable set $E$,
    $\supp f$ the support of $f$, and $\overline f$ denotes the complex conjugate of $f$.
    \item For $\alpha=(\alpha_1,\ldots,\alpha_n)\in\mathbb{N}_0^n$,
    we write $|\alpha|=\sum_{j=1}^n\alpha_j$ and
    $\partial^\alpha=\partial_{x_1}^{\alpha_1}\cdots
    \partial_{x_n}^{\alpha_n}$.
    \item Let $C, c > 0$ denote positive constants that may vary from line to line.
    Unless otherwise stated, they depend only on fixed parameters and exponents.
     We write $A\lesssim B$ if
    $A\leq CB$, and $A\asymp B$ if $A\leq CB$ and $B\leq CA$.
    \item
    Operator pairings $\langle Tf,g\rangle$ are linear in $f$ and
    conjugate-linear in $g$; $T^*$ denotes the adjoint and $I$ the
    identity operator.
    \item Denote by $C^\infty_c(\mathbb{R}^n)$ the space of all compactly supported infinitely differentiable functions on
    $\mathbb{R}^n$.
\end{itemize}

\section{Preliminaries}\label{sec:prelim}\label{sec2}
\subsection{Dunkl operators }
This section reviews some fundamental concepts related to the theory of Dunkl operators (see \cite{BKO,A1,Dunkl,A4,DX,Rosler}).

Let $\langle x, y\rangle = \sum_{j=1}^n x_j y_j$. For $\alpha \in \mathbb{R}^n\setminus\{0\}$, let $\sigma_\alpha$ denote the reflection across the hyperplane $H_\alpha \subset \mathbb{R}^n$ orthogonal to $\alpha$, i.e.,
\[
\sigma_\alpha(x) = x - 2 \frac{\langle\alpha, x\rangle}{\|\alpha\|^2} \alpha, \quad x \in \mathbb{R}^n,
\]
where $\|\alpha\| := \sqrt{\langle \alpha, \alpha\rangle}$. Each reflection $\sigma_\alpha$ is orthogonal with respect to the standard inner product. A finite set $R \subset \mathbb{R}^n\setminus\{0\}$ is called a root system if $\sigma_\alpha(R)=R$ for all $\alpha\in R$ and $R\cap \mathbb{R}\cdot \alpha = \{\pm\alpha\}$. The reflection group (or Coxeter group) $G \subset O(n)$ generated by $\{\sigma_\alpha\}_{\alpha\in R}$ is the group associated with $R$. For any $\beta \in \mathbb{R}^n\setminus \bigcup_{\alpha\in R} H_\alpha$, there exists a positive root system $R_{+} = \{\alpha \in R : \langle\alpha, \beta\rangle > 0\}$ such that for every $\alpha\in R$, either $\alpha\in R_{+}$ or $-\alpha\in R_{+}$.

Throughout the paper, $k:R\to[0,\infty)$ is a $G$-invariant
multiplicity function, that is, $k(g\alpha)=k(\alpha)$ for every
$g\in G$ and $\alpha\in R$. The Dunkl operators $T_\xi$ are the
following deformations of directional derivatives $\partial_\xi$ by
difference operators:
\begin{equation}\label{eq:Dunkl}
    T_\xi f(x)=\partial_\xi f(x)
    +\sum_{\alpha\in R_+}k(\alpha)\langle\alpha,\xi\rangle
    \frac{f(x)-f(\sigma_\alpha x)}{\langle\alpha,x\rangle}.
\end{equation}

The Dunkl operators commute with each other.
If either $f$ or $g$ is $G$-invariant, the product rule reduces to the standard Leibniz rule:

\begin{equation}\label{eq:product}
    T_\xi(fg)=(T_\xi f)g+fT_\xi g.
\end{equation}
We use  $e_j$ to denote the $j$th standard unit vector; moreover,
$\partial_\xi=\sum_{j=1}^n\xi_j\partial_{x_j}$,
$\nabla=(\partial_{x_1},\ldots,\partial_{x_n})$, and
$\Delta=\sum_{j=1}^n\partial_{x_j}^2$ are the ordinary directional
derivative, gradient, and Laplacian, respectively. Writing
$T_j=T_{e_j}$, we have
\begin{equation}\label{eq:Dunkl-Laplacian}
    \begin{aligned}
        \Delta_k f(x)&=\sum_{j=1}^nT_j^2f(x)\\
        &=\Delta f(x)+2\sum_{\alpha\in R_+}k(\alpha)
        \left(\frac{\langle\nabla f(x),\alpha\rangle}{\langle\alpha,x\rangle}
        -\frac{\|\alpha\|^2}{2}
        \frac{f(x)-f(\sigma_\alpha x)}{\langle\alpha,x\rangle^2}\right).
    \end{aligned}
\end{equation}
The $(k,a)$-measure is denoted  as
\begin{equation}\label{eq:weight}
    dw_{k,a}(x)=\|x\|^{a-2}
    \prod_{\alpha\in R_+}|\langle\alpha,x\rangle|^{2k(\alpha)}\,dx,
\end{equation} and let $
 w_{k,a}(E)=\int_Edw_{k,a}(x).$   The weighted Lebesgue spaces are denoted by
\[
L^p_{k,a}(\mathbb{R}^n)=L^p(\mathbb{R}^n,dw_{k,a}),
\qquad 1\leq p\leq\infty.
\]
For $1\leq p<\infty$, the norm is
\begin{equation}\label{eq:Lp-norm-notation}
    \|f\|_p:=\|f\|_{L^p_{k,a}(\mathbb{R}^n)}
    =\left(\int_{\mathbb{R}^n}|f(x)|^p\,dw_{k,a}(x)\right)^{1/p},
\end{equation}
and
\[
\|f\|_\infty:=\|f\|_{L^\infty_{k,a}(\mathbb{R}^n)}
=\operatorname*{ess\,sup}_{x\in\mathbb{R}^n}|f(x)|.
\]
The weak space $L^{1,\infty}_{k,a}(\mathbb{R}^n)$ consists of those
$f$ for which
\begin{equation}\label{eq:weak-norm-notation}
    \|f\|_{L^{1,\infty}_{k,a}(\mathbb{R}^n)}
    =\sup_{\lambda>0}\lambda\,
    w_{k,a}\{x\in\mathbb{R}^n:|f(x)|>\lambda\}<\infty.
\end{equation}

Throughout we impose the standing conditions:
\begin{equation}\label{eq:parameters-minimal}
    n\geq1,\qquad a>0,\qquad k\geq0,\qquad
    n+2\gamma+a-2>0,
    \qquad \gamma=\sum_{\alpha\in R_+}k(\alpha).
\end{equation}
For $1\leq p,q\leq\infty$, the operator norm is fixed by
\begin{equation}\label{eq:operator-norm-notation}
    \|A\|_{p\to q}
    =\sup_{\|f\|_{L^p_{k,a}(\mathbb{R}^n)}\leq1}
    \|Af\|_{L^q_{k,a}(\mathbb{R}^n)}.
\end{equation}
Thus $\|A\|_{2\to2}$ always refers to the weighted space
$L^2_{k,a}(\mathbb{R}^n)$.

For $f,g\in C_c^\infty(\mathbb R^n)$,
\begin{equation}\label{eq:parts}
    \int_{\mathbb R^n}(T_jf)(x)\overline{g(x)}\,dw_{k,2}(x)
    =-\int_{\mathbb R^n}f(x)\overline{(T_jg)(x)}\,dw_{k,2}(x).
\end{equation}

For $f,g\in C_c^\infty(\mathbb R^n)$, applying the identity
\[
dw_{k,a}(x)=\|x\|^{a-2}dw_{k,2}(x)
\]
and the skew-symmetry of the Dunkl operators with respect to
$dw_{k,2}$ gives
\begin{align*}
    &\int_{\mathbb R^n}
    (\mathfrak Lf)(x)\overline{g(x)}\,dw_{k,a}(x)\\
    &\quad=
    \int_{\mathbb R^n}
    \|x\|^{2-a}(\Delta_kf)(x)\overline{g(x)}
    \|x\|^{a-2}\,dw_{k,2}(x)\\
    &\quad=
    \int_{\mathbb R^n}
    (\Delta_kf)(x)\overline{g(x)}\,dw_{k,2}(x)\\
    &\quad=
    -\sum_{j=1}^n\int_{\mathbb R^n}
    (T_jf)(x)\overline{(T_jg)(x)}\,dw_{k,2}(x)\\
    &\quad=
    \int_{\mathbb R^n}
    f(x)\overline{(\Delta_kg)(x)}\,dw_{k,2}(x)\\
    &\quad=
    \int_{\mathbb R^n}
    f(x)\overline{(\mathfrak Lg)(x)}\,dw_{k,a}(x).
\end{align*}
Thus $\mathfrak L$ is formally symmetric with respect to
$dw_{k,a}$.

The associated nonnegative energy form is
\begin{equation}\label{eq:Dirichlet}
    \begin{aligned}
        \mathcal E(f,f)={}&\int_{\R^n}|\nabla f(x)|^2\,dw_{k,2}(x)\\
        &+\frac12\sum_{\alpha\in R_+}k(\alpha)\|\alpha\|^2
        \int_{\R^n}\frac{|f(x)-f(\sigma_\alpha x)|^2}
        {\langle\alpha,x\rangle^2}\,dw_{k,2}(x).
    \end{aligned}
\end{equation}

The intertwining operator $V_k$ satisfies
$T_\xi V_k=V_k\partial_\xi$ on polynomials and admits a positive
integral representation when the multiplicities are nonnegative. This
background is recalled in \cite{A4,RV}. The operator $V_k$ determines
the ordinary Dunkl kernel $E_k$ through
\[
T_{\xi,x}E_k(x,y)=\langle\xi,y\rangle E_k(x,y)\qquad \mathrm{and}
\qquad E_k(0,y)=1.
\]

We next recall the spherical decomposition and radial formulas.
Background on orthogonal polynomials and special functions is
given in \cite{A222,DX}; the reflection-group spherical structure
originates in \cite{A1}.
Integration in polar coordinates in the original variables yields
\begin{equation}\label{eq:polar-measure}
    dw_{k,a}(r\theta)=r^{n+2\gamma+a-3}
    w_{k,2}(\theta)\,dr\,d\sigma(\theta),
    \qquad r>0,\quad\theta\in S^{n-1},
\end{equation}
where $d\sigma$ is surface measure on $S^{n-1}$. In particular, for
nonnegative measurable $f$,
\[
\int_{\mathbb R^n}f(x)\,dw_{k,a}(x)
=\int_0^\infty\int_{S^{n-1}}f(r\theta)w_{k,2}(\theta)
r^{n+2\gamma+a-3}\,d\sigma(\theta)\,dr.
\]
The angular operator $\Delta_{k,S}$ is given by
\begin{equation}\label{eq:polar-Dunkl}
    \Delta_k=\partial_r^2+\frac{n+2\gamma-1}{r}\partial_r
    +r^{-2}\Delta_{k,S}.
\end{equation}
In what follows, \(\nabla_S\) denotes the spherical gradient.
In the transformed coordinates \(\widehat{x}=\rho\theta\) defined
in \eqref{eq:phi}, we have
\[
|\nabla_{\widehat{x}}f(x)|^2
=|\partial_\rho f(x)|^2+\rho^{-2}|\nabla_S f(x)|^2.
\]

Let $\mathcal P_\ell$ be the space of ordinary homogeneous
polynomials of degree $\ell$. The restrictions to the sphere of
$\mathcal P_\ell\cap\ker\Delta_k$ form $\mathcal H_k^\ell$. These give
the orthogonal decomposition
\[
L^2(S^{n-1},w_{k,2}\,d\sigma)
=\widehat{\bigoplus}_{\ell\geq0}\mathcal H_k^\ell,
\]
and
\begin{equation}\label{eq:spherical-eigenvalue}
    \Delta_{k,S}Y=-\ell(\ell+n+2\gamma-2)Y,
    \qquad Y\in\mathcal H_k^\ell.
\end{equation}
The hat over $\bigoplus$ denotes a Hilbert-space orthogonal direct
sum. These are the $k$-spherical harmonics; see
\cite{BKO,A1,DX}.

For $h\in\mathcal P_\ell\cap\ker\Delta_k$, put
$Y=h|_{S^{n-1}}$. For a radial scalar function $v$, the radial action
recalled in \cite[Lemma~3.6]{BKO} can be written as
\begin{equation}\label{eq:radial-action-background}
    \Delta_k\{v(r)Y(\theta)\}
    =\left(v''+\frac{n+2\gamma-1}{r}v'
    -\frac{\ell(\ell+n+2\gamma-2)}{r^2}v\right)Y(\theta).
\end{equation}
Equivalently,
\[
\Delta_k(\|x\|^\nu h(x))
=\nu(\nu+2\ell+n+2\gamma-2)\|x\|^{\nu-2}h(x)
\]
away from the origin. These formulas are only used once the corresponding closed radial realizations are established.

The modified Bessel function is normalized by
\begin{equation}\label{eq:Bessel-series}
    I_\lambda(z)=\sum_{j=0}^\infty
    \frac{(z/2)^{2j+\lambda}}{j!\Gamma(j+\lambda+1)},
    \qquad z>0,\quad\lambda>-1.
\end{equation}
In the radial setting, we define the operator
$B_\zeta = \partial_\rho^2 + (\zeta-1)\rho^{-1}\partial_\rho$, where $\zeta>0$.
 Its kernel $b_t^\zeta(\rho,s)$ relative to
$s^{\zeta-1}ds$ on $(0,\infty)$ and the boundary realization are
described in Section~\ref{sec:heat-definition}.
%The factor $4t$ in
%the Laplace formula in Section~\ref{sec:kernels} comes from the
%normalization $e^{tB_\zeta}=e^{(2t)(B_\zeta/2)}$.

\subsection{Coordinates, distances and volumes}\label{sec:coordinate-geometry}

We write $x = u\theta$, where $u = \|x\|$ and $\theta \in S^{n-1}$.
The radial second-order component of $\mathfrak{L}$ is $u^{2-a}\partial_u^2$.
To normalize its diffusion coefficient, we choose an increasing radial coordinate $\rho = \rho(u)$ with $\rho(0) = 0$ satisfying
\begin{equation}\label{eq:radial-coordinate-choice}
    u^{2-a}\left(\frac{\mathrm{d}\rho}{\mathrm{d}u}\right)^2 = 1.
\end{equation}
Taking the positive square root and using $\rho(0) = 0$, we obtain
\[
\frac{\mathrm{d}\rho}{\mathrm{d}u} = u^{a/2 - 1} \qquad \mathrm{and}
\qquad \rho(u) = \int_0^u s^{a/2 - 1}\,\mathrm{d}s =
\frac{2}{a}u^{a/2}.
\]
Preserving the angular variable yields the coordinate map $\widehat{x} = \rho\theta$.
Thus, with the multiplicity sum $\gamma$ fixed above, we set
\begin{equation}\label{eq:phi}
    \widehat{x} := \phi(x) = \frac{2}{a}\|x\|^{a/2 - 1}x \qquad  \mathrm{and} \qquad \phi(0) = 0,
\end{equation}
and denote the homogeneous dimension in these coordinates by
\begin{equation}\label{eq:homogeneous-dimension}
    n_{k,a} = \frac{2(n + 2\gamma + a - 2)}{a}.
\end{equation}
When $a = 2$, the coordinate map reduces to the identity map and $n_{k,2} = n + 2\gamma$.
The natural heat kernel length scale in the transformed coordinates is $\sqrt{t}$.

The inverse coordinate map is
\[
\phi^{-1}(z) = \left(\frac{a}{2}\right)^{2/a}\|z\|^{2/a - 1}z.
\]
Accordingly, we define
\begin{equation}\label{eq:balls}
    B_a(x, r) = \phi^{-1}\big(B(\widehat{x}, r)\big)
    = \bigl\{ y \in \mathbb{R}^n : |\widehat{x} - \widehat{y}| < r \bigr\}.
\end{equation}
We write $cB_a(x, r) = B_a(x, cr)$ and $\mathcal{O}(E) = \bigcup_{g \in G} gE$.
The orbit distance is
\begin{equation}\label{eq:distances}
    d(x, y) = \min_{g \in G} \|\widehat{x} - g\widehat{y}\|.
\end{equation}
Since the identity element belongs to $G$,
\[
0 \leq d(x, y) \leq \|\widehat{x} - \widehat{y}\|\qquad \mathrm{and}
\qquad d(x, y) = d(y, x).
\]In the special case \(a=2\), \(d(x,y)\) is precisely the Dunkl distance.
For the volume estimates, we use
\begin{equation}\label{eq:V}
    V(x, y, r) = \max\bigl\{ w_{k,a}\big(B_a(x, r)\big), w_{k,a}\big(B_a(y, r)\big) \bigr\}.
\end{equation}
By coordinate derivatives we mean
\[
\partial_{\widehat{x}}^{\alpha} f(x)
= \bigl[\partial^\alpha (f \circ \phi^{-1})\bigr](\widehat{x}).
\]

\subsection{\texorpdfstring{The $(k,a)$-generalized Fourier transform}{The (k,a)-generalized Fourier transform}}
\label{sec:Fourier}

We adopt the normalization of \cite[Section~5]{BKO}, with
normalization constant given by
\begin{equation}\label{eq:Fourier-constant}
    \begin{aligned}
        c_{k,a}^{-1}
        &=\int_{\mathbb{R}^n} e^{-\|x\|^a/a}\,\mathrm{d}w_{k,a}(x)\\
        &=a^{\frac{n+2\gamma+a-2}{a}-1}
        \Gamma\left(\frac{n+2\gamma+a-2}{a}\right)
        \int_{S^{n-1}} w_{k,2}(\theta)\,\mathrm{d}\sigma(\theta).
    \end{aligned}
\end{equation}
The second equality follows by polar integration and the substitution
$v=r^a/a$. The constant is positive and finite under
\eqref{eq:parameters-minimal}.

For $\operatorname{Re} z > 0$, the normalized oscillator kernel is defined by
\begin{equation}\label{eq:oscillator-normalization}
    e^{z\Delta_{k,a}/a} f(x)
    = c_{k,a} \int_{\mathbb{R}^n} \Lambda_{k,a}(x,y;z) f(y)\,\mathrm{d}w_{k,a}(y).
\end{equation}

The $(k,a)$-generalized Fourier transform and its normalized kernel are defined by
\begin{equation}\label{s2.44}
    \mathscr{F}_{k,a} f(\xi)
    = c_{k,a} \int_{\mathbb{R}^n} B_{k,a}(\xi,x) f(x)\,\mathrm{d}w_{k,a}(x),
\end{equation}
where the normalized kernel satisfies
\begin{equation}\label{s2.10}
    B_{k,a}(x,y)
    = e^{i\pi(n+2\gamma+a-2)/(2a)}
    \Lambda_{k,a}\left(x,y;\frac{i\pi}{2}\right).
\end{equation}

The inverse transform is defined by
\begin{equation}\label{eq:Fourier-inverse}
    \mathscr{F}_{k,a}^{-1} g(x)
    = c_{k,a} \int_{\mathbb{R}^n} \overline{B_{k,a}(\xi,x)} g(\xi)\,\mathrm{d}w_{k,a}(\xi).
\end{equation}

The kernel is symmetric and satisfies the following positive scaling symmetry:
\[
B_{k,a}(x,y) = B_{k,a}(y,x)\qquad \mathrm{and}  \qquad
B_{k,a}(\lambda x,y) = B_{k,a}(x,\lambda y),\quad \lambda > 0.
\]

On a common invariant core, the spectral relation reads
\begin{equation}\label{eq:Fourier-spectral}
    \mathscr{F}_{k,a}(\mathfrak{L} f)(\xi)
    = -\|\xi\|^a \mathscr{F}_{k,a} f(\xi).
\end{equation}
The sectorial realizations used in Section~\ref{sec:heat-definition} are consistent with this relation. When the spectral kernel formula is justified, \eqref{s2.44} and \eqref{eq:Fourier-inverse} yield
\begin{equation}\label{eq:heat-spectral-background}
    h_t(x,y) = c_{k,a}^{\,2} \int_{\mathbb{R}^n} e^{-t\|\xi\|^a}
    \overline{B_{k,a}(\xi,x)} B_{k,a}(\xi,y)\,\mathrm{d}w_{k,a}(\xi).
\end{equation}
Otherwise, this identity holds after spectral truncation and is interpreted in the weak sense. The factor $c_{k,a}^2$ arises from the composition of two normalized transforms.

For $a=2$, one has $B_{k,2}(x,y) = E_k(x,-iy)$ with the normalization adopted here. Consequently, the Fourier kernel is oscillatory, whereas the heat kernel is recovered from a spectral integral involving $E_k$ at purely imaginary arguments.

\subsection{Semigroups and function spaces}\label{sec:analytic-symbols}

We write $H_t = e^{t\mathfrak{L}}$ and $P_t = e^{-t\sqrt{-\mathfrak{L}}}$, $t>0$, for the heat and Poisson semigroups associated
 with the  self-adjoint realization in Section~\ref{sec:heat-definition}. Their kernels with respect to $\mathrm{d}w_{k,a}$ are
 denoted by $h_t(x,y)$ and $p_t(x,y)$, respectively, such that
\[
H_t f(x) = \int_{\mathbb{R}^n} h_t(x,y) f(y)\,\mathrm{d}w_{k,a}(y)
\]
and
\[
P_t f(x) = \int_{\mathbb{R}^n} p_t(x,y) f(y)\,\mathrm{d}w_{k,a}(y).
\]
We write $q_t=-t\partial_t p_t$ for the Poisson cancellation kernel.

Let
\[
\mathscr D_a=\{F\circ\phi:F\in C_c^\infty(\mathbb R^n)\}.
\]
We equip $\mathscr D_a$ with the test-function topology transported
from $C_c^\infty(\mathbb R^n)$ by $F\mapsto F\circ\phi$.
Its continuous dual is denoted by $\mathscr D_a'$.
For $0<\eta\leq1$ and a continuous function $f$, write
\begin{equation}\label{eq:Holder-notation}
    [f]_\eta=\sup_{x\ne y}
    \frac{|f(x)-f(y)|}{\|\widehat x-\widehat y\|^\eta}.
\end{equation}
We write $C^\eta(\mathbb R^n)$ for the space of continuous
functions with finite seminorm $[f]_\eta$.
The seminorm $[f]_\eta$ is used in the quantitative estimates;
it does not define the test-function topology above.

The space $\operatorname{BMO}(\mathbb{R}^n, \mathrm{d}w_{k,a})$
consists of locally integrable functions with bounded mean
oscillation over the balls $B_a$. Its seminorm is
\begin{equation}\label{eq:BMO}
    \|b\|_{\operatorname{BMO}(\mathbb{R}^n, \mathrm{d}w_{k,a})}
    = \sup_{B = B_a(x,r)} \frac{1}{w_{k,a}(B)}
    \int_B |b(y) - b_B|\,\mathrm{d}w_{k,a}(y),
\end{equation}
where
\[
b_B = \frac{1}{w_{k,a}(B)} \int_B b(y)\,\mathrm{d}w_{k,a}(y).
\]
We write $\|b\|_{\operatorname{BMO}}$ for the seminorm in \eqref{eq:BMO}. Two $\operatorname{BMO}$ functions are identified if they differ by a constant.

The distributions $T1$ and $T^*1$ are tested against functions of integral zero. All identities involving these distributions follow this convention.

The distribution kernel of $T$ is denoted by $\mathcal{K}$, and its function representative away from the diagonal is denoted by $K(x,y)$. The constant $C_K$ is the universal constant appearing in the size and difference conditions of Definition~\ref{def:CZ}.

The operator pairing $\langle Tf, g\rangle$ is linear in $f$ and conjugate-linear in $g$, whereas the distributional pairing $\langle \mathcal{K}, F\rangle$ is linear in the joint test function $F$. Thus, with respect to the product measure $\mathrm{d}w_{k,a}(x)\,\mathrm{d}w_{k,a}(y)$,
\[
\langle Tf, g\rangle
= \langle \mathcal{K}, f(y)\overline{g(x)} \rangle.
\]
The identity operator is $I$, and the adjoint convention is $\langle Tf, g\rangle = \langle f, T^*g\rangle$.

Weak boundedness is denoted by $T \in \operatorname{WBP}$; the
joint-test definition and the minimal   constant
$\|T\|_{\operatorname{WBP}}$  are given in
Section~\ref{sec:constant-tests}.

For $0<s<1$, let $X_s$ be the completion of $\mathscr D_a$
with respect to the norm
\[
\|f\|_{L^2_{k,a}(\mathbb R^n)}+[f]_s.
\]

Dyadic cubes are the pullbacks under $\phi$ of cubes in the coordinate space; for such a cube $Q$, $\ell(Q)$ denotes its side length in that space.

We use the real-variable notation of \cite{Hd,HMY,35,THHLL}. For $j \in \mathbb{Z}$, $r_j = 2^{-j}$ is a coordinate scale; $S_j$ is the normalized Coifman averaging operator and $D_j = S_j - S_{j-1}$. Given a positive integer $M$, write
\[
D_j^M = \sum_{|h| \leq M} D_{j+h},\qquad U_M = \sum_{j \in
\mathbb{Z}} D_j D_j^M,\qquad \mathrm{and} \qquad  R_M = I - U_M,
\]
where $h \in \mathbb{Z}$. These formulas fix the notation; convergence and invertibility are proved in Section~\ref{sec:T1-proof}. In the construction of $S_j$, $L_j$ denotes the unnormalized cutoff average and $b_j = (L_j 1)^{-1}$ the corresponding normalizing function.

The classes $\mathcal{M}(\beta,\delta,r,x_0)$ and
$\widetilde{\mathcal{M}}(\beta,\delta,r,x_0)$ are the strong and
weak smooth-molecule spaces introduced in
Definitions~\ref{def:molecules} and~\ref{def:weak-molecules}. Here $\beta$ is the  H\"{o}lder
exponent, $\delta$ the decay exponent, $r$ the scale and $x_0$ the
center. Their norms are the minimal constants in the corresponding
size and difference conditions; the fixed root multiplicity sum
$\gamma$ is not counted as a molecule exponent.

For discrete reconstruction, $\mathcal{Q}_j$ is a sufficiently fine partition into coordinate dyadic cubes at scale $r_j$, and $x_Q$ denotes a selected sample point of $Q \in \mathcal{Q}_j$.

The notation $\widetilde{D}_j$ is reserved for dual analysis
operators  and not for adjoints of $D_j$.

The paraproduct associated with a $\operatorname{BMO}$ function $b$ is denoted by $\Pi_b$ and is defined by \eqref{eq:Pi-definition}.

We write $\mathscr{M}$ for the   Hardy--Littlewood maximal operator:
\[
\mathscr{M} f(x) = \sup_{B_a \ni x}
\frac{1}{w_{k,a}(B_a)} \int_{B_a} |f(y)|\,\mathrm{d}w_{k,a}(y).
\]
Its definition is independent of the integer $M$ appearing in $D_j^M$ and $U_M$.

\section{\texorpdfstring{$(k,a)$-heat kernel and $(k,a)$-Poisson kernel}{(k,a)-heat kernel and (k,a)-Poisson kernel}}\label{sec:kernels}
\subsection{Geometry and the heat semigroup}\label{sec:heat-definition}
\begin{proposition}\label{prop:geometry}
    The map $\phi$ in \eqref{eq:phi} is an equivariant homeomorphism and
    \[
    \phi^{-1}(z)=\left(\frac a2\right)^{2/a}\|z\|^{2/a-1}z,
    \qquad \phi(gx)=g\phi(x).
    \]
    For every nonnegative measurable function $F$,
    \begin{equation}\label{eq:change-measure}
        \int_{\mathbb R^n} F(\widehat x)\,dw_{k,a}(x)
        =\left(\frac a2\right)^{n_{k,a}-1}
        \int_{\mathbb R^n} F(z)\|z\|^{n_{k,a}-n-2\gamma}w_{k,2}(z)\,dz.
    \end{equation}
    For $F\in C^2(\mathbb R^n\setminus\{0\})$ and
    $z=\rho\theta\ne0$, the differential expression satisfies
    \begin{equation}\label{eq:Lhat}
        [\mathfrak L(F\circ\phi)](\phi^{-1}z)
        =\left(\partial_\rho^2+\frac{n_{k,a}-1}{\rho}\partial_\rho
        +\frac4{a^2\rho^2}\Delta_{k,S}\right)F(z).
    \end{equation}
    Moreover,
    \begin{equation}\label{eq:ball-volume}
        \begin{aligned}
            w_{k,a}(B_a(x,r))\asymp{}&r^n(\|\widehat x\|+r)^{n_{k,a}-n-2\gamma}\times\prod_{\alpha\in R_+}
            (|\langle\alpha,\widehat x\rangle|+r)^{2k(\alpha)},
        \end{aligned}
    \end{equation}
    and for $0<r<s$,
    \begin{equation}\label{eq:volume-growth}
        c(s/r)^{\min(n,n_{k,a})}\leq
        \frac{w_{k,a}(B_a(x,s))}{w_{k,a}(B_a(x,r))}
        \leq C(s/r)^{\max(n+2\gamma,n_{k,a})}.
    \end{equation}
\end{proposition}
\begin{proof}
    Set $x=u\theta$ and $\rho=(2/a)u^{a/2}$. Then
    \[
    \frac{d\rho}{du}=u^{a/2-1},\qquad
    u=(a\rho/2)^{2/a},
    \]
    and
    \[
    u^{n+2\gamma+a-3}\,du
    =u^{n+2\gamma+a/2-2}\,d\rho
    =\left(\frac a2\right)^{n_{k,a}-1}\rho^{n_{k,a}-1}\,d\rho.
    \]
    This proves \eqref{eq:change-measure}. Then
    \begin{align*}
        u^{2-a}\partial_u^2
        &=\partial_\rho^2+\frac{1-2/a}{\rho}\partial_\rho,\\
        (n+2\gamma-1)u^{1-a}\partial_u
        &=\frac{2(n+2\gamma-1)}{a\rho}\partial_\rho,\\
        u^{-a}\Delta_{k,S}&=\frac{4}{a^2\rho^2}\Delta_{k,S},
     \end{align*}
     and \begin{align*}    1-\frac{2}{a}+\frac{2(n+2\gamma-1)}{a}&=n_{k,a}-1.
    \end{align*}
    Thus \eqref{eq:Lhat} follows.

    If $r\leq\|\widehat x\|/2$, the radial factor in
    \eqref{eq:change-measure} is comparable to
    $\|\widehat x\|^{n_{k,a}-n-2\gamma}$.
    The standard linear-factor estimate gives
    \[
    \int_{B(z,r)}w_{k,2}(v)\,dv
    \asymp r^n\prod_{\alpha\in R_+}
    (|\langle\alpha,z\rangle|+r)^{2k(\alpha)}.
    \]
    For the lower bound, remove strips of width $\epsilon r$ about the
    relevant hyperplanes. For sufficiently small fixed $\epsilon$, the
    remaining set has volume at least $cr^n$, and each linear factor is
    bounded below by the corresponding factor on the right. Thus
    \eqref{eq:ball-volume} holds.

    If $r>\|\widehat x\|/2$, by scaling we have
    \[
    w_{k,a}(B_a(x,r))
    =\left(\frac a2\right)^{n_{k,a}-1}r^{n_{k,a}}
    \int_{B(\widehat x/r,1)}\|v\|^{n_{k,a}-n-2\gamma}w_{k,2}(v)\,dv,
    \]
    where the integrand  is locally integrable, since its radial power is
    $n_{k,a}-1>-1$.
    The last integral is a continuous positive function with the center $z$
    on $\{\|z\|\leq2\}$. Hence
    \[
    0<c\leq\int_{B(z,1)}\|v\|^{n_{k,a}-n-2\gamma}w_{k,2}(v)\,dv\leq C,
    \qquad \|z\|\leq2.
    \]
    Thus $w_{k,a}(B_a(x,r))\asymp r^{n_{k,a}}$, which is exactly
    \eqref{eq:ball-volume}.

    For $r\leq s\leq\|\widehat x\|/2$, comparison of the linear factors
    gives
    \[
    c(s/r)^n\leq
    \frac{w_{k,a}(B_a(x,s))}{w_{k,a}(B_a(x,r))}
    \leq C(s/r)^{n+2\gamma}.
    \]

    For $\|\widehat x\|/2\leq r\leq s$ the ratio is comparable to
    $(s/r)^{n_{k,a}}$.
    Splitting at $\|\widehat x\|/2$ proves \eqref{eq:volume-growth}.
\end{proof}

In particular, $\mathbb{R}^n$, equipped with the distance
$(x,y)\mapsto\|\widehat x-\widehat y\|$ and the measure
$dw_{k,a}$, is a doubling and reverse-doubling metric measure space.
The function $d$ is the associated orbit pseudometric.
Invariance of $dw_{k,a}$ under $G$ gives
\begin{equation}\label{eq:orbit-ball}
    w_{k,a}(B_a(x,r))\leq w_{k,a}(\{y:d(x,y)<r\})
    \leq|G|w_{k,a}(B_a(x,r)).
\end{equation}
If $d(x,y)\leq Ar$, then
\begin{equation}\label{eq:center-comparison}
    w_{k,a}(B_a(x,r))\asymp_A w_{k,a}(B_a(y,r)).
\end{equation}
Indeed, choose $g\in G$ with
$\|\widehat x-g\widehat y\|\leq Ar$ and use
$B_a(gy,r)\subset B_a(x,(A+1)r)$. Also, for every measurable set $E$
and $s>0$,
\begin{equation}\label{eq:homogeneity}
    w_{k,a}(\phi^{-1}(sE))=s^{n_{k,a}}w_{k,a}(\phi^{-1}(E)).
\end{equation}

The following annular estimates will be used.
\begin{lemma}\label{lem:annuli}
    For every $x\in\mathbb R^n$, $r>0$ and $\delta>0$,
    \begin{equation}\label{eq:inner-annuli}
        \int_{\{y\in\mathbb R^n:0<d(x,y)<r\}}
        \frac{d(x,y)^\delta}{w_{k,a}(B_a(x,d(x,y)))}\,dw_{k,a}(y)
        \leq Cr^\delta.
    \end{equation}
    The corresponding exterior estimate is
    \begin{equation}\label{eq:outer-annuli}
        \int_{\{y\in\mathbb R^n:d(x,y)>r\}}
        \frac{1}{w_{k,a}(B_a(x,d(x,y)))}
        \left(\frac{r}{d(x,y)}\right)^\delta\,dw_{k,a}(y)
        \leq C.
    \end{equation}
    Moreover,
    \begin{equation}\label{eq:annuli}
        \int_{\mathbb R^n}\frac{1}{V(x,y,r+d(x,y))}
        \left(\frac{r}{r+d(x,y)}\right)^\delta\,dw_{k,a}(y)
        \leq C.
    \end{equation}
    The constant is independent of $x$ and $r$. The same estimates hold
    with $x$ and $y$ interchanged.
\end{lemma}
\begin{proof}
    We first prove \eqref{eq:inner-annuli}. Decomposing the region
    $0<d(x,y)<r$ into orbit annuli gives
    \begin{align*}
        &\int_{\{y\in\mathbb R^n:0<d(x,y)<r\}}
        \frac{d(x,y)^\delta}{w_{k,a}(B_a(x,d(x,y)))}\,dw_{k,a}(y)\\
        &=\sum_{j=0}^{\infty}
        \int_{\{y\in\mathbb R^n:2^{-j-1}r\leq d(x,y)<2^{-j}r\}}
        \frac{d(x,y)^\delta}{w_{k,a}(B_a(x,d(x,y)))}\,dw_{k,a}(y)\\
        &\leq\sum_{j=0}^{\infty}
        \frac{(2^{-j}r)^\delta}
        {w_{k,a}(B_a(x,2^{-j-1}r))}
        w_{k,a}\bigl(\{y\in\mathbb R^n:d(x,y)<2^{-j}r\}\bigr).
    \end{align*}
    By \eqref{eq:orbit-ball} and the doubling property,
    \[
    w_{k,a}\bigl(\{y\in\mathbb R^n:d(x,y)<2^{-j}r\}\bigr)
    \leq |G|w_{k,a}(B_a(x,2^{-j}r))
    \leq Cw_{k,a}(B_a(x,2^{-j-1}r)).
    \]
    Consequently,
    \[
    \int_{\{y\in\mathbb R^n:0<d(x,y)<r\}}
    \frac{d(x,y)^\delta}{w_{k,a}(B_a(x,d(x,y)))}\,dw_{k,a}(y)
    \leq Cr^\delta\sum_{j=0}^{\infty}2^{-j\delta}
    =\frac{Cr^\delta}{1-2^{-\delta}}.
    \]
    This proves \eqref{eq:inner-annuli}.

    For \eqref{eq:outer-annuli}, we decompose the region $\{y\in\mathbb R^n:d(x,y)>r\}$ as
    follows:
    \begin{align*}
        &\int_{\{y\in\mathbb R^n:d(x,y)>r\}}
        \frac{1}{w_{k,a}(B_a(x,d(x,y)))}
        \left(\frac{r}{d(x,y)}\right)^\delta\,dw_{k,a}(y)\\
        &=\sum_{j=0}^{\infty}
        \int_{\{y\in\mathbb R^n:2^jr<d(x,y)\leq2^{j+1}r\}}
        \frac{1}{w_{k,a}(B_a(x,d(x,y)))}
        \left(\frac{r}{d(x,y)}\right)^\delta\,dw_{k,a}(y)\\
        &\leq\sum_{j=0}^{\infty}
        \frac{2^{-j\delta}}
        {w_{k,a}(B_a(x,2^jr))}
        w_{k,a}\bigl(\{y\in\mathbb R^n:d(x,y)\leq2^{j+1}r\}\bigr)\\
        &\leq C\sum_{j=0}^{\infty}2^{-j\delta}
        =\frac{C}{1-2^{-\delta}},
    \end{align*}
  where the last estimate follows from \eqref{eq:orbit-ball} and doubling,
    as in the preceding calculation.

    It remains to prove \eqref{eq:annuli}. On the one hand, monotonicity of
    ball measures gives
    \begin{align*}
        &\int_{\{y\in\mathbb R^n:d(x,y)<r\}}
        \frac{1}{V(x,y,r+d(x,y))}
        \left(\frac{r}{r+d(x,y)}\right)^\delta\,dw_{k,a}(y)\\
        &\qquad\leq
        \frac{w_{k,a}\bigl(\{y\in\mathbb R^n:d(x,y)<r\}\bigr)}
        {w_{k,a}(B_a(x,r))}
        \leq |G|.
    \end{align*}
    On the other hand, we obtain
    \begin{align*}
        &\int_{\{y\in\mathbb R^n:d(x,y)\geq r\}}
        \frac{1}{V(x,y,r+d(x,y))}
        \left(\frac{r}{r+d(x,y)}\right)^\delta\,dw_{k,a}(y)\\
        &=\sum_{j=0}^{\infty}
        \int_{\{y\in\mathbb R^n:2^jr\leq d(x,y)<2^{j+1}r\}}
        \frac{1}{V(x,y,r+d(x,y))}
        \left(\frac{r}{r+d(x,y)}\right)^\delta\,dw_{k,a}(y)\\
        &\leq\sum_{j=0}^{\infty}
        \frac{2^{-j\delta}}
        {w_{k,a}(B_a(x,2^jr))}
        w_{k,a}\bigl(\{y\in\mathbb R^n:d(x,y)<2^{j+1}r\}\bigr)\\
        &\leq C\sum_{j=0}^{\infty}2^{-j\delta}
        =\frac{C}{1-2^{-\delta}}.
    \end{align*}
    Adding the last two estimates proves \eqref{eq:annuli}.
Finally, interchanging the fixed and integration variables in
    the preceding calculations gives all three transposed estimates.
\end{proof}

We first verify that \eqref{eq:Dirichlet} is closable on
$L^2_{k,a}(\mathbb R^n)$.
Let $(f_m)\subset C_c^\infty(\mathbb R^n)$ be
$\mathcal E$-Cauchy and satisfy $f_m\to0$ in $L^2_{k,a}$.
On every compact subset of
\[
\Omega=\mathbb R^n\setminus
\left(\{0\}\cup\bigcup_{\alpha\in R}H_\alpha\right),
\]
the densities of $dw_{k,a}$ and $dw_{k,2}$ are bounded above
and below by positive constants. The $L^2(dw_{k,2})$ limit of
$\nabla f_m$ is therefore zero on $\Omega$, by distributional
differentiation, and hence is zero almost everywhere.
After passing to a subsequence, $f_m\to0$ almost everywhere.
Reflection invariance implies
$f_m-f_m\circ\sigma_\alpha\to0$ almost everywhere for each
$\alpha\in R_+$.
For each root with $k(\alpha)>0$, these differences are Cauchy
in the weighted $L^2$ space occurring in the jump part of
\eqref{eq:Dirichlet}, so their limits are also zero.
Consequently, $\mathcal E(f_m,f_m)\to0$.

Let $\mathcal E$ denote the closure of \eqref{eq:Dirichlet} on
$C_c^\infty(\mathbb R^n)$ in the norm
\[
\bigl(\|f\|_{L^2_{k,a}(\mathbb R^n)}^2+\mathcal E(f,f)\bigr)^{1/2}.
\]
This is the   Dunkl Dirichlet form with reference measure
$dw_{k,a}$; see \cite[Section~6.2]{FOT}. The closure is taken from
$C_c^\infty(\mathbb R^n)$, without requiring the test functions to
vanish near the origin. Every normal contraction $C$ with $C(0)=0$
and $|C(s)-C(t)|\leq|s-t|$ satisfies
\[
\mathcal E(C\circ f,C\circ f)\leq\mathcal E(f,f).
\]
Thus the form is Markovian (see \cite[Section~1.4]{FOT}).
We denote by $\mathfrak L$ its generator, which is nonpositive
and self-adjoint, and set $H_t=e^{t\mathfrak L}$.
We record that $\mathscr D_a$ is also a core for the closed form
$\mathcal E$. Let $\chi_\epsilon$ be a smooth radial cutoff in the
original variables, equal to zero on $\{\|x\|\leq\epsilon\}$ and to
one on $\{\|x\|\geq2\epsilon\}$, with
$|\nabla\chi_\epsilon|\leq C\epsilon^{-1}$, and set
\[
f_\epsilon=\chi_\epsilon f+(1-\chi_\epsilon)f(0).
\]
If $f=F\circ\phi\in\mathscr D_a$, then
$f_\epsilon\in C_c^\infty(\mathbb R^n)$. Near the origin, with
$u=\|x\|$, we have
$|f(x)-f(0)|=O(u^{a/2})$, $|\nabla f(x)|=O(u^{a/2-1})$, and
\[
|f(x)-f(\sigma_\alpha x)|
\leq C u^{a/2-1}|\langle\alpha,x\rangle|.
\]
Consequently, the expression in \eqref{eq:Dirichlet} gives
\[
\mathcal E(f-f_\epsilon,f-f_\epsilon)
=O(\epsilon^{n+2\gamma+a-2})\longrightarrow0,
\qquad \|f-f_\epsilon\|_2\longrightarrow0.
\]
Thus $\mathscr D_a\subset D(\mathcal E)$. Conversely, if
$f\in C_c^\infty(\mathbb R^n)$, then the same $f_\epsilon$ belongs
to $\mathscr D_a$, since it is constant near the origin and smooth
in the transformed coordinates elsewhere. In this case,
\[
\mathcal E(f-f_\epsilon,f-f_\epsilon)
=O(\epsilon^{n+2\gamma})\longrightarrow0,
\qquad \|f-f_\epsilon\|_2\longrightarrow0.
\]
The two cores therefore have the same closure in the form norm.
Averaging the approximants over $G$ preserves convergence in this
norm, so the same conclusion holds on the $G$-invariant subspace.

For $\zeta>0$, write
\begin{equation}\label{eq:Bessel-generator}
    B_\zeta=\partial_\rho^2+\frac{\zeta-1}{\rho}\partial_\rho.
\end{equation}
The conservative Bessel kernel with respect to $s^{\zeta-1}ds$ is
\begin{equation}\label{eq:radial-Bessel-kernel}
    b_t^\zeta(\rho,s)=\frac{(\rho s)^{1-\zeta/2}}{2t}
    e^{-(\rho^2+s^2)/(4t)}
    I_{\zeta/2-1}\left(\frac{\rho s}{2t}\right),
    \qquad \rho,s,t>0,
\end{equation}
where $I_{\zeta/2-1}$ is the modified Bessel function; see \cite{RY}.
For $\lambda\geq0$, the nonnegative series \eqref{eq:Bessel-series}
and Tonelli's theorem give
\begin{align}
    &\int_0^\infty e^{-\lambda s^2}b_t^\zeta(\rho,s)s^{\zeta-1}\,ds
    =\frac{e^{-\rho^2/(4t)}}{2t}
    \sum_{j=0}^\infty
    \frac{\rho^{2j}(4t)^{-2j-\zeta/2+1}}
    {j!\Gamma(j+\zeta/2)}
    \int_0^\infty s^{2j+\zeta-1}
    e^{-(\lambda+1/(4t))s^2}\,ds. \label{eq:Laplace-series}
\end{align}
Letting $v=(\lambda+1/(4t))s^2$, we have
\begin{align*}
    \int_0^\infty s^{2j+\zeta-1}e^{-(\lambda+1/(4t))s^2}\,ds
    &=\frac{(\lambda+1/(4t))^{-j-\zeta/2}}{2}
    \int_0^\infty v^{j+\zeta/2-1}e^{-v}\,dv\\
    &=\frac{\Gamma(j+\zeta/2)}{2}
    (\lambda+1/(4t))^{-j-\zeta/2}.
\end{align*}
Consequently,
\begin{align}
    &\int_0^\infty e^{-\lambda s^2}b_t^\zeta(\rho,s)s^{\zeta-1}\,ds
    \notag\\
    &=e^{-\rho^2/(4t)}(1+4t\lambda)^{-\zeta/2}
    \sum_{j=0}^\infty\frac1{j!}
    \left(\frac{\rho^2}{4t(1+4t\lambda)}\right)^j\notag\\
    &=(1+4t\lambda)^{-\zeta/2}
    \exp\left(-\frac{\rho^2\lambda}{1+4t\lambda}\right).
    \label{eq:Bessel-Laplace}
\end{align}
Setting $\lambda=0$ in \eqref{eq:Bessel-Laplace}, we obtain

$$  \int_0^\infty b_t^\zeta(\rho,s)s^{\zeta-1}\,ds=1.$$
    Differentiating the same formula at $\lambda=0$ gives
\begin{align}\label{eq:Bessel-mass-moment}
    \int_0^\infty s^2b_t^\zeta(\rho,s)s^{\zeta-1}\,ds
    &=-\left.\frac{d}{d\lambda}\left[
    (1+4t\lambda)^{-\zeta/2}
    e^{-\rho^2\lambda/(1+4t\lambda)}\right]\right|_{\lambda=0}\notag\\
    &=\rho^2+2\zeta t.
\end{align}

The constant angular sector of \eqref{eq:Lhat} is $B_{n_{k,a}}$. For
$Y\in\mathcal H_k^\ell$, $\ell\geq1$, direct differentiation yields
\begin{equation}\label{eq:sector-ground}
    \begin{aligned}
        &\left(B_{n_{k,a}}-\frac{4\ell(\ell+n+2\gamma-2)}{a^2\rho^2}\right)
        (\rho^{2\ell/a}v(\rho))=\rho^{2\ell/a}B_{n_{k,a}+4\ell/a}v(\rho),
    \end{aligned}
\end{equation}
because
\[
\frac{2\ell}{a}\left(\frac{2\ell}{a}+n_{k,a}-2\right)
=\frac4{a^2}\ell(\ell+n+2\gamma-2).
\]
For $v\in C_c^\infty(0,\infty)$, the corresponding form reads
\begin{equation}\label{eq:sector-form}
    \begin{aligned}
        &\int_0^\infty\left(
        |(\rho^{2\ell/a}v)'|^2+
        \frac{4\ell(\ell+n+2\gamma-2)}{a^2\rho^2}
        |\rho^{2\ell/a}v|^2\right)\rho^{n_{k,a}-1}\,d\rho\\
        &\hspace{25mm}=\int_0^\infty|v'|^2\rho^{n_{k,a}+4\ell/a-1}\,d\rho.
    \end{aligned}
\end{equation}
For $\ell\geq1$, multiplication by $\rho^{2\ell/a}$ defines
a unitary map
\[
U_\ell:L^2\bigl((0,\infty),\rho^{n_{k,a}+4\ell/a-1}\,d\rho\bigr)
\longrightarrow L^2\bigl((0,\infty),\rho^{n_{k,a}-1}\,d\rho\bigr).
\]
Equation~\eqref{eq:sector-form} identifies the forms on
$C_c^\infty(0,\infty)$. Moreover,
\[
n_{k,a}+\frac{4\ell}{a}
=2+\frac{2(n+2\gamma+2\ell-2)}a>2.
\]
Thus the origin has zero capacity for the transformed Bessel form.
Functions vanishing near zero are dense in its form domain.
For a test function in $C_c^\infty(\mathbb R^n)$, its $\ell$-th
spherical coefficient is $O(u^\ell)$ as $u\downarrow0$, by Taylor's
formula and spherical orthogonality. Therefore its transformed
radial coefficient after division by $\rho^{2\ell/a}$ is bounded
near zero. Multiplication by cutoffs vanishing near zero
approximates this coefficient in the Bessel form norm.
It follows that $U_\ell$ identifies the closed $\ell$-sector
form with the natural Bessel form.

For $n_{k,a}>2$ a radial cutoff at zero has energy $O(r^{n_{k,a}-2})$; for $n_{k,a}=2$
a logarithmic cutoff has energy $O(|\log r|^{-1})$. Thus zero has
zero capacity when $n_{k,a}\geq2$. If $n_{k,a}<2$, then $n=1$, and
\[
|v(\rho)-v(0)|^2
\leq\left(\int_0^\rho|v'(s)|^2s^{n_{k,a}-1}\,ds\right)
\frac{\rho^{2-n_{k,a}}}{2-n_{k,a}}.
\]
For \(0<n_{k,a}<2\), the origin is not negligible.
Taking the closure of the Dirichlet form from
\(C_c^\infty(\mathbb{R}^n)\), rather than from functions vanishing
near the origin, selects the conservative Bessel realization in
the even sector and the zero-trace realization in the odd sector.
The mass identity above yields \(H_t1=1\), with \(H_t\) acting on
bounded functions.

\begin{definition}\label{def:heat}
    Let $\mu_t(x,dy)$ be the transition kernel of the Markov semigroup
    $H_t$. The heat kernel $h_t(x,y)$ is the Radon--Nikodym density of
    $\mu_t(x,dy)$ with respect to $dw_{k,a}$:
    \[
    \mu_t(x,dy)=h_t(x,y)\,dw_{k,a}(y),
    \]
    equivalently,
    \[
    H_tf(x)=\int_{\R^n}h_t(x,y)f(y)\,dw_{k,a}(y).
    \]
    The existence of this density is established in
    Proposition~\ref{prop:Gaussian} below.
\end{definition}
For each fixed choice of the positive time and scale parameters,
kernel identities and estimates below are understood to hold
almost everywhere in the spatial variables, unless a particular
pointwise representative is explicitly specified.

Positivity and symmetry follow from the Dirichlet form
\eqref{eq:Dirichlet}. Applying the conservative radial semigroup to
an increasing sequence of radial cutoffs gives $H_t1=1$.
Consequently, once the density exists,
$$
    h_t(x,y)=h_t(y,x)\geq0,$$
and the mass identity is
    $$\int_{\mathbb R^n}h_t(x,y)\,dw_{k,a}(y)=1.$$
The semigroup property takes the form
    \begin{align}\label{eq:heat-semigroup}
    h_{t+s}(x,y)=\int_{\mathbb R^n}h_t(x,u)h_s(u,y)\,dw_{k,a}(u).
\end{align}
Homogeneity also gives, for $r>0$,
\begin{equation}\label{eq:heat-scaling}
    h_t(rx,ry)=r^{-n-2\gamma-a+2}h_{t/r^a}(x,y).
\end{equation}

Since the heat kernel is initially defined only almost everywhere,
we choose its value at $x=0$ to agree with the conservative radial
Bessel kernel. Letting $\rho\downarrow0$ in
\eqref{eq:radial-Bessel-kernel} and using
\eqref{eq:Fourier-constant}, we set
\begin{equation}\label{eq:origin-heat}
    h_t(0,y)=c_{k,a}(at)^{-n_{k,a}/2}
    \exp\left(-\frac{\|y\|^a}{a^2t}\right).
\end{equation}
With $u=\|y\|^a/(a^2t)$, polar integration gives
\begin{align*}
    \int_{\mathbb R^n} h_t(0,y)\,dw_{k,a}(y)
    &=c_{k,a}(at)^{-n_{k,a}/2}
    \left(\int_{S^{n-1}}w_{k,2}(\theta)\,d\sigma(\theta)\right)
    \frac{(a^2t)^{n_{k,a}/2}}a\int_0^\infty u^{n_{k,a}/2-1}e^{-u}\,du\\
    &=c_{k,a}a^{n_{k,a}/2-1}\Gamma(n_{k,a}/2)
    \int_{S^{n-1}}w_{k,2}(\theta)\,d\sigma(\theta)=1.
\end{align*}

\subsection{Heat kernel estimates}\label{sec:heat-estimates}
We first establish the orbit-Gaussian bound.
Let $\mathcal C$ be a closed Weyl chamber. On invariant functions, the
reflection part of the energy vanishes. In the coordinates
$\widehat x=\rho\theta$, the chamber energy is
\begin{equation}\label{eq:chamber-energy}
    \int_{\mathcal C}\left(|\partial_\rho f(x)|^2+
    \frac4{a^2\rho^2}|\nabla_Sf(x)|^2\right)\,dw_{k,a}(x),
    \qquad \rho=\|\widehat x\|.
\end{equation}
For fixed $a$, this form is comparable to
$\int_{\mathcal C}|\nabla_{\widehat x}f(x)|^2\,dw_{k,a}(x)$.

\begin{lemma}\label{lem:chamber-PI}
    For \(x\in\mathcal C\) and \(B=B_a(x,r)\cap\mathcal C\), with \(\widehat{x}\) and \(\widehat{y}\) as in \eqref{eq:phi},
    \begin{equation}\label{eq:chamber-vol}
        |G|^{-1}w_{k,a}(B_a(x,r))\leq w_{k,a}(B)
        \leq w_{k,a}(B_a(x,r)),
    \end{equation}
    and, for any  $f\in \mathscr D_a$,
    \begin{equation}\label{eq:PI}
        \int_B|f(y)-f_B|^2\,dw_{k,a}(y)
        \leq Cr^2\int_{B_a(x,9r)\cap\mathcal C}
        |\nabla_{\widehat y}f(y)|^2\,dw_{k,a}(y).
    \end{equation}
\end{lemma}
\begin{proof}
    Put $F=f\circ\phi^{-1}$. Folding into $\mathcal C$ decreases the distance from $x$ and has
    multiplicity at most $|G|$. This proves \eqref{eq:chamber-vol} and
    the doubling property on the chamber.

    Suppose $\|\widehat x\|>4r$. On $B(\widehat x,r)\cap\mathcal C$,
    \[
    \|z\|^{n_{k,a}-n-2\gamma}\asymp\|\widehat x\|^{n_{k,a}-n-2\gamma}\qquad \mathrm{and} \qquad
    \nabla^2\log w_{k,2}(z)
    =-2\sum_{\alpha\in R_+}k(\alpha)
    \frac{\alpha\otimes\alpha}{\langle\alpha,z\rangle^2}\leq0.
    \]
    The convex-domain weighted Poincar\'e inequality
    \cite[Theorem~1.1]{FNT}, with regularization at the walls if needed,
    gives
    \[
    \inf_c\int_{B(\widehat x,r)\cap\mathcal C}|F(z)-c|^2
    \|z\|^{n_{k,a}-n-2\gamma}w_{k,2}(z)\,dz
    \leq Cr^2\int_{B(\widehat x,r)\cap\mathcal C}
    |\nabla F(z)|^2\|z\|^{n_{k,a}-n-2\gamma}w_{k,2}(z)\,dz.
    \]
    This proves \eqref{eq:PI} in the first case.

    For $B_0=B_a(0,L)\cap\mathcal C$, reflection extension and
    \eqref{eq:spherical-eigenvalue} give, when $n\geq2$,
    \begin{equation}\label{eq:angular-PI}
        \int_{\mathcal C\cap S^{n-1}}
        |v(\theta)-v_{\mathcal C\cap S^{n-1}}|^2w_{k,2}(\theta)\,d\sigma(\theta)
        \leq\frac1{n+2\gamma-1}
        \int_{\mathcal C\cap S^{n-1}}|\nabla_Sv(\theta)|^2
        w_{k,2}(\theta)\,d\sigma(\theta).
    \end{equation}
    For the radial part, write
    \[
    v_{(0,L)}=\frac{n_{k,a}}{L^{n_{k,a}}}
    \int_0^L v(\rho)\rho^{n_{k,a}-1}\,d\rho.
    \]
    Comparison with $v(L/2)$ gives
    \begin{align*}
        \int_0^{L/2}|v(\rho)-v(L/2)|^2\rho^{n_{k,a}-1}\,d\rho
        &\leq\left(\int_0^{L/2}|v'(s)|^2s^{n_{k,a}-1}\,ds\right)
        \int_0^{L/2}\rho^{n_{k,a}-1}\int_\rho^{L/2}s^{1-n_{k,a}}\,ds\,d\rho\\
        &=\frac{L^2}{8n_{k,a}}\int_0^{L/2}|v'(s)|^2s^{n_{k,a}-1}\,ds.
    \end{align*}
    On $[L/2,L]$, the weight is comparable to $L^{n_{k,a}-1}$.
    Therefore,
    \begin{equation}\label{eq:radial-PI}
        \int_0^L|v(\rho)-v_{(0,L)}|^2\rho^{n_{k,a}-1}\,d\rho
        \leq CL^2\int_0^L|v'(\rho)|^2\rho^{n_{k,a}-1}\,d\rho.
    \end{equation}
    Set
    \[
    \bar f(\rho)=
    \frac{\displaystyle\int_{\mathcal C\cap S^{n-1}}
        f(\phi^{-1}(\rho\theta))w_{k,2}(\theta)\,d\sigma(\theta)}
    {\displaystyle\int_{\mathcal C\cap S^{n-1}}
        w_{k,2}(\theta)\,d\sigma(\theta)}.
    \]
    By \eqref{eq:angular-PI}, \eqref{eq:radial-PI}, and Jensen's inequality,
    \begin{align*}
        &\int_{B_0}|f(y)-f_{B_0}|^2\,dw_{k,a}(y)\\
        &=\int_{B_0}|f(y)-\bar f(\|\widehat y\|)|^2\,dw_{k,a}(y)
        +\int_{B_0}|\bar f(\|\widehat y\|)-f_{B_0}|^2\,dw_{k,a}(y)\\
        &\leq C L^2\int_{B_0}
        \left(\|\widehat y\|^{-2}|\nabla_S f(y)|^2
        +|\partial_\rho f(y)|^2\right)\,dw_{k,a}(y).
    \end{align*}
    For a one-dimensional reflected chamber, \eqref{eq:radial-PI}
    suffices.

    For the empty root system in dimension one,
    $n_{k,a}=2-2/a\in(0,2)$. The common trace at zero gives, on each
    half-line,
    \[
    \int_0^L|v(\rho)-v(0)|^2\rho^{n_{k,a}-1}\,d\rho
    \leq\frac{L^2}{2(2-n_{k,a})}\int_0^L|v'(\rho)|^2\rho^{n_{k,a}-1}\,d\rho.
    \]
    This yields the centered-interval inequality. Finally, for
    $\|\widehat x\|\leq4r$,
    \[
    B(\widehat x,r)\cap\mathcal C\subset B(0,5r)\cap\mathcal C
    \subset B(\widehat x,9r)\cap\mathcal C.
    \]
    Using the mean on the middle ball as a comparison constant proves
    \eqref{eq:PI}. By the core property established above, approximation in the chamber
    form norm extends \eqref{eq:PI} to the domain of the closed chamber form.
\end{proof}

\begin{proposition}\label{prop:Gaussian}
    The heat semigroup has a density satisfying
    \begin{equation}\label{eq:Gaussian}
        0\leq h_t(x,y)\leq\frac C{V(x,y,\sqrt t)}e^{-c d(x,y)^2/t}.
    \end{equation}
\end{proposition}
\begin{proof}
    The invariant form is a regular strongly local form on the closed
    chamber, with the natural reflecting walls. Its diffusion matrix in
    coordinate variables is
    \[
    P_{\mathrm{rad}}+\frac4{a^2}P_{\mathrm{tan}}\] and \[
    \min(1,4/a^2)\|\xi\|^2\leq
    \langle(P_{\mathrm{rad}}+4a^{-2}P_{\mathrm{tan}})\xi,\xi\rangle
    \leq\max(1,4/a^2)\|\xi\|^2.
    \]
  Functions that are Lipschitz in the coordinate metric on the chamber
  belong locally to the form domain. Indeed, extend them by the folding
  map and mollify in the transformed coordinates with a $G$-invariant
  radial mollifier. The core property above and the uniform bounds on
  the gradients justify this approximation in the local form norm.
  The ellipticity bounds therefore make the intrinsic metric comparable
  to the coordinate metric. It is complete and induces the topology
  of the closed chamber.
  The fixed-dilation Poincar\'e inequality \eqref{eq:PI} transfers to
  intrinsic balls by metric comparability. Together with doubling,
  it yields the local scale-invariant Sobolev inequality.
  The Gaussian upper-bound theorem for strongly local Dirichlet forms
  \cite{SturmII,SturmIII} therefore gives
    \[
    h_t^{\mathcal C}(x,y)\leq
    \frac{C e^{-c\|\widehat x-\widehat y\|^2/t}}
    {[w_{k,a}(B_a(x,\sqrt t)\cap\mathcal C)
        w_{k,a}(B_a(y,\sqrt t)\cap\mathcal C)]^{1/2}}.
    \]
    The chamber volume estimate and doubling give
    \begin{align*}
        &[w_{k,a}(B_a(x,\sqrt t)\cap\mathcal C)
        w_{k,a}(B_a(y,\sqrt t)\cap\mathcal C)]^{-1/2}\\
        &\qquad\leq\frac{C}{V(x,y,\sqrt t)}
        \left(1+\frac{\|\widehat x-\widehat y\|}{\sqrt t}\right)^{\max(n+2\gamma,n_{k,a})/2}.
    \end{align*}
    The polynomial factor is absorbed by the Gaussian.

    For each $g\in G$, reflection invariance of the measure and
    \eqref{eq:Dirichlet} give
    \begin{equation}\label{eq:form-equivariance}
        \|f\circ g\|_2=\|f\|_2\qquad \mathrm{ and} \qquad
        \mathcal E(f\circ g,f\circ g)=\mathcal E(f,f).
    \end{equation}
    The core is $G$-invariant. Polarization and closure therefore give
    \begin{equation}\label{eq:semigroup-equivariance}
        H_t(f\circ g)=(H_tf)\circ g,\qquad g\in G.
    \end{equation}
 Let $\pi:\mathbb R^n\to\mathcal C$ send each point to the
 representative of its $G$-orbit in $\mathcal C$.
 For invariant functions, the full-space energy and measure are
 $|G|$ times their chamber counterparts. Hence
 \[
 H_t(F\circ\pi)=(H_t^{\mathcal C}F)\circ\pi
 \]
 for bounded measurable $F$ on $\mathcal C$.
 For each fixed $t>0$, this gives the following identity of measures
 for almost every $x$:
 \[
 \mu_t(x,\pi^{-1}(A))
 =\int_A h_t^{\mathcal C}(\pi(x),z)\,dw_{k,a}(z),
 \qquad A\subset\mathcal C\ \text{Borel}.
 \]
 A common exceptional set can be chosen by first using a countable
 generating class and then applying the monotone class theorem.
 If $N\subset\mathbb R^n$ is a Borel null set, then
 \[
 \pi(N)=\mathcal C\cap\bigcup_{g\in G}gN
 \]
 is also a null set. Since $N\subset\pi^{-1}(\pi(N))$, the measure
 identity implies $\mu_t(x,N)=0$.
 Thus $\mu_t(x,\cdot)$ has a density $h_t(x,\cdot)$ with respect
 to $dw_{k,a}$. Decomposing the full space into chambers now yields
 \begin{equation}\label{eq:folding}
    h_t^{\mathcal C}(x^+,y^+)
    =\sum_{g\in G}h_t(x^+,gy^+)\qquad \mathrm{and}
    \qquad \|\widehat{x^+}-\widehat{y^+}\|=d(x,y),
 \end{equation}
 almost everywhere, where $x^+,y^+$ are the orbit representatives
 in $\mathcal C$.
 By \eqref{eq:semigroup-equivariance},
    \begin{equation}\label{eq:heat-equivariance}
        h_t(gx,gy)=h_t(x,y),\qquad g\in G,
    \end{equation}
    almost everywhere. Since $G$ is finite, the equalities may be imposed
    on a common full-measure set for each $t>0$. After applying a
    suitable reflection, $h_t(x,y)$ is one of the nonnegative summands in
    \eqref{eq:folding}. The chamber bound and \eqref{eq:chamber-vol}
    prove \eqref{eq:Gaussian}.
\end{proof}

The moment estimate follows from the radial and first angular sectors.
\begin{lemma}\label{lem:inverse-moment}
    If $0<\eta<\zeta$, then
    \begin{equation}\label{eq:inverse-moment}
        \int_0^\infty s^{-\eta}b_t^\zeta(\rho,s)s^{\zeta-1}\,ds
        \leq C_{\zeta,\eta}\rho^{-\eta},\qquad \rho,t>0.
    \end{equation}
\end{lemma}
\begin{proof}
    Using the identity $s^{-\eta}=\Gamma(\eta/2)^{-1}
    \int_0^\infty\lambda^{\eta/2-1}e^{-\lambda s^2}\,d\lambda$,
    together with Tonelli's theorem and \eqref{eq:Bessel-Laplace}, we
    obtain
    \begin{align}
        &\int_0^\infty s^{-\eta}b_t^\zeta(\rho,s)s^{\zeta-1}\,ds\notag\\
        &=\frac1{\Gamma(\eta/2)}\int_0^\infty
        \lambda^{\eta/2-1}(1+4t\lambda)^{-\zeta/2}
        e^{-\rho^2\lambda/(1+4t\lambda)}\,d\lambda\notag\\
        &=\frac1{\Gamma(\eta/2)}\int_0^{1/(4t)}
        v^{\eta/2-1}(1-4tv)^{(\zeta-\eta)/2-1}e^{-\rho^2v}\,dv,
        \label{eq:inverse-integral}
    \end{align}
where we have substituted   $v=\lambda/(1+4t\lambda)$ in the last
line. For $0<v<1/(8t)$, the factor $(1-4tv)^{(\zeta-\eta)/2-1}$ is
bounded, and hence
    \[
    \int_0^{1/(8t)}v^{\eta/2-1}
    (1-4tv)^{(\zeta-\eta)/2-1}e^{-\rho^2v}\,dv
    \leq C\int_0^\infty v^{\eta/2-1}e^{-\rho^2v}\,dv
    =C\rho^{-\eta}.
    \]
    On the remaining interval,
    \begin{align*}
        &\int_{1/(8t)}^{1/(4t)}v^{\eta/2-1}
        (1-4tv)^{(\zeta-\eta)/2-1}e^{-\rho^2v}\,dv\\
        &\quad\leq Ct^{-\eta/2}e^{-\rho^2/(8t)}
        \int_{1/2}^1u^{\eta/2-1}(1-u)^{(\zeta-\eta)/2-1}\,du
        \leq C\rho^{-\eta}.
    \end{align*}
    The Beta integral is finite because $\zeta>\eta>0$.
\end{proof}

\begin{lemma}\label{lem:displacement}
 There exists a radial scalar function \(F_t\) such that
 \begin{equation}\label{eq:first-sector}
H_t(\widehat y_j)(x)=\widehat x_jF_t(\|\widehat x\|)
    \qquad\text{and}\qquad
    H_t(\|\widehat y\|^2)(x)=\|\widehat x\|^2+2n_{k,a}t,
 \end{equation}
 where \(\widehat{x}\) and \(\widehat{y}\) are as in \eqref{eq:phi}.
Under the standing assumptions \eqref{eq:parameters-minimal},
the estimate \eqref{eq:moment} holds with a constant independent of
\(x\in\mathbb R^n\) and \(t>0\).
\end{lemma}
\begin{proof}
    The radial identity follows from \eqref{eq:Bessel-mass-moment}.
The Cauchy–Schwarz inequality gives
    \[
    \int_{\mathbb R^n}|\widehat y_j|h_t(x,y)\,dw_{k,a}(y)
    \leq\bigl(\|\widehat x\|^2+2n_{k,a}t\bigr)^{1/2}.
    \]
    Thus radial truncation justifies the first-sector identity.

    Write $\rho=\|\widehat x\|>0$. Equation~\eqref{eq:sector-ground}, with
    $\ell=1$ gives
    \begin{equation}\label{eq:F-ground}
        F_t(\rho)=\rho^{2/a-1}\int_0^\infty
        b_t^{n_{k,a}+4/a}(\rho,s)s^{1-2/a}s^{n_{k,a}+4/a-1}\,ds.
    \end{equation}
    Put
    \[
    \kappa=\left(\frac2a-1\right)
    \left(n_{k,a}+\frac2a-1\right) \qquad \mathrm{and} \qquad
    n_{k,a}+\frac2a-1=1+\frac2a(n+2\gamma-1)>0.
    \]
    Then
    \begin{equation}\label{eq:first-sector-power}
        B_{n_{k,a}+4/a}s^{1-2/a}=-\kappa s^{-1-2/a}.
    \end{equation}
    For fixed $u,\rho>0$, the Bessel series gives
    \[
    b_u^{n_{k,a}+4/a}(\rho,s)=O(1) \qquad \mathrm{and} \qquad
    \partial_sb_u^{n_{k,a}+4/a}(\rho,s)=O(s)\quad \mathrm{as} \quad
    s\downarrow 0.
    \]
    Green's formula on $[\epsilon,R]$ implies
    \begin{align*}
        &\int_\epsilon^R
        \left\{(B_{n_{k,a}+4/a}b_u^{n_{k,a}+4/a}(\rho,\cdot))(s)
        s^{1-2/a}
        -b_u^{n_{k,a}+4/a}(\rho,s)B_{n_{k,a}+4/a}(s^{1-2/a})\right\}
        s^{n_{k,a}+4/a-1}\,ds\\
        &=\left[s^{n_{k,a}+4/a-1}\left\{
        (\partial_sb_u^{n_{k,a}+4/a}(\rho,s))s^{1-2/a}
        -b_u^{n_{k,a}+4/a}(\rho,s)(1-2/a)s^{-2/a}\right\}\right]_{s=\epsilon}^{s=R}.
    \end{align*}
As $s\downarrow0$, the expression in square brackets is
$O(s^{n_{k,a}+2/a+1})+O(s^{n_{k,a}+2/a-1})$ and hence tends to zero,
since $n_{k,a}+2/a-1>0$. As $s\to\infty$, the Bessel kernel and its
derivative have Gaussian decay in $s-\rho$, up to polynomial factors,
uniformly for $u$ in any compact subinterval of $(0,\infty)$.
Thus both boundary terms vanish. The heat equation and integration
    over
    $u\in[\tau,t]$, $0<\tau<t$, give
    \begin{equation}\label{eq:F-positive-time}
        F_t(\rho)-F_\tau(\rho)
        =-\kappa\rho^{2/a-1}\int_\tau^t\int_0^\infty
        b_u^{n_{k,a}+4/a}(\rho,s)s^{-1-2/a}s^{n_{k,a}+4/a-1}\,ds\,du.
    \end{equation}
    Lemma~\ref{lem:inverse-moment}, applied with
    \[
    \zeta=n_{k,a}+4/a,\qquad \eta=1+2/a,\qquad \mathrm{and} \qquad
    \zeta-\eta=n_{k,a}+2/a-1>0,
    \]
    implies
    \begin{equation}\label{eq:F-integrable-derivative}
        \int_0^\infty b_u^{n_{k,a}+4/a}(\rho,s)
        s^{-1-2/a}s^{n_{k,a}+4/a-1}\,ds
        \leq C\rho^{-1-2/a},\qquad u>0.
    \end{equation}

    Fix $t_0>0$. For $0<\tau\leq t_0$ and $0<\epsilon<\rho<R$,
    \begin{align}
        &\int_0^\epsilon b_\tau^{n_{k,a}+4/a}(\rho,s)
        s^{1-2/a}s^{n_{k,a}+4/a-1}\,ds\notag\\
        &\qquad\leq\epsilon^2\int_0^\infty b_\tau^{n_{k,a}+4/a}(\rho,s)
        s^{-1-2/a}s^{n_{k,a}+4/a-1}\,ds
        \leq C\epsilon^2\rho^{-1-2/a},\label{eq:F-small-tail}\\
        &\int_R^\infty b_\tau^{n_{k,a}+4/a}(\rho,s)
        s^{1-2/a}s^{n_{k,a}+4/a-1}\,ds\notag\\
        &\qquad\leq R^{-1-2/a}\int_0^\infty s^2
        b_\tau^{n_{k,a}+4/a}(\rho,s)s^{n_{k,a}+4/a-1}\,ds\notag\\
        &\qquad\leq R^{-1-2/a}\bigl(\rho^2+2(n_{k,a}+4/a)t_0\bigr).
        \label{eq:F-large-tail}
    \end{align}
    Here we used \eqref{eq:F-integrable-derivative} and
    \eqref{eq:Bessel-mass-moment}. For every $\zeta>0$,
    differentiating \eqref{eq:Bessel-Laplace} twice gives
    \begin{equation}\label{eq:Bessel-concentration}
        \int_0^\infty(s^2-\rho^2)^2
        b_\tau^\zeta(\rho,s)s^{\zeta-1}\,ds
        =8\rho^2\tau+4\zeta(\zeta+2)\tau^2\longrightarrow 0\ \ \mathrm{as} \ \ \tau\rightarrow 0.
    \end{equation}
    The Bessel transition measures therefore converge weakly to the point mass at $\rho$.
   Applying this to the truncated functions  and using
    \eqref{eq:F-small-tail}--\eqref{eq:F-large-tail}, we obtain
    \[
    \lim_{\tau\downarrow0}\int_0^\infty
    b_\tau^{n_{k,a}+4/a}(\rho,s)s^{1-2/a}s^{n_{k,a}+4/a-1}\,ds
    =\rho^{1-2/a} \qquad \mathrm{and} \ \ \
    \lim_{\tau\downarrow0}F_\tau(\rho)=1.
    \]
    Letting $\tau\to0$ in \eqref{eq:F-positive-time} and using
    \eqref{eq:F-integrable-derivative}  give
    \begin{equation}\label{eq:F-difference}
        1-F_t(\rho)=\kappa\rho^{2/a-1}
        \int_0^t\int_0^\infty
        b_u^{n_{k,a}+4/a}(\rho,s)s^{-1-2/a}s^{n_{k,a}+4/a-1}\,ds\,du.
    \end{equation}

    The inner integral is nonnegative. If $\kappa\le 0, $ then $1-F_t\le 0$ and the claim is trivial;
 if $\kappa>0$,  \eqref{eq:F-integrable-derivative}  gives
    \begin{equation}\label{eq:F-estimate}
        \rho^2(1-F_t(\rho))\leq C\max(\kappa,0)t.
    \end{equation}
    Consequently,
    \begin{align}
        \int_{\mathbb R^n}\|\widehat x-\widehat y\|^2h_t(x,y)\,dw_{k,a}(y)
        &=\|\widehat x\|^2+H_t(\|\widehat y\|^2)(x)
        -2\sum_{j=1}^n\widehat x_jH_t(\widehat y_j)(x)\notag\\
        &=2n_{k,a}t+2\rho^2(1-F_t(\rho))\notag\\
        &\leq [2n_{k,a}+C\max(\kappa,0)]t.
        \label{eq:displacement-identity}
    \end{align}
    At $\rho=0$ the radial identity gives $2n_{k,a}t$. This proves the lemma.
\end{proof}

\begin{proof}[Proof of Theorem~\ref{thm:heat}]
    It remains to prove \eqref{eq:heat-improved}. By \eqref{eq:moment},
    for every $r>0$,
    \begin{align}\label{eq:moment-tail}
        &\int_{\{u\in\mathbb R^n:\|\widehat x-\widehat u\|\geq r\}}
        h_t(x,u)\,dw_{k,a}(u)\notag\\
        &\qquad\leq\frac{1}{r^2}\int_{\mathbb R^n}
        \|\widehat x-\widehat u\|^2h_t(x,u)\,dw_{k,a}(u)
        \leq\frac{Ct}{r^2}.
    \end{align}
    Suppose first that $\|\widehat x-\widehat y\|\geq\sqrt t$. The triangle inequality gives
    \[
    \mathbb R^n=
    \left\{u:\|\widehat x-\widehat u\|\geq\frac{\|\widehat x-\widehat y\|}{2}\right\}
    \cup
    \left\{u:\|\widehat u-\widehat y\|\geq\frac{\|\widehat x-\widehat y\|}{2}\right\}.
    \]
    Thus the semigroup identity and positivity imply
    \begin{align*}
        h_{2t}(x,y)
        &\leq\int_{\|\widehat x-\widehat u\|\geq\|\widehat x-\widehat y\|/2}
        h_t(x,u)h_t(u,y)\,dw_{k,a}(u)+\int_{\|\widehat u-\widehat y\|\geq\|\widehat x-\widehat y\|/2}
        h_t(x,u)h_t(u,y)\,dw_{k,a}(u).
    \end{align*}
    For the first integral, \eqref{eq:Gaussian} and \eqref{eq:moment-tail}
    yield
    \begin{align*}
        &\int_{\|\widehat x-\widehat u\|\geq\|\widehat x-\widehat y\|/2}
        h_t(x,u)h_t(u,y)\,dw_{k,a}(u)\\
        &\qquad\leq\frac{C}{w_{k,a}(B_a(y,\sqrt t))}
        \int_{\|\widehat x-\widehat u\|\geq\|\widehat x-\widehat y\|/2}
        h_t(x,u)\,dw_{k,a}(u)\\
        &\qquad\leq\frac{Ct}
        {w_{k,a}(B_a(y,\sqrt t))\,\|\widehat x-\widehat y\|^2}.
    \end{align*}
    For the second integral, symmetry gives
    \begin{align*}
        &\int_{\|\widehat u-\widehat y\|\geq\|\widehat x-\widehat y\|/2}
        h_t(x,u)h_t(u,y)\,dw_{k,a}(u)\\
        &\qquad\leq\frac{C}{w_{k,a}(B_a(x,\sqrt t))}
        \int_{\|\widehat y-\widehat u\|\geq\|\widehat x-\widehat y\|/2}
        h_t(y,u)\,dw_{k,a}(u)\\
        &\qquad\leq\frac{Ct}
        {w_{k,a}(B_a(x,\sqrt t))\,\|\widehat x-\widehat y\|^2}.
    \end{align*}
    Combining these two bounds, we obtain
    \begin{equation}\label{eq:pre-interpolation}
        h_{2t}(x,y)\leq\frac{Ct}{\|\widehat x-\widehat y\|^2}
        \left\{\frac{1}{w_{k,a}(B_a(x,\sqrt t))}
        +\frac{1}{w_{k,a}(B_a(y,\sqrt t))}\right\}.
    \end{equation}
    Choose $g\in G$ with $\|\widehat x-g\widehat y\|=d(x,y)$.
    The inclusion
    \[
    B_a(gy,\sqrt t)\subset B_a(x,\sqrt t+d(x,y))
    \]
    and \eqref{eq:volume-growth} give
    \[
    w_{k,a}(B_a(y,\sqrt t))
    \leq C\left(1+\frac{d(x,y)}{\sqrt t}\right)^{\max(n+2\gamma,n_{k,a})}
    w_{k,a}(B_a(x,\sqrt t)).
    \]
    Interchanging $x$ and $y$ gives the reverse comparison. Consequently,
    \begin{align}\label{eq:volume-loss}
        &\frac{1}{w_{k,a}(B_a(x,\sqrt t))}
        +\frac{1}{w_{k,a}(B_a(y,\sqrt t))}\leq\frac{C}{V(x,y,\sqrt t)}
        \left(1+\frac{d(x,y)}{\sqrt t}\right)^{\max(n+2\gamma,n_{k,a})}.
    \end{align}
    For $0<\sigma<2$, combine \eqref{eq:pre-interpolation} to the power
    $\sigma/2$ with \eqref{eq:Gaussian} at time $2t$ to the power
    $1-\sigma/2$. By doubling and \eqref{eq:volume-loss},
    \begin{align}\label{eq:interpolation}
        h_{2t}(x,y)
        &\leq\frac{C}{V(x,y,\sqrt t)}
        \left(\frac{\sqrt t}{\|\widehat x-\widehat y\|}\right)^\sigma
        \left(1+\frac{d(x,y)}{\sqrt t}\right)^{\sigma\max(n+2\gamma,n_{k,a})/2}
        \notag\\
        &\qquad\times
        \exp\left(-c(1-\sigma/2)\frac{d(x,y)^2}{t}\right)\notag\\
        &\leq\frac{C_\sigma}{V(x,y,\sqrt t)}
        \left(\frac{\sqrt t}{\|\widehat x-\widehat y\|}\right)^\sigma
        \exp\left(-c_\sigma\frac{d(x,y)^2}{t}\right).
    \end{align}
    The last inequality absorbs the polynomial factor into the Gaussian.
    For $\|\widehat x-\widehat y\|<\sqrt t$, the factor
    $(1+\|\widehat x-\widehat y\|/\sqrt t)^{-\sigma}$ is bounded below by
    $2^{-\sigma}$; hence the factor $(1+\|\widehat x-\widehat y\|/\sqrt t)^{-\sigma}$
can be inserted into \eqref{eq:Gaussian} at the cost of a constant.
Replacing $2t$ by $t$ and
    using doubling proves \eqref{eq:heat-improved}.
\end{proof}

\begin{remark}
 For the endpoint estimates in the case $a=2$, see \cite{DH}.
\end{remark}

\subsection{Poisson kernel estimates}\label{sec:poisson-estimates}
\begin{definition}\label{def:poisson}
    The Poisson semigroup is $P_t=e^{-t\sqrt{-\mathfrak L}}$. Its kernel
    with respect to $dw_{k,a}$ is given by
    \begin{equation}\label{eq:subordination}
        p_t(x,y)=\frac{t}{2\sqrt\pi}\int_0^\infty
        u^{-3/2}e^{-t^2/(4u)}h_u(x,y)\,du.
    \end{equation}
\end{definition}
The spectral theorem gives \eqref{eq:subordination}. By Tonelli's
theorem and \eqref{eq:heat-semigroup}, we obtain
\[
\begin{gathered}
    p_t(x,y)=p_t(y,x)\geq0,\qquad
    \int_{\mathbb R^n}p_t(x,y)\,dw_{k,a}(y)=1,\end{gathered}\]
 and  \[\begin{gathered}  p_{t+s}(x,y)=\int_{\mathbb R^n}p_t(x,u)p_s(u,y)\,dw_{k,a}(u).
\end{gathered}
\]

\begin{proof}[Proof of Theorem~\ref{thm:poisson}]
    Since $d(x,y)\leq\|\widehat x-\widehat y\|$,
    \[
    t^2+d(x,y)^2\asymp(t+d(x,y))^2.
    \]
    Theorem~\ref{thm:heat} and \eqref{eq:subordination} give
    \begin{align}\label{eq:Poisson-start}
        p_t(x,y)
        &\leq Ct\int_0^{(t+d(x,y))^2}
        \frac{u^{-3/2}}{V(x,y,\sqrt u)}
        \left(1+\frac{\|\widehat x-\widehat y\|}{\sqrt u}\right)^{-\sigma}
        e^{-c(t+d(x,y))^2/u}\,du\notag\\
        &\quad+Ct\int_{(t+d(x,y))^2}^{\infty}
        \frac{u^{-3/2}}{V(x,y,\sqrt u)}
        \left(1+\frac{\|\widehat x-\widehat y\|}{\sqrt u}\right)^{-\sigma}
        e^{-c(t+d(x,y))^2/u}\,du\notag\\
        &=:I_1+I_2.
    \end{align}
    For $0<u\leq(t+d(x,y))^2$, the upper volume growth gives
    \[
    \frac{1}{V(x,y,\sqrt u)}
    \leq\frac{C}{V(x,y,t+d(x,y))}
    \left(\frac{t+d(x,y)}{\sqrt u}\right)^{\max(n+2\gamma,n_{k,a})}.
    \]
    In the same range,
    \[
    \left(1+\frac{\|\widehat x-\widehat y\|}{\sqrt u}\right)^{-\sigma}
    \leq\left(\frac{t+d(x,y)}{t+d(x,y)+\|\widehat x-\widehat y\|}\right)^\sigma
    \leq\left(\frac{t+d(x,y)}{t+\|\widehat x-\widehat y\|}\right)^\sigma.
    \]
    Applying these two estimates and setting $v=(t+d(x,y))^2/u$, we obtain
    \begin{align}\label{eq:Poisson-small}
        I_1
        &\leq\frac{Ct}{V(x,y,t+d(x,y))}
        \left(\frac{t+d(x,y)}{t+\|\widehat x-\widehat y\|}\right)^\sigma\notag\\
        &\quad\times\int_0^{(t+d(x,y))^2}
        u^{-3/2}\left(\frac{t+d(x,y)}{\sqrt u}\right)^{\max(n+2\gamma,n_{k,a})}
        e^{-c(t+d(x,y))^2/u}\,du\notag\\
        &=\frac{C}{V(x,y,t+d(x,y))}\frac{t}{t+d(x,y)}
        \left(\frac{t+d(x,y)}{t+\|\widehat x-\widehat y\|}\right)^\sigma
        \int_1^\infty v^{(\max(n+2\gamma,n_{k,a})-1)/2}e^{-cv}\,dv\notag\\
        &\leq\frac{C}{V(x,y,t+d(x,y))}\frac{t}{t+d(x,y)}
        \left(\frac{t+d(x,y)}{t+\|\widehat x-\widehat y\|}\right)^\sigma.
    \end{align}
    For $u\geq(t+d(x,y))^2$, the lower volume growth implies
    \[
    \frac{1}{V(x,y,\sqrt u)}
    \leq\frac{C}{V(x,y,t+d(x,y))}
    \left(\frac{t+d(x,y)}{\sqrt u}\right)^{\min(n,n_{k,a})}.
    \]
    Furthermore, $\sqrt u\geq t+d(x,y)\geq t$ yields
    \[
    \left(1+\frac{\|\widehat x-\widehat y\|}{\sqrt u}\right)^{-\sigma}
    =\frac{u^{\sigma/2}}{(\sqrt u+\|\widehat x-\widehat y\|)^\sigma}
    \leq\frac{u^{\sigma/2}}{(t+\|\widehat x-\widehat y\|)^\sigma}.
    \]
    We therefore have
    \begin{align}\label{eq:Poisson-large}
        I_2
        &\leq\frac{Ct(t+d(x,y))^{\min(n,n_{k,a})}}
        {V(x,y,t+d(x,y))(t+\|\widehat x-\widehat y\|)^\sigma}
        \int_{(t+d(x,y))^2}^{\infty}
        u^{(\sigma-\min(n,n_{k,a})-3)/2}\,du\notag\\
        &=\frac{2C}{1+\min(n,n_{k,a})-\sigma}
        \frac{1}{V(x,y,t+d(x,y))}\frac{t}{t+d(x,y)}
        \left(\frac{t+d(x,y)}{t+\|\widehat x-\widehat y\|}\right)^\sigma.
    \end{align}
    The integral converges because $\sigma<1+\min(n,n_{k,a})$.
    Combining \eqref{eq:Poisson-small} with \eqref{eq:Poisson-large}
    proves \eqref{eq:poisson-improved} when $m=0$.

    To estimate the time derivatives, put $v=t^2/(4u)$. A straightforward induction gives
    \begin{equation}\label{eq:sub-derivative-polynomial}
        \partial_t^m(te^{-v})=t^{1-m}P_m(v)e^{-v}
        \qquad \mathrm{and} \qquad P_0(v)=1.
    \end{equation}
    The polynomials satisfy
    \[
    P_{m+1}(v)=(1-m-2v)P_m(v)+2vP_m'(v).
    \]
    Since $|P_m(v)|e^{-v/2}\leq C_m$ for $v\geq0$,
    \begin{equation}\label{eq:sub-derivative}
        |\partial_t^m(te^{-t^2/(4u)})|
        \leq C_mt^{1-m}e^{-t^2/(8u)}.
    \end{equation}
    Differentiation under the subordination integral gives
    \[
    \partial_t^mp_t(x,y)=\frac{1}{2\sqrt\pi}\int_0^\infty
    u^{-3/2}\partial_t^m(te^{-t^2/(4u)})h_u(x,y)\,du.
    \]
    For $t$ in a compact subinterval of $(0,\infty)$,
    \eqref{eq:sub-derivative} supplies an integrable majorant. Hence
    \begin{align}\label{eq:Poisson-derivative-domination}
        |t^m\partial_t^mp_t(x,y)|
        &\leq C_mt\int_0^\infty
        u^{-3/2}e^{-t^2/(8u)}h_u(x,y)\,du\notag\\
        &\leq C_mp_{t/\sqrt2}(x,y).
    \end{align}
    Applying the size estimate at $t/\sqrt2$ and using the doubling condition
    prove
    \eqref{eq:poisson-improved}. The choice $\sigma=1$ is admissible
    and gives \eqref{eq:poisson-one}.

    For $q_t=-t\partial_tp_t$, \eqref{eq:Poisson-derivative-domination}
    and the mass identity give
    \begin{align}\label{eq:q-Schur}
        \operatorname*{ess\,sup}_{x\in\mathbb R^n}
        \int_{\mathbb R^n}|q_t(x,y)|\,dw_{k,a}(y)
        &\leq C,\notag\\
        \intertext{and, by symmetry,}
        \operatorname*{ess\,sup}_{y\in\mathbb R^n}
        \int_{\mathbb R^n}|q_t(x,y)|\,dw_{k,a}(x)
        &\leq C.
    \end{align}
 To justify differentiation of the mass identity, fix $t_0>0$.
 For $t\in[t_0/2,2t_0]$, \eqref{eq:sub-derivative} with $m=1$
 and the subordination formula give
 \[
 |\partial_t p_t(x,y)|
 \leq \frac{C}{t_0}p_{t_0/(2\sqrt2)}(x,y).
 \]
 The right-hand side is independent of $t$ and has integral
 $C/t_0$ in $y$. Differentiation under the integral sign is
 therefore valid:
    \begin{align}\label{eq:q-cancel}
        \int_{\mathbb R^n}q_t(x,y)\,dw_{k,a}(y)
        &=-t\partial_t\int_{\mathbb R^n}p_t(x,y)\,dw_{k,a}(y)
        =-t\partial_t1=0.\notag\\
        \intertext{Interchanging $x$ and $y$, we also obtain}
        \int_{\mathbb R^n}q_t(x,y)\,dw_{k,a}(x)
        &=-t\partial_t\int_{\mathbb R^n}p_t(x,y)\,dw_{k,a}(x)=0.
    \end{align}
    The two Schur bounds give uniform boundedness on
    $L^1_{k,a}(\mathbb R^n)$ and $L^\infty_{k,a}(\mathbb R^n)$.
    Interpolation then yields uniform boundedness on  $L^p_{k,a}(\mathbb R^n)$  for all $1<p<\infty$.
\end{proof}
\section{\texorpdfstring{$(k,a)$-Calder\'{o}n-Zygmund singular integral operators}{(k,a)-Calderon-Zygmund singular integral operators}}\label{sec:CZ}
We establish the kernel estimates in Definition~\ref{def:CZ} and
prove Theorem~\ref{thm:Lp}.

\subsection{Kernel estimates and singular integral operators}
The following lemma gives the required scale-integral estimates.

\begin{lemma}\label{lem:kernel-integration}
Here and below, for every \(z\in\mathbb R^n\),
we write \(\widehat z=\phi(z)\), as in \eqref{eq:phi}.
Let $\{\Psi_t(x,y)\}_{t>0}$ be a measurable family satisfying
    \begin{enumerate}[label=\textup{(\roman*)},ref=\thetheorem(\roman*),leftmargin=2.2em,labelsep=0.55em,itemsep=7pt,topsep=6pt,parsep=0pt]
        \item\label{eq:scale-size}
        $\displaystyle
        |\Psi_t(x,y)|\leq
        \frac{Ct}{V(x,y,t+d(x,y))(t+\|\widehat x-\widehat y\|)};
        $

        \item\label{eq:scale-y}
        $\displaystyle
        \begin{aligned}[t]
            |\Psi_t(x,y)-\Psi_t(x,y')|
            &\leq C\frac{\|\widehat y-\widehat{y'}\|}{t}
            \biggl\{\frac{t}{V(x,y,t+d(x,y))(t+\|\widehat x-\widehat y\|)}\\
            &\hspace{13mm}{}+
            \frac{t}{V(x,y',t+d(x,y'))(t+\|\widehat x-\widehat{y'}\|)}\biggr\},
        \end{aligned}
        $
        \par\smallskip
        whenever $\|\widehat y-\widehat{y'}\|\leq t$;

        \item\label{eq:scale-x}
        $\displaystyle
        \begin{aligned}[t]
            |\Psi_t(x,y)-\Psi_t(x',y)|
            &\leq C\frac{\|\widehat x-\widehat{x'}\|}{t}
            \biggl\{\frac{t}{V(x,y,t+d(x,y))(t+\|\widehat x-\widehat y\|)}\\
            &\hspace{13mm}{}+
            \frac{t}{V(x',y,t+d(x',y))(t+\|\widehat{x'}-\widehat y\|)}\biggr\},
        \end{aligned}
        $
        \par\smallskip
        whenever $\|\widehat x-\widehat{x'}\|\leq t$.
    \end{enumerate}
    Then the integral
    \begin{equation}\label{eq:integrated-kernel}
        K(x,y)=\int_0^\infty \Psi_t(x,y)\,\frac{dt}{t},\qquad d(x,y)>0,
    \end{equation}
    is absolutely convergent. For $0<\varepsilon<\min\{1,n,n_{k,a}\}$,
    its size satisfies
    \begin{equation}\label{eq:integrated-size}
        |K(x,y)|\leq\frac{C}{V(x,y,d(x,y))}
        \left(\frac{d(x,y)}{\|\widehat x-\widehat y\|}\right)^\varepsilon,
    \end{equation}
    For the second variable,
    \begin{equation}\label{eq:integrated-y}
        |K(x,y)-K(x,y')|\leq\frac{C}{V(x,y,d(x,y))}
        \left(\frac{\|\widehat y-\widehat{y'}\|}
        {\|\widehat x-\widehat y\|}\right)^\varepsilon
    \end{equation}
    for $\|\widehat y-\widehat{y'}\|\leq d(x,y)/2$, and
    \begin{equation}\label{eq:integrated-x}
        |K(x,y)-K(x',y)|\leq\frac{C}{V(x,y,d(x,y))}
        \left(\frac{\|\widehat x-\widehat{x'}\|}
        {\|\widehat x-\widehat y\|}\right)^\varepsilon
    \end{equation}
    for $\|\widehat x-\widehat{x'}\|\leq d(x,y)/2$.

\end{lemma}
\begin{proof}
    Fix $x,y\in\mathbb R^n$ with $0<d(x,y)\leq\|\widehat x-\widehat y\|$.
    By \eqref{eq:center-comparison},
    \begin{equation}\label{eq:vol-at-d}
        w_{k,a}(B_a(x,d(x,y)))\asymp w_{k,a}(B_a(y,d(x,y)))\asymp V(x,y,d(x,y)).
    \end{equation}
    Choose $\varepsilon<\delta<\min\{1,n,n_{k,a}\}$.
    For $0<t\leq d(x,y)$, monotonicity gives
    \[
    \frac{1}{V(x,y,t+d(x,y))}\leq\frac{1}{V(x,y,d(x,y))}.
    \]
    For $t\geq d(x,y)$, the lower bound in \eqref{eq:volume-growth} yields
    \[
    V(x,y,t+d(x,y))\geq V(x,y,t)
    \geq C^{-1}\left(\frac{t}{d(x,y)}\right)^\delta V(x,y,d(x,y)).
    \]
    We first prove the size estimate. By condition \textup{(i)},
    \begin{align}
        \int_0^\infty|\Psi_t(x,y)|\,\frac{dt}{t}
        &\leq C\int_0^{d(x,y)}
        \frac{dt}{V(x,y,t+d(x,y))(t+\|\widehat x-\widehat y\|)}\notag\\
        &\quad+C\int_{d(x,y)}^{\|\widehat x-\widehat y\|}
        \frac{dt}{V(x,y,t+d(x,y))(t+\|\widehat x-\widehat y\|)}\notag\\
        &\quad+C\int_{\|\widehat x-\widehat y\|}^{\infty}
        \frac{dt}{V(x,y,t+d(x,y))(t+\|\widehat x-\widehat y\|)}\notag\\
        &=:C(I_1+I_2+I_3).
        \label{eq:kernel-three-ranges}
    \end{align}
    For $I_1$, we use $t+\|\widehat x-\widehat y\|\geq\|\widehat x-\widehat y\|$ to obtain
    \begin{align}
        I_1
        &\leq\frac{1}{\|\widehat x-\widehat y\|\,V(x,y,d(x,y))}\int_0^{d(x,y)}1\,dt\notag\\
        &=\frac{1}{V(x,y,d(x,y))}\frac{d(x,y)}{\|\widehat x-\widehat y\|}.
        \label{eq:kernel-I1}
    \end{align}
    For $I_2$, the preceding volume comparison gives
    \begin{align}
        I_2
        &\leq\frac{Cd(x,y)^\delta}{\|\widehat x-\widehat y\|\,V(x,y,d(x,y))}
        \int_{d(x,y)}^{\|\widehat x-\widehat y\|}t^{-\delta}\,dt\notag\\
        &=\frac{C}{(1-\delta)V(x,y,d(x,y))}
        \left\{\left(\frac{d(x,y)}{\|\widehat x-\widehat y\|}\right)^\delta
        -\frac{d(x,y)}{\|\widehat x-\widehat y\|}\right\}\notag\\
        &\leq\frac{C}{(1-\delta)V(x,y,d(x,y))}
        \left(\frac{d(x,y)}{\|\widehat x-\widehat
        y\|}\right)^\delta,
        \label{eq:kernel-I2}
    \end{align}
where $\delta<1$. If $d(x,y)=\|\widehat x-\widehat y\|$, then $I_2=0$.
    For the last integral, $t+\|\widehat x-\widehat y\|\geq t$ and $\delta>0$ imply
    \begin{align}
        I_3
        &\leq\frac{Cd(x,y)^\delta}{V(x,y,d(x,y))}
        \int_{\|\widehat x-\widehat y\|}^{\infty}t^{-1-\delta}\,dt\notag\\
        &=\frac{C}{\delta V(x,y,d(x,y))}
        \left(\frac{d(x,y)}{\|\widehat x-\widehat y\|}\right)^\delta.
        \label{eq:kernel-I3}
    \end{align}
    Since $0<d(x,y)/\|\widehat x-\widehat y\|\leq1$ and $\varepsilon<\delta<1$,
    \eqref{eq:kernel-I1}--\eqref{eq:kernel-I3} prove the absolute convergence
    and \eqref{eq:integrated-size}.

    We next prove \eqref{eq:integrated-y}. Suppose
    $0<\|\widehat y-\widehat{y'}\|\leq d(x,y)/2$; when the increment is zero, the assertion is immediate.
    The triangle inequality gives
    \[
    |d(x,y')-d(x,y)|\leq\|\widehat y-\widehat{y'}\|
    \qquad \mathrm{and} \qquad\bigl|\|\widehat x-\widehat{y'}\|-\|\widehat x-\widehat y\|\bigr|\leq\|\widehat y-\widehat{y'}\|.
    \]
    Consequently,
    \begin{equation}\label{eq:kernel-distance-comparison}
        \begin{aligned}
            \frac{d(x,y)}{2}&\leq d(x,y')\leq\frac{3d(x,y)}{2},\\
            \frac{\|\widehat x-\widehat y\|}{2}&\leq\|\widehat x-\widehat{y'}\|\leq\frac{3\|\widehat x-\widehat y\|}{2}.
        \end{aligned}
    \end{equation}
    Ball inclusion and the doubling condition then imply
    \[
    V(x,y',t+d(x,y'))\asymp V(x,y,t+d(x,y)),\qquad t>0.
    \]
    We decompose the scale integral into three parts:
    \begin{align}
        |K(x,y)-K(x,y')|
        &\leq\int_0^{\|\widehat y-\widehat{y'}\|}|\Psi_t(x,y)-\Psi_t(x,y')|\,\frac{dt}{t}\notag\\
        &\quad+\int_{\|\widehat y-\widehat{y'}\|}^{\|\widehat x-\widehat y\|}|\Psi_t(x,y)-\Psi_t(x,y')|\,\frac{dt}{t}\notag\\
        &\quad+\int_{\|\widehat x-\widehat y\|}^{\infty}|\Psi_t(x,y)-\Psi_t(x,y')|\,\frac{dt}{t}\notag\\
        &=:J_1+J_2+J_3.
        \label{eq:kernel-difference-ranges}
    \end{align}
    For the first part, condition \textup{(i)} gives
    \begin{align*}
        J_1
        &\leq C\int_0^{\|\widehat y-\widehat{y'}\|}
        \frac{dt}{V(x,y,t+d(x,y))(t+\|\widehat x-\widehat y\|)}\\
        &\quad+C\int_0^{\|\widehat y-\widehat{y'}\|}
        \frac{dt}{V(x,y',t+d(x,y'))(t+\|\widehat x-\widehat{y'}\|)}\\
        &\leq\frac{C}{\|\widehat x-\widehat y\|\,V(x,y,d(x,y))}\int_0^{\|\widehat y-\widehat{y'}\|} \,dt\\
        &=\frac{C}{V(x,y,d(x,y))}\frac{\|\widehat y-\widehat{y'}\|}{\|\widehat x-\widehat y\|}.
    \end{align*}
    For $t\geq\|\widehat y-\widehat{y'}\|$, condition \textup{(ii)} and the volume comparison yield
    \[
    |\Psi_t(x,y)-\Psi_t(x,y')|\leq
    C\left(\frac{\|\widehat y-\widehat{y'}\|}{t}\right)^\varepsilon
    \frac{t}{V(x,y,t+d(x,y))(t+\|\widehat x-\widehat y\|)}.
    \]
    Using this estimate on the middle interval, we have
    \begin{align*}
        J_2
        &\leq C\|\widehat y-\widehat{y'}\|^\varepsilon
        \int_{\|\widehat y-\widehat{y'}\|}^{\|\widehat x-\widehat y\|}
        \frac{t^{-\varepsilon}}{V(x,y,t+d(x,y))(t+\|\widehat x-\widehat y\|)}\,dt\\
        &\leq\frac{C\|\widehat y-\widehat{y'}\|^\varepsilon}{\|\widehat x-\widehat y\|\,V(x,y,d(x,y))}
        \int_{\|\widehat y-\widehat{y'}\|}^{\|\widehat x-\widehat y\|}t^{-\varepsilon}\,dt\\
        &=\frac{C}{(1-\varepsilon)V(x,y,d(x,y))}
        \left\{\left(\frac{\|\widehat y-\widehat{y'}\|}{\|\widehat x-\widehat y\|}\right)^\varepsilon
        -\frac{\|\widehat y-\widehat{y'}\|}{\|\widehat x-\widehat y\|}\right\}\\
        &\leq\frac{C}{(1-\varepsilon)V(x,y,d(x,y))}
        \left(\frac{\|\widehat y-\widehat{y'}\|}{\|\widehat x-\widehat y\|}\right)^\varepsilon.
    \end{align*}
    On the last interval, we use $t+\|\widehat x-\widehat y\|\geq t$ to obtain
    \begin{align*}
        J_3
        &\leq C\|\widehat y-\widehat{y'}\|^\varepsilon
        \int_{\|\widehat x-\widehat y\|}^{\infty}
        \frac{t^{-\varepsilon}}{V(x,y,t+d(x,y))(t+\|\widehat x-\widehat y\|)}\,dt\\
        &\leq\frac{C\|\widehat y-\widehat{y'}\|^\varepsilon}{V(x,y,d(x,y))}
        \int_{\|\widehat x-\widehat y\|}^{\infty}t^{-1-\varepsilon}\,dt\\
        &=\frac{C}{\varepsilon V(x,y,d(x,y))}
        \left(\frac{\|\widehat y-\widehat{y'}\|}{\|\widehat x-\widehat y\|}\right)^\varepsilon.
    \end{align*}
    Since $\|\widehat y-\widehat{y'}\|/\|\widehat x-\widehat y\|\leq1/2$, adding the three estimates proves
    \eqref{eq:integrated-y}.

    We next prove \eqref{eq:integrated-x}. Suppose
    $0<\|\widehat x-\widehat{x'}\|\leq d(x,y)/2$; when the increment is zero, the assertion is immediate.
    The triangle inequality gives
    \[
    |d(x',y)-d(x,y)|\leq\|\widehat x-\widehat{x'}\|
    \qquad \mathrm{and} \qquad\bigl|\|\widehat{x'}-\widehat y\|-\|\widehat x-\widehat y\|\bigr|\leq\|\widehat x-\widehat{x'}\|.
    \]
    Consequently,
    \[
    \begin{aligned}
        \frac{d(x,y)}{2}&\leq d(x',y)\leq\frac{3d(x,y)}{2}
       \end{aligned}\]
     and \[\begin{aligned}  \frac{\|\widehat x-\widehat y\|}{2}&\leq\|\widehat{x'}-\widehat y\|\leq\frac{3\|\widehat x-\widehat y\|}{2}.
    \end{aligned}
    \]
    Ball inclusion and the doubling condition then imply
    \[
    V(x',y,t+d(x',y))\asymp V(x,y,t+d(x,y)),\qquad t>0.
    \]
    We decompose the scale integral into three parts:
    \begin{align*}
        |K(x,y)-K(x',y)|
        &\leq\int_0^{\|\widehat x-\widehat{x'}\|}|\Psi_t(x,y)-\Psi_t(x',y)|\,\frac{dt}{t}\\
        &\quad+\int_{\|\widehat x-\widehat{x'}\|}^{\|\widehat x-\widehat y\|}|\Psi_t(x,y)-\Psi_t(x',y)|\,\frac{dt}{t}\\
        &\quad+\int_{\|\widehat x-\widehat y\|}^{\infty}|\Psi_t(x,y)-\Psi_t(x',y)|\,\frac{dt}{t}\\
        &=:J_1'+J_2'+J_3'.
    \end{align*}
    For the first part, condition \textup{(i)} gives
    \begin{align*}
        J_1'
        &\leq C\int_0^{\|\widehat x-\widehat{x'}\|}
        \frac{dt}{V(x,y,t+d(x,y))(t+\|\widehat x-\widehat y\|)}\\
        &\quad+C\int_0^{\|\widehat x-\widehat{x'}\|}
        \frac{dt}{V(x',y,t+d(x',y))(t+\|\widehat{x'}-\widehat y\|)}\\
        &\leq\frac{C}{\|\widehat x-\widehat y\|\,V(x,y,d(x,y))}\int_0^{\|\widehat x-\widehat{x'}\|}1\,dt\\
        &=\frac{C}{V(x,y,d(x,y))}\frac{\|\widehat x-\widehat{x'}\|}{\|\widehat x-\widehat y\|}.
    \end{align*}
    For $t\geq\|\widehat x-\widehat{x'}\|$, condition \textup{(iii)} and the volume comparison yield
    \[
    |\Psi_t(x,y)-\Psi_t(x',y)|\leq
    C\left(\frac{\|\widehat x-\widehat{x'}\|}{t}\right)^\varepsilon
    \frac{t}{V(x,y,t+d(x,y))(t+\|\widehat x-\widehat y\|)}.
    \]
    Using this estimate on the middle interval, we have
    \begin{align*}
        J_2'
        &\leq C\|\widehat x-\widehat{x'}\|^\varepsilon
        \int_{\|\widehat x-\widehat{x'}\|}^{\|\widehat x-\widehat y\|}
        \frac{t^{-\varepsilon}}{V(x,y,t+d(x,y))(t+\|\widehat x-\widehat y\|)}\,dt\\
        &\leq\frac{C\|\widehat x-\widehat{x'}\|^\varepsilon}{\|\widehat x-\widehat y\|\,V(x,y,d(x,y))}
        \int_{\|\widehat x-\widehat{x'}\|}^{\|\widehat x-\widehat y\|}t^{-\varepsilon}\,dt\\
        &=\frac{C}{(1-\varepsilon)V(x,y,d(x,y))}
        \left\{\left(\frac{\|\widehat x-\widehat{x'}\|}{\|\widehat x-\widehat y\|}\right)^\varepsilon
        -\frac{\|\widehat x-\widehat{x'}\|}{\|\widehat x-\widehat y\|}\right\}\\
        &\leq\frac{C}{(1-\varepsilon)V(x,y,d(x,y))}
        \left(\frac{\|\widehat x-\widehat{x'}\|}{\|\widehat x-\widehat y\|}\right)^\varepsilon.
    \end{align*}
    On the last interval, we use $t+\|\widehat x-\widehat y\|\geq t$ to obtain
    \begin{align*}
        J_3'
        &\leq C\|\widehat x-\widehat{x'}\|^\varepsilon
        \int_{\|\widehat x-\widehat y\|}^{\infty}
        \frac{t^{-\varepsilon}}{V(x,y,t+d(x,y))(t+\|\widehat x-\widehat y\|)}\,dt\\
        &\leq\frac{C\|\widehat x-\widehat{x'}\|^\varepsilon}{V(x,y,d(x,y))}
        \int_{\|\widehat x-\widehat y\|}^{\infty}t^{-1-\varepsilon}\,dt\\
        &=\frac{C}{\varepsilon V(x,y,d(x,y))}
        \left(\frac{\|\widehat x-\widehat{x'}\|}{\|\widehat x-\widehat y\|}\right)^\varepsilon.
    \end{align*}
    Since $\|\widehat x-\widehat{x'}\|/\|\widehat x-\widehat y\|\leq1/2$, adding the three estimates proves
    \eqref{eq:integrated-x}.
    This completes the proof.
\end{proof}

Lemma~\ref{lem:kernel-integration} yields the kernel estimates in
Definition~\ref{def:CZ}; the operator representation must be established
separately. Taking $\Psi_t=q_t=-t\partial_t p_t$,
Theorem~\ref{thm:poisson} verifies condition~\ref{eq:scale-size}, whereas
conditions~\ref{eq:scale-y} and~\ref{eq:scale-x} remain additional
assumptions.

The set $\{(x,y):d(x,y)=0\}$ has product measure zero. For $r,r_0>0$,
\begin{align*}
    &\int_{\substack{\|\widehat x-\widehat y\|\geq r_0\\0<d(x,y)<r}}
    |K(x,y)|\,dw_{k,a}(y)\\
    &\qquad\leq Cr_0^{-\varepsilon}
    \int_{0<d(x,y)<r}
    \frac{d(x,y)^\varepsilon}{w_{k,a}(B_a(x,d(x,y)))}\,dw_{k,a}(y)
    \leq C(r/r_0)^\varepsilon.
\end{align*}
Thus $K$ is locally integrable off the ordinary diagonal, and
\eqref{eq:offdiagonal-pairing} applies to disjoint supports even
when their reflected copies intersect.

\subsection{Boundedness of singular integral operators}
\begin{proof}[Proof of Theorem~\ref{thm:Lp}]
    We first prove weak type $(1,1)$. By normalization, we may assume
    \[
    C_K+\|T\|_{2\to2}\leq1.
    \]
    Fix $f\in L^1_{k,a}(\mathbb{R}^n)\cap L^2_{k,a}(\mathbb{R}^n)$ and $\lambda>0$.
    Let $\{Q_j\}_{j\geq1}$ be the maximal coordinate dyadic cubes such that
    \begin{equation}\label{eq:dyadic-selection}
        \lambda<\frac1{w_{k,a}(Q_j)}\int_{Q_j}|f(y)|\,dw_{k,a}(y).
    \end{equation}
The ancestor averages tend to zero owing to the $L^1_{k,a}(\mathbb
R^n)$-integrability of $f$ and  the volume growth of the dyadic
ancestors. Hence, the maximal cubes are well-defined and pairwise
disjoint. If $\widetilde Q_j$ denotes the dyadic parent, then the
doubling condition and maximality imply
\begin{equation}\label{eq:dyadic-average-upper}
    \frac{1}{w_{k,a}(Q_j)}\int_{Q_j}|f(y)|\,dw_{k,a}(y)
    \leq \frac{w_{k,a}(\widetilde Q_j)}{w_{k,a}(Q_j)}\lambda
    \leq C\lambda.
\end{equation}
    Here $C$ depends on the doubling constant.

    Put $\Omega=\bigcup_jQ_j$ and define
    \begin{equation}\label{eq:good-bad-definition}
        \begin{aligned}
            g(x)&=f(x)\mathbf1_{\Omega^c}(x)
            +\sum_{j=1}^\infty f_{Q_j}\mathbf1_{Q_j}(x),\\
            b_j(x)&=(f(x)-f_{Q_j})\mathbf1_{Q_j}(x),
            \qquad \mathrm{and } \qquad b=\sum_jb_j=f-g.
        \end{aligned}
    \end{equation}
   By construction   and \eqref{eq:dyadic-average-upper}, we have
    \begin{equation}\label{eq:CZ-decomposition}
        f=g+\sum_{j\geq1}b_j,\qquad \supp b_j\subset Q_j
        \qquad \mathrm{and } \qquad\int_{Q_j}b_j(y)\,dw_{k,a}(y)=0.
    \end{equation}
    The good part satisfies $|g|\leq C\lambda$ almost everywhere and
    \[
    \|g\|_1\leq\|f\|_1 \qquad \mathrm{and } \qquad
    \|g\|_2^2\leq\|g\|_\infty\|g\|_1
    \leq C\lambda\|f\|_1.
    \]
    By disjointness and \eqref{eq:dyadic-selection},
    \[
    \sum_{j\geq1}w_{k,a}(Q_j)
    \leq\lambda^{-1}\sum_{j\geq1}\int_{Q_j}|f(y)|\,dw_{k,a}(y)
    \leq\lambda^{-1}\|f\|_1.
    \]
    For each bad part,
    \[
    \|b_j\|_1\leq2\int_{Q_j}|f(y)|\,dw_{k,a}(y)
    \leq C\lambda w_{k,a}(Q_j).
    \]
    Summing over the disjoint cubes also gives
    \[
    \sum_{j\geq1}\|b_j\|_1\leq2\|f\|_1.
    \]
    Jensen's inequality furthermore gives
    \[
    \sum_{j=1}^\infty\int_{Q_j}|b_j(y)|^2\,dw_{k,a}(y)
    \leq4\int_\Omega|f(y)|^2\,dw_{k,a}(y)<\infty.
    \]
    Hence the bad-function series converges in $L^2_{k,a}(\mathbb{R}^n)$ as well as in $L^1_{k,a}(\mathbb{R}^n)$.

    Let $x_j$ be the coordinate center of $Q_j$, and choose $r_j$
    comparable to its side length so that
    $Q_j\subset B_j=B_a(x_j,r_j)$ and
    $w_{k,a}(B_j)\leq Cw_{k,a}(Q_j)$. Define the reflected exceptional
    set
    \begin{equation}\label{eq:exceptional-set}
        \Omega^*=\bigcup_j\mathcal O(8B_j).
    \end{equation}
    Reflection invariance and doubling yield
    \begin{align}
        w_{k,a}(\Omega^*)&\leq
        \sum_{j=1}^\infty\sum_{\sigma\in G}w_{k,a}(8\sigma B_j)
        \leq C|G|\sum_jw_{k,a}(Q_j)\notag\\
        &\leq C\lambda^{-1}\|f\|_1.
        \label{eq:exceptional-measure}
    \end{align}
    For $x\notin\Omega^*$ and $y\in Q_j$,
    \[
    d(x,x_j)>8r_j,\qquad
    \|\widehat y-\widehat{x_j}\|\leq r_j,\qquad \mathrm{and} \qquad
    d(x,y)\asymp d(x,x_j).
    \]
    Using cancellation of $b_j$ in the kernel representation,
    \begin{equation}\label{eq:bad-cancellation}
        Tb_j(x)=\int_{Q_j}[K(x,y)-K(x,x_j)]b_j(y)\,dw_{k,a}(y).
    \end{equation}
    For   $b_j$, the identity follows by $L^2_{k,a}(\mathbb R^n)$
    approximation and local square integrability of the kernel on
    these separated sets.

    The regularity estimate gives
    \begin{align}
        |Tb_j(x)|
        &\leq\int_{Q_j}
        \frac{C}{V(x,x_j,d(x,x_j))}
        \left(\frac{\|\widehat y-\widehat{x_j}\|}
        {d(x,x_j)}\right)^\varepsilon
        |b_j(y)|\,dw_{k,a}(y)\notag\\
        &\leq\frac{C\|b_j\|_1}{w_{k,a}(B_a(x_j,d(x,x_j)))}
        \left(\frac{r_j}{d(x,x_j)}\right)^\varepsilon.
        \label{eq:bad-pointwise}
    \end{align}
    Integrating over the orbit annuli
    $2^\ell r_j<d(x,x_j)\leq2^{\ell+1}r_j$, $\ell\geq3$,
    \begin{align}
        \int_{(\Omega^*)^c}|Tb_j(x)|\,dw_{k,a}(x)
        &\leq C\|b_j\|_1
        \sum_{\ell=3}^\infty
        \frac{2^{-\ell\varepsilon}
            w_{k,a}\{x:d(x,x_j)\leq2^{\ell+1}r_j\}}
        {w_{k,a}(B_a(x_j,2^\ell r_j))}\notag\\
        &\leq C\|b_j\|_1\sum_{\ell=3}^\infty2^{-\ell\varepsilon}
        \leq C\|b_j\|_1.
        \label{eq:bad-annuli}
    \end{align}
For finite sums, the estimate follows from the triangle inequality.
Since the partial sums $\sum_{j=1}^N Tb_j$ converge to $Tb$ in
$L^2_{k,a}(\mathbb{R}^n)$, a subsequence converges almost everywhere.
Applying Fatou's lemma along this subsequence gives
    \begin{equation}\label{eq:bad-integral}
        \int_{(\Omega^*)^c}|Tb(x)|\,dw_{k,a}(x)
        \leq C\sum_{j=1}^\infty\|b_j\|_1\leq C\|f\|_1.
    \end{equation}
Combining the estimates for the good and bad parts, we obtain
    \begin{align}
        w_{k,a}\{|Tf|>\lambda\}
        &\leq w_{k,a}(\Omega^*)+
        w_{k,a}\{|Tg|>\lambda/2\}\notag\\
        &\quad+w_{k,a}\{x\notin\Omega^*:|Tb(x)|>\lambda/2\}\notag\\
        &\leq C\lambda^{-1}\|f\|_1
        +4\lambda^{-2}\|Tg\|_2^2
        +2\lambda^{-1}\int_{(\Omega^*)^c}|Tb(x)|\,dw_{k,a}(x)\notag\\
        &\leq C\lambda^{-1}\|f\|_1.
        \label{eq:weak11}
    \end{align}
Approximation in $L^1_{k,a}(\mathbb{R}^n)$ yields the weak-type extension.

The Marcinkiewicz interpolation theorem gives
$L^p_{k,a}(\mathbb{R}^n)$ boundedness for $1<p<2$ from the weak-type
$(1,1)$ estimate and the $L^2$ bound. The adjoint has kernel
$\overline{K(y,x)}$. Symmetry of $V$ and $d$ shows that it satisfies
the same conditions. Thus the preceding argument applies to $T^*$
and, for $2<p<\infty$,
\[
\|Tf\|_p
=\sup_{\|g\|_{p'}\leq1}|\langle f,T^*g\rangle|
\leq C_p\|f\|_p.
\]
Density yields the unique bounded extensions. Restoring the normalization proves the asserted $L^p_{k,a}(\mathbb{R}^n)$ conclusions.

    Let $f\in L^\infty_{k,a}(\mathbb R^n)$. Choose $x_*$ outside the
    exceptional null set of the kernel estimates. For $N>1$, set
    \[
    f_N(y)=f(y)\mathbf1_{\{d(x_*,y)<N\}} \qquad \mathrm{and} \qquad
    c_N(f)=\int_{1<d(x_*,y)<N}K(x_*,y)f(y)\,dw_{k,a}(y).
    \]
    The truncation has bounded support, so $f_N\in L^2_{k,a}(\mathbb{R}^n)$. Moreover,
    \begin{align*}
        |c_N(f)|
        &\leq C C_K\|f\|_\infty
        \sum_{0\leq j\leq\lceil\log_2N\rceil}
        \frac{w_{k,a}\{y:d(x_*,y)<2^{j+1}\}}
        {w_{k,a}(B_a(x_*,2^j))}\\
        &\leq C C_K(1+\log N)\|f\|_\infty<\infty.
    \end{align*}

    Let $R>1$, $N>4R$, and let
    $f^{(R)}(y)=f(y)\mathbf1_{\{d(x_*,y)<4R\}}$. For
    $x\in B_a(x_*,R)$, the separated representation gives
    \begin{align}
        Tf_N(x)-c_N(f)
        ={}&Tf^{(R)}(x)
        -\int_{1<d(x_*,y)<4R}K(x_*,y)f(y)\,dw_{k,a}(y)\notag\\
        &+\int_{4R<d(x_*,y)<N}
        [K(x,y)-K(x_*,y)]f(y)\,dw_{k,a}(y).
        \label{eq:BMO-truncation}
    \end{align}
    For $M>N>4R$ and almost every $x\in B_a(x_*,R)$,
    \begin{align*}
        &|[Tf_M(x)-c_M(f)]-[Tf_N(x)-c_N(f)]|\\
        &=\left|\int_{N\leq d(x_*,y)<M}
        [K(x,y)-K(x_*,y)]f(y)\,dw_{k,a}(y)\right|\\
        &\leq C C_K\|f\|_\infty
        \int_{d(x_*,y)\geq N}
        \frac{(R/d(x_*,y))^\varepsilon}
        {w_{k,a}(B_a(x_*,d(x_*,y)))}\,dw_{k,a}(y)\\
        &\leq C C_K\|f\|_\infty(R/N)^\varepsilon
        \sum_{j=0}^\infty2^{-j\varepsilon}.
    \end{align*}
    Thus $Tf_N-c_N(f)$ converges in
    $L^2_{\mathrm{loc}}(\mathbb R^n,dw_{k,a})$, with uniformly
    convergent tail differences. Its limit is a representative of $Tf$ and we show next that it lies in  BMO.

    Now let $B=B_a(x_0,r)$ be arbitrary, set
    $E=\mathcal O(8B)$, and split
    \[
    f=f_1+f_2,\qquad f_1=f\mathbf1_E,
    \qquad f_2=f\mathbf1_{E^c}.
    \]
    For $x\in B$ define the far part by
    \begin{equation}\label{eq:BMO-far-definition}
        T_Bf_2(x)=\int_{E^c}[K(x,y)-K(x_0,y)]f(y)\,dw_{k,a}(y).
    \end{equation}
    The above integral converges absolutely. Indeed,
    \begin{align}
        |T_Bf_2(x)|
        &\leq C C_K\|f\|_\infty
        \int_{d(x_0,y)>8r}
        \frac{(r/d(x_0,y))^\varepsilon}
        {w_{k,a}(B_a(x_0,d(x_0,y)))}\,dw_{k,a}(y)\notag\\
        &\leq C C_K\|f\|_\infty
        \sum_{j=3}^\infty2^{-j\varepsilon}
        \leq C C_K\|f\|_\infty.
        \label{eq:BMO-far-bound}
    \end{align}
    Comparing the finite truncations in \eqref{eq:BMO-truncation}
    with this decomposition, then passing to the limit, gives
    $Tf=Tf_1+T_Bf_2+c_B$ on $B$, where $c_B$ is constant.

    Reflection invariance and the doubling condition imply
    \[
    w_{k,a}(E)\leq |G|w_{k,a}(8B)\leq Cw_{k,a}(B).
    \]
    The local term therefore satisfies
    \begin{align}
        \left(\frac1{w_{k,a}(B)}\int_B|Tf_1(x)|^2\,dw_{k,a}(x)\right)^{1/2}
        &\leq\frac{\|T\|_{2\to2}}{w_{k,a}(B)^{1/2}}
        \|f_1\|_2\notag\\
        &\leq\|T\|_{2\to2}
        \left(\frac{w_{k,a}(E)}{w_{k,a}(B)}\right)^{1/2}\|f\|_\infty
        \leq C\|T\|_{2\to2}\|f\|_\infty.
        \label{eq:BMO-local-bound}
    \end{align}
  By \eqref{eq:BMO-far-bound} and
    \eqref{eq:BMO-local-bound},
    \begin{align*}
        \left(\frac1{w_{k,a}(B)}\int_B|Tf(x)-c_B|^2\,dw_{k,a}(x)\right)^{1/2}
        &\leq C(\|T\|_{2\to2}+C_K)\|f\|_\infty,\\   \intertext{and}
        \frac1{w_{k,a}(B)}\int_B|Tf(x)-(Tf)_B|\,dw_{k,a}(x)
        &\leq2\left(\frac1{w_{k,a}(B)}
        \int_B|Tf(x)-c_B|^2\,dw_{k,a}(x)\right)^{1/2}.
    \end{align*}
    Taking the supremum over $B$ proves
    \begin{equation}\label{eq:BMO-operator-bound}
        \|Tf\|_{\BMO(\mathbb{R}^n,dw_{k,a})}
        \leq C(\|T\|_{2\to2}+C_K)\|f\|_\infty.
    \end{equation}
    This completes the proof of Theorem~\ref{thm:Lp}.
\end{proof}

The same conclusions hold under the following integral regularity
conditions.
\begin{corollary}\label{cor:Hormander}
    The conclusions of Theorem~\ref{thm:Lp} remain valid for an $L^2_{k,a}(\mathbb{R}^n)$-bounded operator admitting the kernel representation and size estimate above, provided that the two difference conditions are replaced by
    \begin{equation}\label{eq:Hormander-y}
        \sup_{y,y'}\int_{d(x,y)>2\|\widehat y-\widehat{y'}\|}
        |K(x,y)-K(x,y')|\,dw_{k,a}(x)<\infty
    \end{equation}
    and
    \begin{equation}\label{eq:Hormander-x}
        \sup_{x,x'}\int_{d(x,y)>2\|\widehat x-\widehat{x'}\|}
        |K(x,y)-K(x',y)|\,dw_{k,a}(y)<\infty.
    \end{equation}
    The suprema are interpreted as essential suprema for kernel formulations.
    The operator norms may also depend on the constants in
    \eqref{eq:Hormander-y} and \eqref{eq:Hormander-x}.
\end{corollary}

\begin{proof}
    Since the conditions use essential suprema, the reference points
    must be chosen by Fubini's theorem.
    For each cube $Q_j$, choose $z_j\in Q_j$ such that
    \eqref{eq:Hormander-y} holds with $y'=z_j$ for almost every $y$.
    Keep the geometric centers used in the definition of $B_j$ and
    $\Omega^*$ unchanged. If $x\notin\Omega^*$ and $y\in Q_j$, then
    \[
    d(x,y)>7r_j
    \qquad \mathrm{and} \qquad \|\widehat y-\widehat{z_j}\|\leq2r_j.
    \]
    Using $\int b_j\,dw_{k,a}=0$, we obtain
    \[
    Tb_j(x)=\int_{Q_j}[K(x,y)-K(x,z_j)]b_j(y)\,dw_{k,a}(y).
    \]
    Tonelli's theorem and \eqref{eq:Hormander-y} give
    \[
    \int_{(\Omega^*)^c}|Tb_j(x)|\,dw_{k,a}(x)
    \leq C\|b_j\|_1.
    \]
    Together with \eqref{eq:CZ-decomposition} and the $L^2_{k,a}$
    estimate for the good part, this proves the weak type $(1,1)$ bound.
    Applying the same argument to $T^*$, using \eqref{eq:Hormander-x},
    and then using interpolation and duality gives $L^p_{k,a}$
    boundedness for every $1<p<\infty$.

    We next prove the BMO estimate and identify a global representative.
    Fix $f\in L^\infty_{k,a}(\mathbb R^n)$.
    By Fubini's theorem, choose $a_*$ such that
    \eqref{eq:Hormander-x} holds for $(x,x')=(x,a_*)$ for almost every $x$.
    Fix a coordinate ball $B_*$ centered at $a_*$ and put
    $E_*=\mathcal O(8B_*)$.
    For any $B=B_a(x_0,r)$, put $E_B=\mathcal O(8B)$ and choose
    $a_B\in B$ so that \eqref{eq:Hormander-x} holds both for
    $(x,a_B)$ for almost every $x$ and for the pair $(a_B,a_*)$.
    Set
    \[
    F_B(x)=T(f\mathbf1_{E_B})(x)
    +\int_{E_B^c}[K(x,y)-K(a_B,y)]f(y)\,dw_{k,a}(y),
    \qquad x\in B.
    \]
    For $x\in B$ and $y\notin E_B$,
    \[
    d(x,y)>7r
    \qquad \mathrm{and} \qquad\|\widehat x-\widehat{a_B}\|\leq2r.
    \]
    Thus the second term in $F_B$ is absolutely convergent for almost
    every $x\in B$ and has $L^\infty(B)$ norm at most $C\|f\|_\infty$.
    The local term has the same mean-oscillation bound as in the
    proof of Theorem~\ref{thm:Lp}.

    For completeness, let $A_R=B_a(0,R)$, $f_R=f\mathbf1_{A_R}$, and
    \[
    c_R=\int_{A_R\cap E_*^c}K(a_*,y)f(y)\,dw_{k,a}(y).
    \]
    The difference
    \[
    \int_{A_R\cap E_B^c}K(a_B,y)f(y)\,dw_{k,a}(y)-c_R
    \]
    converges to a finite constant $c_B$ as $R\to\infty$.
    Indeed, the far tail is absolutely integrable by the condition
    for $(a_B,a_*)$, while the remaining bounded region is separated
    from the relevant reference orbits and is controlled by the size
    estimate. For large $R$, $A_R$ contains $E_B$.
    The difference tails defining $F_B$ converge almost everywhere
    and are uniformly bounded by the integral regularity condition.
    Dominated convergence therefore gives
    \[
    Tf_R-c_R\longrightarrow F_B+c_B\quad\text{in }L^2(B,dw_{k,a}).
    \]
    These limits agree on overlaps and define a global function
    modulo constants. Its mean oscillation over every coordinate ball
    is bounded by $C\|f\|_\infty$, as required.
\end{proof}

\section{\texorpdfstring{The $T1$ theorem}{The T1 theorem}}\label{sec:T1}
In this section we prove Theorem~\ref{thm:T1-intro}. We first establish the Meyer-type commutation estimates and the boundedness of $T$
on smooth molecules. These estimates yield almost orthogonality, and the proof is completed via reproducing formulae
 and paraproducts. We follow the real-variable argument of \cite[Section~2]{THHLL}; see also \cite{Hd,35} for the classical background.
 Throughout this section, $T$ is given  as in Definition~\ref{def:CZ}, and $C_K$ denotes its kernel constant. Distances and volumes are those from Section~\ref{sec:coordinate-geometry}, and all integrals are taken with respect to $dw_{k,a}$.

\subsection{Meyer-type commutation estimates}\label{sec:constant-tests}
We begin with the weak boundedness property and the constant   test
functions.
\begin{definition}\label{def:WBP}
    Fix $0<\eta\leq1$. An operator is said to have the weak boundedness
    property if
    \begin{equation}\label{eq:WBP-joint}
        |\langle\mathcal K,F\rangle|
        \leq C_W\max\{w_{k,a}(B_a(x_0,r)),w_{k,a}(B_a(y_0,r))\}
    \end{equation}
   for every \(F(x,y)=\Phi(\widehat x,\widehat y)\) with
   \(\Phi\in C_c^\infty(\mathbb R^n\times\mathbb R^n)\),
   \(\operatorname{supp}F\subset B_a(x_0,r)\times B_a(y_0,r)\) and satisfying
    \begin{enumerate}[label=\textup{(\roman*)},ref=\thetheorem(\roman*),leftmargin=2.2em,labelsep=0.55em,itemsep=6pt,topsep=6pt,parsep=0pt]
        \item\label{eq:WBP-joint-normalization} $\|F\|_\infty\leq1$;
        \item $\displaystyle \sup_y[F(\cdot,y)]_\eta\leq r^{-\eta}$;
        \item $\displaystyle \sup_x[F(x,\cdot)]_\eta\leq r^{-\eta}$.
    \end{enumerate}
    We write $T\in\WBP$; the least constant for the fixed exponent
    is denoted by $\|T\|_{\WBP}$.
\end{definition}
We use the formulation of \cite[Definition~1.2]{THHLL}.
The joint pairing is initially defined on smooth test functions.
It extends by continuity to their closure in the displayed
separate H\"older norms. Extensions to less regular functions
are constructed below using the local formula.
In particular, let $\varphi,\psi$ be coordinate-smooth test
functions satisfying
\[
\operatorname{supp}\varphi\subset B_a(y_0,r) \qquad \mathrm{and}
\qquad \operatorname{supp}\psi\subset B_a(x_0,r);
\]
\[
\|\varphi\|_\infty,\|\psi\|_\infty\leq1 \qquad \mathrm{and} \qquad
[\varphi]_\eta,[\psi]_\eta\leq r^{-\eta}.
\]
Then
\begin{equation}\label{eq:WBP-product}
    |\langle T\varphi,\psi\rangle|
    \leq C\|T\|_{\WBP}
    \max\{w_{k,a}(B_a(x_0,r)),w_{k,a}(B_a(y_0,r))\}.
\end{equation}
The requirement does not include an $L^2_{k,a}(\mathbb{R}^n)$ bound.

Let $g$ be a mean-zero test function supported in $B_a(x_0,r)$.
Choose a cutoff $\chi$, smooth in the transformed coordinates, such that
\[
0\leq\chi\leq1,\qquad \chi=1\ \text{on }\mathcal
O(B_a(x_0,4r)),\qquad \mathrm{and} \qquad
\operatorname{supp}\chi\subset\mathcal O(B_a(x_0,8r)).
\]
Define
\begin{equation}\label{eq:T1-definition}
    \begin{aligned}
        \langle T1,g\rangle={}&\langle T\chi,g\rangle +\int_{\mathbb R^n}\overline{g(x)}\int_{\mathbb R^n}
        [K(x,y)-K(x_0,y)](1-\chi(y))\,dw_{k,a}(y)\,dw_{k,a}(x).
    \end{aligned}
\end{equation}
The tail is absolutely convergent, since
\begin{align}
    &\sup_{x\in B_a(x_0,r)}
    \int_{d(x_0,y)>4r}|K(x,y)-K(x_0,y)|\,dw_{k,a}(y)\notag\\
    &\quad\leq C C_K\int_{d(x_0,y)>4r}
    \frac{(r/d(x_0,y))^\varepsilon}
    {w_{k,a}(B_a(x_0,d(x_0,y)))}\,dw_{k,a}(y)
    \leq C C_K\sum_{j=2}^\infty2^{-j\varepsilon}.
    \label{eq:T1-tail}
\end{align}
For two cutoffs, their difference vanishes near the support and its reflections. The separated kernel representation and $\int_{\mathbb R^n} g=0$ therefore cancel the difference of the two expressions in
\eqref{eq:T1-definition}. This proves independence of $\chi$.
The distribution $T^*1$ is defined in the same way, on mean-zero test functions.

We first establish the necessity of these testing conditions.
\begin{proposition}\label{prop:necessity}
    If $T$ in Definition~\ref{def:CZ} is bounded on $L^2_{k,a}(\mathbb{R}^n)$, then
    \begin{equation}\label{eq:necessity-bound}
        \|T1\|_{\BMO}+\|T^*1\|_{\BMO}+\|T\|_{\WBP}
        \leq C(\|T\|_{2\to2}+C_K).
    \end{equation}
\end{proposition}
\begin{proof}
    The BMO assertions follow from Theorem~\ref{thm:Lp} applied to
    $T$ and $T^*$ with input $1$. Explicitly, on $B=B_a(x_0,r)$ the distribution in
    \eqref{eq:T1-definition} has the following representative:
    \[
    T1=T\chi+F_B+c_B \qquad \mathrm{and} \qquad
    F_B(x)=\int_{\mathbb
    R^n}[K(x,y)-K(x_0,y)](1-\chi(y))\,dw_{k,a}(y),
    \]
where $c_B$ is constant. Equation \eqref{eq:T1-tail} and the
$L^2_{k,a}(\mathbb{R}^n)$ bound give
    \begin{align*}
        \|F_B\|_{L^\infty(B,dw_{k,a})}&\leq C C_K,\\
        \frac1{w_{k,a}(B)}\int_B|T\chi|\,dw_{k,a}
        &\leq \|T\|_{2\to2}
        \left(\frac{w_{k,a}(\supp\chi)}{w_{k,a}(B)}\right)^{1/2}
        \leq C\|T\|_{2\to2}.
    \end{align*}
    Subtracting the mean and taking the supremum proves the first
    two bounds.

It remains to verify the weak boundedness property for joint test functions. We first show that, for a smooth compactly supported function $\mathrm{V}$ which vanishes on the ordinary diagonal,
\begin{equation}\label{eq:L2-diagonal-identity}
    \langle\mathcal K,\mathrm{V}\rangle
    =\iint_{\mathbb R^n\times\mathbb R^n} K(x,y)\mathrm{V}(x,y)\,dw_{k,a}(y)\,dw_{k,a}(x).
\end{equation}
We work in coordinate variables on a fixed compact set and introduce a cutoff function to ensure $\|\widehat x-\widehat y\|<2h$.
 A bounded-overlap partition by balls of radius $h$ reduces this part to $O(1)$ neighboring product balls for each center.
  Since $\mathrm{V}(x,x)=0$, its rescaled pieces, divided by $h$, have uniformly bounded derivatives of every fixed order.
  After a further fixed smooth cutoff, they admit Fourier series expansions on a fixed product
  cube:
\[
h^{-1}\mathrm{V}_i(\phi^{-1}(\widehat{x_i}+hX),\phi^{-1}(\widehat{x_i}+hY))
=\sum_{p,q\in\mathbb Z^n}c_{i,p,q}e^{ip\cdot X}e^{iq\cdot Y}\] and
\[ \qquad |c_{i,p,q}|\leq C_N(1+\|p\|+\|q\|)^{-N}.
\]
The constant is independent of $h$ and $i$. Each exponential factor together with its cutoff has $L^2_{k,a}(\mathbb{R}^n)$ norm at most a fixed multiple of the square root of the local ball measure. For $N>2n+2$,
\[
|\langle\mathcal K,\mathrm{V}_i\rangle|
\leq C\|T\|_{2\to2}h\,w_{k,a}(CB_i)
\sum_{p,q}(1+\|p\|+\|q\|)^{-N}.
\]

    Bounded overlap gives a total $O(h)$ tending to zero. Away from
    the ordinary diagonal, the kernel representation applies. The
    integral converges by the coordinate-decay size estimate, the bound
    $|\mathrm{V}(x,y)|\leq C\|\widehat x-\widehat y\|^\theta$ on compact sets,
    $0<\theta<\varepsilon$, and \eqref{eq:inner-annuli}.
    This proves \eqref{eq:L2-diagonal-identity}.

    Now let $F$ satisfy conditions (i)--(iii) of Definition~\ref{def:WBP}. Write
    $R_0=\|\widehat{x_0}-\widehat{y_0}\|$ and $d_0=d(x_0,y_0)$.
    If $R_0\leq10r$, choose a cutoff $\chi$, smooth in the transformed
    coordinates, equal to one on both supports, and supported in
    $B_a(x_0,20r)$, with
    \[
    0\leq\chi\leq1\qquad \mathrm{and} \qquad [\chi]_\eta\leq Cr^{-\eta}.
    \]
    Then
    \[
    F(x,y)=F(x,x)\chi(y)+\mathrm{V}(x,y)\qquad \mathrm{and} \qquad \mathrm{V}(x,x)=0.
    \]
    Choose $0<\theta<\min\{\eta,\varepsilon\}$. Separate H\"older
    regularity gives, on its support,
    \[
    |\mathrm{V}(x,y)|\leq C\left(\frac{\|\widehat x-\widehat y\|}{r}\right)^\theta.
    \]
    The diagonal product is bounded by
    \[
    |\langle T\chi,\overline{F(\cdot,\cdot)|_{\mathrm{diag}}}\rangle|
    \leq C\|T\|_{2\to2}w_{k,a}(B_a(x_0,r)).
    \]
    For the remaining part, \eqref{eq:L2-diagonal-identity} yields
    \begin{align*}
        |\langle\mathcal K,\mathrm{V}\rangle|
        &\leq\frac{C C_K}{r^\theta}
        \int_{B_a(x_0,Cr)}\int_{d(x,y)<Cr}
        \frac{d(x,y)^\theta}{w_{k,a}(B_a(x,d(x,y)))}
        \,dw_{k,a}(y)\,dw_{k,a}(x)\\
        &\leq C C_K w_{k,a}(B_a(x_0,r)).
    \end{align*}
    If $R_0>10r$, the two supports are disjoint. When $d_0\leq4r$,
    \begin{align*}
        |\langle\mathcal K,F\rangle|
        &\leq C C_K R_0^{-\varepsilon}
        \int_{B_a(x_0,r)}\int_{d(x,y)<6r}
        \frac{d(x,y)^\varepsilon}{w_{k,a}(B_a(x,d(x,y)))}
        \,dw_{k,a}(y)\,dw_{k,a}(x)\\
        &\leq C C_K(r/R_0)^\varepsilon w_{k,a}(B_a(x_0,r)).
    \end{align*}
    When $d_0>4r$, $d(x,y)\asymp d_0$ on the supports, and
    \begin{align*}
        |\langle\mathcal K,F\rangle|
        &\leq C C_K\frac{w_{k,a}(B_a(x_0,r))w_{k,a}(B_a(y_0,r))}
        {V(x_0,y_0,d_0)}\\
        &\leq C C_K\max\{w_{k,a}(B_a(x_0,r)),w_{k,a}(B_a(y_0,r))\}.
    \end{align*}
    These three cases prove \eqref{eq:WBP-joint}. They also show that
    the estimates depend only on the displayed separate H\"older norms,
    so the joint pairing extends to limits of smooth joint test functions  in
    these norms.
\end{proof}

For the sufficiency argument we normalize
\begin{equation}\label{eq:testing-normalization}
    C_K+\|T\|_{\WBP}\leq1.
\end{equation}
 Smooth orbit cutoffs
may be constructed from a fixed radial cutoff $\vartheta$ by
\[
1-\prod_{g\in G}
\left[1-\vartheta\left(\frac{\|\widehat x-g\widehat{x_0}\|}{r}\right)\right].
\]
They equal one on a prescribed orbit ball, are supported in a fixed
enlargement, and have derivatives bounded at scale $r$.

We next establish the commutation identity needed for the local
boundedness estimate.
\begin{lemma}\label{lem:commutator}
    For compactly supported functions \(f,g,\chi\) that are smooth in the transformed
    coordinates \(\widehat x=\phi(x)\),
    \begin{equation}\label{eq:commutator}
        \begin{aligned}
            &\langle T(\chi f),g\rangle-\langle T\chi,\overline f g\rangle\\
            &\quad=\iint_{\mathbb R^n\times\mathbb R^n} K(x,y)\chi(y)[f(y)-f(x)]\overline{g(x)}
            \,dw_{k,a}(y)\,dw_{k,a}(x).
        \end{aligned}
    \end{equation}
   Moreover,  the above integral is absolutely convergent.
\end{lemma}
\begin{proof}
    For $0<\theta<\varepsilon$, the  decay size estimate gives
    \begin{equation}\label{eq:commutator-domination}
        \begin{aligned}
            |K(x,y)|\,|f(y)-f(x)|
            &\leq\frac C{w_{k,a}(B_a(x,d(x,y)))}
            \left(\frac{d(x,y)}{\|\widehat x-\widehat y\|}\right)^\theta
            [f]_\theta\|\widehat x-\widehat y\|^\theta\\
            &=\frac{C[f]_\theta d(x,y)^\theta}
            {w_{k,a}(B_a(x,d(x,y)))}.
        \end{aligned}
    \end{equation}
    The inner-annulus estimate proves local integrability.
Let $\lambda_h(x,y)$ be a smooth cutoff equal to one when
    $\|\widehat x-\widehat y\|\leq h$ and zero when
    $\|\widehat x-\widehat y\|\geq2h$.
    The part with $1-\lambda_h$ is represented by the integral
    and converges by \eqref{eq:commutator-domination}.

To treat the remaining distributional term, take a bounded-overlap partition $\eta_j(y)$ at coordinate scale $h$. Its nonzero pieces are supported in product sets $B_a(x_j,4h)\times B_a(x_j,2h)$.

    For
    \[
    F_j(x,y)=\lambda_h(x,y)\eta_j(y)\chi(y)
    [f(y)-f(x)]\overline{g(x)},
    \]
    smoothness and the product rule give
    \begin{align*}
        \|F_j\|_\infty&\leq C_fh
       \end{align*}
     and \begin{align*}  \sup_y[F_j(\cdot,y)]_\eta+
        \sup_x[F_j(x,\cdot)]_\eta&\leq C_{f,g,\chi}h^{1-\eta}.
    \end{align*}
    Indeed, for a difference shorter than $h$, the gradient is uniformly
    bounded, and for a longer difference one uses the supremum bound.
    Applying joint WBP to $F_j/(Ch)$ implies
    \begin{equation}\label{eq:diagonal-error}
        \left|\sum_j\langle\mathcal K,F_j\rangle\right|
        \leq Ch\sum_jw_{k,a}(B_a(x_j,4h))
        \leq Ch\,w_{k,a}(K_0)\longrightarrow0.
    \end{equation}
    The compact set $K_0$ is fixed. The last inequality uses bounded
    overlap and doubling, not boundedness of the density near the
    origin. This proves \eqref{eq:commutator}.
\end{proof}

The constant test and weak boundedness give the following estimate
for an orbit cutoff.
\begin{lemma}\label{lem:cutoff}
   Suppose $T1=0$. Let $C_0\geq4$ be fixed, and let $\chi$ be a
   cutoff that equals one on $\mathcal O(B_a(x_0,C_0r))$, is supported
   in a fixed enlargement of this set, and is normalized at scale
   $r$. Then
    \begin{equation}\label{eq:cutoff-bound}
        \|T\chi\|_{L^\infty(B_a(x_0,r),dw_{k,a})}\leq C.
    \end{equation}
\end{lemma}
\begin{proof}
    For mean-zero $g$ supported in $B=B_a(x_0,r)$, the definition of
    $T1$ and the cancellation assumption give
$$
        \langle T\chi,g\rangle
        =-\int_{\mathbb R^n}\overline{g(x)}\int_{\mathbb R^n}
        [K(x,y)-K(x_0,y)](1-\chi(y))\,dw_{k,a}(y)\,dw_{k,a}(x), $$ and $$
        |\langle T\chi,g\rangle|\leq C\|g\|_1.
$$
    Thus $T\chi=F+c_B$ on $B$ with $\|F\|_\infty\leq C$.
    Choose $0\leq\psi\leq1$, supported in $B$, normalized at scale
    $r$, and with $\int_{B}\psi\,dw_{k,a}\geq c\,w_{k,a}(B)$.
    A smooth partition of the orbit cutoff has at most $|G|$ terms
    centered on the orbit, with uniformly bounded scale norms.
    Joint WBP applied to their products with $\psi$ gives $|\langle T\chi,\psi\rangle|\le
    Cw_{k,a}(B)$ and
    \begin{align*}
        |c_B|\int_{B}\psi\,dw_{k,a}
        &\leq|\langle T\chi,\psi\rangle|
        +\int_{B}|F\psi|\,dw_{k,a}\\
        &\leq C\sum_{g\in G}
        \max\{w_{k,a}(B_a(gx_0,Cr)),w_{k,a}(B)\}
        \leq Cw_{k,a}(B).
    \end{align*}
Hence $|c_B|\leq C$, which yields the estimate.
\end{proof}

We are now in a position to prove the Meyer-type commutation
estimate, corresponding to \cite[Lemma~2.2]{THHLL}.
\begin{lemma}\label{lem:Meyer}
    Suppose $T1=0$ and $T\in\WBP$. If $\supp f\subset B_a(x_0,r)$,
    $\|f\|_\infty\leq1$ and $[f]_\theta\leq r^{-\theta}$ for some
    $0<\theta<\varepsilon$, then for each fixed $M>1$,
    \begin{equation}\label{eq:Meyer}
        \|Tf\|_{L^\infty(B_a(x_0,Mr),dw_{k,a})}\leq C_M.
    \end{equation}
    The statement holds initially for smooth functions and extends
    by local H\"older approximation.
\end{lemma}
\begin{proof}
    Choose $\chi=1$ on an orbit enlargement of
    $B_a(x_0,4Mr)$ that also contains $\supp f$. The commutation identity gives the
    local representative
    \begin{equation}\label{eq:local-T-formula}
        Tf(x)=\int_{\mathbb R^n} K(x,y)\chi(y)[f(y)-f(x)]\,dw_{k,a}(y)+f(x)T\chi(x).
    \end{equation}
    For $x\in B_a(x_0,Mr)$,
    \begin{align*}
        |Tf(x)|
        &\leq C_M+\frac C{r^\theta}
        \int_{d(x,y)<C_Mr}
        \frac{d(x,y)^\theta}{w_{k,a}(B_a(x,d(x,y)))}\,dw_{k,a}(y)\\
        &\leq C_M+\frac C{r^\theta}
        \sum_{j=0}^\infty(2^{-j}C_Mr)^\theta
        \leq C_M.
    \end{align*}
    For the approximation, first establish the estimate at a H\"{o}lder exponent strictly smaller than the one in the hypothesis.
   The local integral is continuous
    in that seminorm by \eqref{eq:commutator-domination}, and the cutoff
    term is continuous in the supremum norm. This defines the same
    local distributional extension and proves \eqref{eq:Meyer}.
\end{proof}

Assume in addition that $T1=T^*1=0$.
The local formula also gives the extension needed in the molecule
and reproducing arguments. Fix $0<s<\varepsilon$ and a compact
output test $g$. Choose an orbit cutoff $\chi$ equal to one near
$\supp g$ and all its reflected copies. For $f\in X_s$, put
\[
\langle Tf,g\rangle
=\langle T(\chi f),g\rangle+
\int_{\mathbb R^n}\overline{g(x)}\int_{\mathbb R^n} K(x,y)(1-\chi(y))f(y)
\,dw_{k,a}(y)\,dw_{k,a}(x).
\]
The first term is defined by \eqref{eq:local-T-formula}.
It is continuous when $f$ converges locally uniformly and
in the seminorm $[\cdot]_s$.
For the second term, orbit annuli and reverse doubling give
\begin{align}
    \int_{d(x,y)>r}|K(x,y)|^2\,dw_{k,a}(y)
    &\leq C\sum_{j=0}^\infty
    \frac{w_{k,a}\{y:d(x,y)<2^{j+1}r\}}
    {w_{k,a}(B_a(x,2^jr))^2}\notag\\
    &\leq\frac C{w_{k,a}(B_a(x,r))}
    \sum_{j=0}^\infty2^{-j\min(n,n_{k,a})}<\infty.
    \label{eq:tail-extension}
\end{align}
On a fixed compact set, its right side is uniformly bounded.
By the Cauchy--Schwarz inequality, the tail is bounded by a constant
multiple of \(\|f\|_2\).
Also, comparison with a fixed local average gives
\[
\|f\|_{L^\infty(E,dw_{k,a})}\leq C_E(\|f\|_2+[f]_s)
\]
for compact $E$. The resulting extension is independent of $\chi$,
by the separated representation applied to the difference of two
cutoffs. Thus $T$ and $T^*$ extend continuously from $X_s$ to local
test distributions, without an assumed $L^2_{k,a}(\mathbb{R}^n)$ bound for $T$.

\subsection{Boundedness on smooth molecules}
As in the classical theory, smooth molecules supply the estimates
required for almost orthogonality. We use the two classes introduced
in \cite[Definitions~2.4 and 2.5]{THHLL}, adapted to the coordinates and
measure of the present setting. The first class has decay in the
coordinate distance; the second has decay in the orbit distance.

\begin{definition}\label{def:molecules}
    Let \(0<\beta\leq1\), \(\delta>0\), \(r>0\), and \(x_0\in\mathbb R^n\),
    with \(\widehat{x}\) as in \eqref{eq:phi}.
    A function $f$ is said to be a smooth molecule if there exists a
    constant $C$ such that
    \begin{enumerate}[label=\textup{(\roman*)},ref=\thetheorem(\roman*),leftmargin=2.2em,labelsep=0.55em,itemsep=6pt,topsep=6pt,parsep=0pt]
        \item\label{eq:strong-molecule-size}
        $\displaystyle
        |f(x)|\leq
        \frac{C}{V(x,x_0,r+d(x,x_0))}
        \left(\frac{r}{r+\|\widehat x-\widehat{x_0}\|}\right)^\delta;
        $

        \item\label{eq:strong-molecule-difference}
        $\displaystyle
        \begin{aligned}[t]
            |f(x)-f(x')|
            &\leq C\left(\frac{\|\widehat x-\widehat{x'}\|}{r}\right)^\beta
            \left\{
            \frac{(r/(r+\|\widehat x-\widehat{x_0}\|))^\delta}
            {V(x,x_0,r+d(x,x_0))}+
            \frac{(r/(r+\|\widehat{x'}-\widehat{x_0}\|))^\delta}
            {V(x',x_0,r+d(x',x_0))}
            \right\};
        \end{aligned}
        $

        \item\label{eq:molecule-cancellation}
        $\displaystyle \int_{\mathbb R^n} f\,dw_{k,a}=0$.
    \end{enumerate}
    We write $f\in\mathcal M(\beta,\delta,r,x_0)$ and define
    $\|f\|_{\mathcal M(\beta,\delta,r,x_0)}$ to be the least constant
    in (i) and (ii).
\end{definition}

Replacing the coordinate distances in the decay factors gives the
weak class. The distance in the H\"older factor is unchanged.
\begin{definition}\label{def:weak-molecules}
    Let $0<\beta\leq1$, $\delta>0$, $r>0$, and $x_0\in\mathbb R^n$.
    A function $f$ belongs to $\widetilde{\mathcal M}(\beta,\delta,r,x_0)$
    if there exists a constant $C$ such that
    \begin{enumerate}[label=\textup{(\roman*)},ref=\thetheorem(\roman*),leftmargin=2.2em,labelsep=0.55em,itemsep=6pt,topsep=6pt,parsep=0pt]
        \item\label{eq:weak-molecule-size}
        $\displaystyle
        |f(x)|\leq
        \frac{C}{V(x,x_0,r+d(x,x_0))}
        \left(\frac{r}{r+d(x,x_0)}\right)^\delta;
        $

        \item\label{eq:weak-molecule-difference}
        $\displaystyle
        \begin{aligned}[t]
            |f(x)-f(x')|
            &\leq C\left(\frac{\|\widehat x-\widehat{x'}\|}{r}\right)^\beta
            \left\{
            \frac{(r/(r+d(x,x_0)))^\delta}{V(x,x_0,r+d(x,x_0))}+
            \frac{(r/(r+d(x',x_0)))^\delta}{V(x',x_0,r+d(x',x_0))}
            \right\};
        \end{aligned}
        $

        \item\label{eq:weak-molecule-cancellation}
        $\displaystyle \int_{\mathbb R^n} f\,dw_{k,a}=0$.
    \end{enumerate}
    The norm $\|f\|_{\widetilde{\mathcal M}(\beta,\delta,r,x_0)}$ is
    the least constant in (i) and (ii).
\end{definition}
It follows from the annular estimates that both classes consist
of integrable functions and that
\[
\mathcal M(\beta,\delta,r,x_0)
\subset\widetilde{\mathcal M}(\beta,\delta,r,x_0).
\]
Moreover,
\[
\begin{aligned}
    \|f\|_1&\leq C\|f\|_{\mathcal M},
    \mathrm{ and}    \|f\|_\infty&\leq
    \frac{C\|f\|_{\mathcal M}}{w_{k,a}(B_a(x_0,r))}.
\end{aligned}
\]
Thus strong molecules belong to $L^2_{k,a}(\mathbb{R}^n)$ and to $X_s$ for every
$s<\beta$, via cutoff and smoothing approximation.

\begin{theorem}\label{thm:molecule}
    Suppose $T1=T^*1=0$ and $T\in\WBP$. For
    $0<\beta<\varepsilon$ and $0<\delta'<\delta<\varepsilon$,
    \begin{equation}\label{eq:molecule-map}
        T:\mathcal M(\beta,\delta,r,x_0)
        \longrightarrow\widetilde{\mathcal M}(\beta,\delta',r,x_0)
    \end{equation}
    is bounded uniformly in $r,x_0$.
\end{theorem}
\begin{proof}
    By homogeneity, we may assume
    \eqref{eq:testing-normalization} and $\|f\|_{\mathcal M}\leq1$.
    We first establish the size estimate in
    Definition~\ref{def:weak-molecules}(i), considering separately
    $d(x,x_0)\leq5r$ and $d(x,x_0)>5r$.

    \emph{Case 1: $d(x,x_0)\leq5r$.} Choose $\chi=1$ for
    $d(y,x_0)\leq10r$ and supported where $d(y,x_0)<20r$. The
    commutation formula and the separated tail give
    \begin{align}
        Tf(x)={}&\int_{\mathbb R^n} K(x,y)\chi(y)[f(y)-f(x)]\,dw_{k,a}(y)
        +f(x)T\chi(x)\notag\\
        &+\int_{\mathbb R^n} K(x,y)(1-\chi(y))f(y)\,dw_{k,a}(y)
        =:I_1+I_2+I_3.
        \label{eq:molecule-near-split}
    \end{align}
    On the first integration set, all volumes at scale $r$ are
    comparable to $w_{k,a}(B_a(x_0,r))$. Hence
    \begin{align*}
        |I_1|
        &\leq\frac C{r^\beta w_{k,a}(B_a(x_0,r))}
        \int_{d(x,y)<25r}
        \frac{d(x,y)^\beta}{w_{k,a}(B_a(x,d(x,y)))}\,dw_{k,a}(y)\\
        &\leq\frac C{w_{k,a}(B_a(x_0,r))}.
    \end{align*}
    By Lemma~\ref{lem:cutoff},
    \[
    |I_2|\leq C|f(x)|\leq\frac C{w_{k,a}(B_a(x_0,r))}.
    \]
    If $1-\chi(y)\ne0$, then $d(x,y)\geq5r$ and
    $d(x,y)\asymp d(y,x_0)$. Therefore
    \begin{align*}
        |I_3|
        &\leq\frac C{w_{k,a}(B_a(x_0,r))}
        \int_{\mathbb R^n}\frac{(r/(r+d(y,x_0)))^\delta}
        {V(y,x_0,r+d(y,x_0))}\,dw_{k,a}(y)\leq\frac C{w_{k,a}(B_a(x_0,r))}.
    \end{align*}
    The desired weak size bound follows because
    $r+d(x,x_0)\asymp r$.

    \emph{Case 2: $D=d(x,x_0)>5r$.} Fix disjoint orbit cutoffs $I,J$,
    where $I$ is one near $\mathcal O(x)$ and supported in
    $d(x,y)<D/8$, and $J$ is one near $\mathcal O(x_0)$ and supported
    in $d(x_0,y)<D/8$. They may be chosen equal to one on the
    corresponding balls of radius $D/16$. Put
    \[
    f_1=If,\qquad f_2=Jf,\qquad \mathrm{and} \qquad f_3=(1-I-J)f.
    \]
    On $\supp I$, $d(y,x_0)\asymp D$, and the relevant volumes are
    comparable to $V(x,x_0,D)$. The following estimates will be used:
    \begin{enumerate}[label=\textup{(\roman*)},ref=\thetheorem(\roman*),leftmargin=2.2em,labelsep=0.55em,itemsep=6pt,topsep=6pt,parsep=0pt]
        \item\label{eq:molecule-f1-properties}
        $\displaystyle \|f_1\|_\infty\leq
        \frac{C}{V(x,x_0,D)}(r/D)^\delta$;

        \item\label{cond:f1-L1}
        $\displaystyle \|f_1\|_1\leq C(r/D)^\delta$;

        \item\label{cond:f1-difference}
        $\displaystyle
        |f_1(y)-f_1(y')|
        \leq\frac{C}{V(x,x_0,D)}(r/D)^\delta
        \left(\frac{\|\widehat y-\widehat{y'}\|}{r}\right)^\beta$.
    \end{enumerate}
    The first two estimates follow from the size condition and
    the volume of $\supp I$. To prove (iii), for two points in a fixed
    enlargement of $\supp I$, apply the product rule to $I f$ and
    condition~\ref{eq:strong-molecule-difference}. If only one lies there and
    their distance exceeds a fixed multiple of $D$, use the supremum
    bound and $D>r$. This covers all pairs. Furthermore, on $\supp f_3$
    both orbit distances exceed $D/16$, so annular integration gives
    \begin{equation}\label{eq:molecule-f3-mass}
        \|f_3\|_1\leq C\int_{d(y,x_0)>D/16}
        \frac{(r/(r+d(y,x_0)))^\delta}
        {V(y,x_0,r+d(y,x_0))}\,dw_{k,a}(y)
        \leq C(r/D)^\delta.
    \end{equation}
    Since $\int_{\mathbb R^n} f \,dw_{k,a}=0$,
    \begin{equation}\label{eq:molecule-f2-mass}
        \left|\int_{\mathbb R^n} f_2\,dw_{k,a}\right|
        =\left|\int_{\mathbb R^n}(f_1+f_3)\,dw_{k,a}\right|
        \leq C(r/D)^\delta.
    \end{equation}

    Choose an orbit cutoff $u$ equal to one on $\supp f_1$, supported
    at scale $D$ about $x$. Then
    \[
    Tf_1(x)=\int_{\mathbb R^n} K(x,y)u(y)[f_1(y)-f_1(x)]\,dw_{k,a}(y)
    +f_1(x)Tu(x).
    \]
    The cutoff term is at most $C(r/D)^\delta/V(x,x_0,D)$.
    On $d(x,y)\leq r$, estimate (iii) above  gives
    \begin{align*}
        &\int_{d(x,y)\leq r}|K(x,y)u(y)|\,|f_1(y)-f_1(x)|\,dw_{k,a}(y)\\
        &\quad\leq\frac{C(r/D)^\delta}{r^\beta V(x,x_0,D)}
        \int_{d(x,y)\leq r}\frac{d(x,y)^\beta}
        {w_{k,a}(B_a(x,d(x,y)))}\,dw_{k,a}(y)
        \leq\frac C{V(x,x_0,D)}(r/D)^\delta.
    \end{align*}
    On $r<d(x,y)<CD$, use the supremum bound. With
    $\delta-\delta'>0$,
    \begin{align*}
        &\int_{r<d(x,y)<CD}|K(x,y)u(y)|\,|f_1(y)-f_1(x)|\,dw_{k,a}(y)\\
        &\quad\leq\frac{C(r/D)^\delta}{V(x,x_0,D)}
        \int_{r<d(x,y)<CD}
        \frac{(d(x,y)/r)^{\delta-\delta'}}
        {w_{k,a}(B_a(x,d(x,y)))}\,dw_{k,a}(y)\\
        &\quad\leq\frac{C(r/D)^\delta}{V(x,x_0,D)}(D/r)^{\delta-\delta'}
        =\frac C{V(x,x_0,D)}(r/D)^{\delta'}.
    \end{align*}
    Consequently,
    \begin{equation}\label{eq:molecule-Tf1}
        |Tf_1(x)|\leq\frac C{V(x,x_0,D)}(r/D)^{\delta'}.
    \end{equation}

    For $f_2$, use cancellation with its possibly nonzero mean:
    \begin{align*}
        Tf_2(x)
        ={}&\int_{\mathbb R^n}[K(x,y)-K(x,x_0)]f_2(y)\,dw_{k,a}(y)
        +K(x,x_0)\int_{\mathbb R^n} f_2\,dw_{k,a}.
    \end{align*}
    The second term is bounded by
    $C(r/D)^\delta/V(x,x_0,D)$ using \eqref{eq:molecule-f2-mass}.
    Split the first integral according to
    $\|\widehat y-\widehat{x_0}\|\leq D/4$ or greater than $D/4$.
    For the first part, the kernel difference condition gives
    \begin{align*}
        &\int_{\|\widehat y-\widehat{x_0}\|\leq D/4}
        |K(x,y)-K(x,x_0)|\,|f_2(y)|\,dw_{k,a}(y)\\
        &\quad\leq\frac C{V(x,x_0,D)D^{\delta'}}
        \int_{\mathbb R^n}\frac{\|\widehat y-\widehat{x_0}\|^{\delta'}}
        {V(y,x_0,r+d(y,x_0))}
        \left(\frac r{r+\|\widehat y-\widehat{x_0}\|}\right)^\delta
        dw_{k,a}(y)\\
        &\quad\leq\frac{C(r/D)^{\delta'}}{V(x,x_0,D)}
        \int_{\mathbb R^n}\frac{(r/(r+d(y,x_0)))^{\delta-\delta'}}
        {V(y,x_0,r+d(y,x_0))}\,dw_{k,a}(y)\\
        &\quad\leq\frac C{V(x,x_0,D)}(r/D)^{\delta'}.
    \end{align*}
    On the other part, $d(x_0,y)<D/8$ but
    $\|\widehat y-\widehat{x_0}\|>D/4$. Both $K(x,y)$ and
    $K(x,x_0)$ are bounded by $C/V(x,x_0,D)$. Extracting
    $(r/D)^{\delta'}$ from the ordinary-distance decay of $f$ gives
    \begin{align*}
        &\int_{\substack{d(y,x_0)<D/8\\
                \|\widehat y-\widehat{x_0}\|>D/4}}
        |K(x,y)-K(x,x_0)|\,|f_2(y)|\,dw_{k,a}(y)\\
        &\quad\leq\frac{C(r/D)^{\delta'}}{V(x,x_0,D)}
        \int_{\mathbb R^n}\frac{(r/(r+d(y,x_0)))^{\delta-\delta'}}
        {V(y,x_0,r+d(y,x_0))}\,dw_{k,a}(y)\\
        &\quad\leq\frac C{V(x,x_0,D)}(r/D)^{\delta'}.
    \end{align*}
    Thus
    \begin{equation}\label{eq:molecule-Tf2}
        |Tf_2(x)|\leq\frac C{V(x,x_0,D)}(r/D)^{\delta'}.
    \end{equation}
    Finally, $d(x,y)>D/16$ on $\supp f_3$, and
    \begin{equation}\label{eq:molecule-Tf3}
        |Tf_3(x)|\leq\frac C{V(x,x_0,D)}\|f_3\|_1
        \leq\frac C{V(x,x_0,D)}(r/D)^\delta.
    \end{equation}
    Combining \eqref{eq:molecule-Tf1}--\eqref{eq:molecule-Tf3}
    with the near-orbit case proves
    \begin{equation}\label{eq:weak-molecule-image-size}
        |Tf(x)|\leq\frac C{V(x,x_0,r+d(x,x_0))}
        \left(\frac r{r+d(x,x_0)}\right)^{\delta'}.
    \end{equation}

    We next prove the regularity condition in
    Definition~\ref{def:weak-molecules}(ii). Write
    $h=\|\widehat x-\widehat{x'}\|$. If $h\geq r/32$,
    \eqref{eq:weak-molecule-image-size} at the two points gives the
    required difference estimate, since $(h/r)^\beta$ is bounded
    below. It remains to consider $0<h<r/32$.

    \emph{Case 1: $d(x,x_0)>10r$.} Set $D=d(x,x_0)$. Use a fixed
    cutoff at scale $D$ about the orbit of $x$ and write
    $f=f_1+f_2$, where $f_1$ is supported in $d(x,y)<D/4$ and
    $f_2$ vanishes for $d(x,y)<D/8$. The three estimates for $f_1$ in Case 2 hold for this cutoff as well. Choose a small orbit cutoff
    $\chi$ at scale $h$ which is one in neighborhoods of both orbits
    of $x,x'$ and for which $d(x,y)>ch$, $c>2$, on $\supp(1-\chi)$.
    Set
    \[
    p(x)=\int_{\mathbb R^n} K(x,y)\chi(y)[f_1(y)-f_1(x)]\,dw_{k,a}(y)
    \]
    and define $p(x')$ with the same cutoff. Then
    \begin{align}
        |p(x)|+|p(x')|
        &\leq\frac{C(r/D)^\delta}{r^\beta V(x,x_0,D)}
        \int_{d(x,y)<Ch}
        \frac{d(x,y)^\beta}{w_{k,a}(B_a(x,d(x,y)))}\,dw_{k,a}(y)\notag\\
        &\leq\frac C{V(x,x_0,D)}(h/r)^\beta(r/D)^\delta.
        \label{eq:molecule-small-cutoff}
    \end{align}
    The assumption $T1=0$ implies
    \[
    T\chi(x)-T\chi(x')
    =-\int_{\mathbb R^n}[K(x,y)-K(x',y)](1-\chi(y))\,dw_{k,a}(y).
    \]
    Subtracting the two local representations therefore gives
    \begin{align}
        &[Tf_1(x)-p(x)]-[Tf_1(x')-p(x')]\notag\\
        &\quad=\int_{\mathbb R^n}[K(x,y)-K(x',y)](1-\chi(y))
        [f_1(y)-f_1(x)]\,dw_{k,a}(y)\notag\\
        &\qquad+[f_1(x)-f_1(x')]T\chi(x').
        \label{eq:molecule-difference-identity}
    \end{align}
    The integral term is at most
    \begin{align*}
        &\frac{Ch^\varepsilon(r/D)^\delta}{r^\beta V(x,x_0,D)}
        \int_{d(x,y)>ch}
        \frac{\|\widehat x-\widehat y\|^{\beta-\varepsilon}}
        {w_{k,a}(B_a(x,d(x,y)))}\,dw_{k,a}(y)\\
        &\quad\leq\frac{Ch^\varepsilon(r/D)^\delta}{r^\beta V(x,x_0,D)}
        \sum_{j=0}^\infty(2^jch)^{\beta-\varepsilon}
        \leq\frac C{V(x,x_0,D)}(h/r)^\beta(r/D)^\delta.
    \end{align*}
    We used $\beta<\varepsilon$ and
    $\|\widehat x-\widehat y\|\geq d(x,y)$. The last term of
    \eqref{eq:molecule-difference-identity} has the same bound by
    Lemma~\ref{lem:cutoff}. Together with
    \eqref{eq:molecule-small-cutoff}, this proves the regularity
    estimate for $Tf_1$.

    For $f_2$, the separated kernel estimate and $\|f_2\|_1\leq C$
    give
    \begin{align*}
        |Tf_2(x)-Tf_2(x')|
        &\leq\int_{\mathbb R^n}|K(x,y)-K(x',y)|\,|f_2(y)|\,dw_{k,a}(y)\\
        &\leq\frac C{V(x,x_0,D)}(h/D)^\varepsilon\|f_2\|_1\\
        &\leq\frac C{V(x,x_0,D)}(h/r)^\beta(r/D)^{\delta'}.
    \end{align*}
    Indeed, $h<r<D$ and $\beta,\delta'<\varepsilon$.

    \emph{Case 2: $d(x,x_0)\leq10r$.} Use the same small cutoff $\chi$ and
    \eqref{eq:molecule-difference-identity} with $f$ in place of $f_1$.
    The global bounds
    \[
    |f(y)|\leq\frac C{w_{k,a}(B_a(x_0,r))} \qquad \mathrm{and } \qquad
    |f(y)-f(x)|\leq\frac C{w_{k,a}(B_a(x_0,r))}
    \left(\frac{\|\widehat y-\widehat x\|}{r}\right)^\beta
    \]
    follow directly from Definition~\ref{def:molecules}. Thus

    $$  |p(x)|+|p(x')|
    \leq\frac{C(h/r)^\beta}{w_{k,a}(B_a(x_0,r))}$$

    and
    \begin{align*}
        \left|\int_{\mathbb R^n}[K(x,y)-K(x',y)](1-\chi(y))[f(y)-f(x)]\,dw_{k,a}(y)\right|
        &\leq\frac{Ch^\varepsilon}{r^\beta w_{k,a}(B_a(x_0,r))}
        \sum_{j=0}^\infty(2^jch)^{\beta-\varepsilon}\\
        &\leq\frac{C(h/r)^\beta}{w_{k,a}(B_a(x_0,r))}.
    \end{align*}
    The cutoff term satisfies the same estimate. Volumes at $x$, $x'$ and
    $x_0$ are comparable at scale $r$. Both cases prove

    \begin{equation}\label{eq:weak-molecule-image-difference}
        \begin{aligned}
            |Tf(x)-Tf(x')|\leq{}&C
            \left(\frac{\|\widehat x-\widehat{x'}\|}{r}\right)^\beta
            \left\{\frac{(r/(r+d(x,x_0)))^{\delta'}}
            {V(x,x_0,r+d(x,x_0))}+
            \frac{(r/(r+d(x',x_0)))^{\delta'}}
            {V(x',x_0,r+d(x',x_0))}\right\}.
        \end{aligned}
    \end{equation}

It remains to prove the cancellation condition. The size
estimate shows that $Tf$ is integrable. For a compactly supported mean-zero
smooth input $f$, take expanding orbit cutoffs $\chi_L$. The
input cancellation and kernel difference condition show that the
far tail of $Tf$ has integral tending to zero. The definition of
$T^*1$ consequently gives
\[
\int_{\mathbb R^n} Tf\,dw_{k,a}
=\lim_{L\to\infty}\langle Tf,\chi_L\rangle
=\langle f,T^*1\rangle=0.
\]
For general $f\in\mathcal M(\beta,\delta,r,x_0)$, fix
$0<\beta_0<\beta$ and $0<\delta_0<\delta$, and let $\chi_L$ be
an ordinary-coordinate cutoff of radius $Lr$ centered at $x_0$.
Choose a fixed smooth $\psi$ supported in $B_a(x_0,r)$ with
$\int_{\mathbb R^n}\psi\,dw_{k,a}=1$, and set
\[
f_L=\chi_L f-\left(\int_{\mathbb R^n}\chi_Lf\,dw_{k,a}\right)\psi.
\]
By extracting the factor $(r/(r+\|\widehat x-\widehat{x_0}\|))^{\delta-\delta_0}$
outside the cutoff,
\[
\int_{\mathbb R^n}|(1-\chi_L)f|\,dw_{k,a}
\leq C L^{-(\delta-\delta_0)}.
\]
For pairs separated by at most $r$, the same extraction in the
product difference, together with
$\|\nabla_{\widehat x}\chi_L\|_\infty\leq C/(Lr)$, controls the
molecule seminorm. For pairs farther apart, the size bound controls
the difference. After smoothing the compact inputs at an arbitrarily
small scale, one obtains mean-zero smooth $f_{L,j}$ with
\[
\lim_{L\to\infty}\lim_{j\to\infty}
\|f-f_{L,j}\|_{\mathcal M(\beta_0,\delta_0,r,x_0)}=0.
\]
The preceding size and regularity calculations applied with smaller
exponents show convergence of the images in a weak-molecule norm,
and hence in $L^1_{k,a}(\mathbb{R}^n)$. Choose $0<s<\min\{\beta_0,\varepsilon\}$.
Convergence in the strong-molecule norm also implies convergence
in $X_s$. By the continuity of the local extension of $T$ on $X_s$,
the limit agrees with that extension.
Passing to the integral gives $\int_{\mathbb R^n} Tf=0$, which finishes
\eqref{eq:molecule-map}.

\end{proof}

\subsection{\texorpdfstring{Proof of the $T1$ theorem}{Proof of the T1 theorem}}\label{sec:T1-proof}
The necessity has been proved in Proposition~\ref{prop:necessity}.
For sufficiency, we first treat the case $T1=T^*1=0$ and then
recover the general case by means of paraproducts. We begin with Coifman's
approximation to the identity, which is employed in the almost
orthogonal estimates and in the reproducing formula.
Choose a nonnegative smooth radial function $\vartheta$ equal to one
on $[0,1]$ and vanishing on $[2,\infty)$.
 For $r_j=2^{-j}$, define
\begin{equation}\label{eq:Coifman-construction}
    \begin{gathered}
        L_jf(x)=\int_{\mathbb R^n}\vartheta(\|\widehat x-\widehat y\|/r_j)f(y)\,dw_{k,a}(y),
        \qquad b_j=(L_j1)^{-1},\\
        S_jf=b_jL_j\left(\frac{L_j(b_jf)}{L_jb_j}\right),
        \qquad D_j=S_j-S_{j-1}.
    \end{gathered}
\end{equation}
The functions $b_j$ serve only as normalizing factors. This is the
Coifman construction in \cite{Hd,HMY,THHLL}, now applied to the geometry
of Proposition~\ref{prop:geometry}.

\begin{lemma}\label{lem:approx}
  Here and below, for every \(z\in\mathbb R^n\),
  \(\widehat z\) is as in \eqref{eq:phi}.  For every $j\in\mathbb Z$, the kernels $S_j$ and $D_j$ have the
    following properties:
    \begin{enumerate}[label=\textup{(\roman*)},ref=\thetheorem(\roman*),leftmargin=2.2em,labelsep=0.55em,itemsep=6pt,topsep=6pt,parsep=0pt]
        \item\label{cond:S-symmetry}
        $\displaystyle S_j(x,y)=S_j(y,x)\geq0$, and $S_j$ is real-valued;

        \item\label{cond:SD-support}
        $\displaystyle S_j(x,y)=D_j(x,y)=0$
        whenever $\|\widehat x-\widehat y\|>8r_j$;

        \item\label{eq:Dj-estimates}
        $\displaystyle
        \begin{aligned}[t]
            &|\partial_{\widehat x}^{\alpha}\partial_{\widehat y}^{\beta}S_j(x,y)|
            +|\partial_{\widehat x}^{\alpha}\partial_{\widehat y}^{\beta}D_j(x,y)|\leq
            \frac{Cr_j^{-|\alpha|-|\beta|}}
            {w_{k,a}(B_a(x,r_j))+w_{k,a}(B_a(y,r_j))},
            \qquad |\alpha|,|\beta|\leq1;
        \end{aligned}
        $

        \item\label{eq:S-normalization}
        $\displaystyle \int_{\mathbb R^n} S_j(x,y)\,dw_{k,a}(x)=1$ for every $y\in\mathbb R^n$;

        \item\label{cond:S-mass-y}
        $\displaystyle \int_{\mathbb R^n} S_j(x,y)\,dw_{k,a}(y)=1$ for every $x\in\mathbb R^n$;

        \item\label{cond:D-cancel-x}
        $\displaystyle \int_{\mathbb R^n} D_j(x,y)\,dw_{k,a}(x)=0$ for every $y\in\mathbb R^n$;

        \item\label{cond:D-cancel-y}
        $\displaystyle \int_{\mathbb R^n} D_j(x,y)\,dw_{k,a}(y)=0$ for every $x\in\mathbb R^n$;

        \item\label{eq:unreflected-AO}
        $\displaystyle \|D_jD_\ell\|_{2\to2}\leq C2^{-|j-\ell|}$;

        \item\label{cond:D-square-function}
        $\displaystyle \sum_j\|D_jf\|_2^2\leq C\|f\|_2^2$.
    \end{enumerate}
\end{lemma}
\begin{proof}
  The   doubling condition gives
    \[
    L_j1(x)\asymp w_{k,a}(B_a(x,r_j)),\qquad
    L_jb_j(x)\asymp1,\qquad \mathrm{and} \qquad
    |\partial_{\widehat x}^{\alpha}b_j(x)|
    \leq Cr_j^{-|\alpha|}b_j(x)\quad(|\alpha|\leq1).
    \]
    The kernel is
    \begin{equation}\label{eq:S-kernel}
        S_j(x,y)=b_j(x)b_j(y)\int_{\mathbb R^n}
        \frac{\vartheta(\|\widehat x-\widehat u\|/r_j)
            \vartheta(\|\widehat u-\widehat y\|/r_j)}{L_jb_j(u)}
        \,dw_{k,a}(u).
    \end{equation}
    Differentiating the smooth factors proves (iii), without
    differentiating the measure. Symmetry and the support condition
    follow from the same formula. Moreover,
    \[
    S_j1=b_jL_j\left(\frac{L_jb_j}{L_jb_j}\right)=b_jL_j1=1.
    \]
    Symmetry gives the other mass identity, and subtraction at
    consecutive scales gives the two cancellations of $D_j$. This
    proves (iv)--(vii).

Suppose $r_j\leq r_\ell$. The zero-mean property of $D_j$ gives
    \begin{align*}
        D_jD_\ell(x,y)
        &=\int_{\mathbb R^n} D_j(x,u)[D_\ell(u,y)-D_\ell(x,y)]\,dw_{k,a}(u) \\ \intertext{and}
        |D_jD_\ell(x,y)|
        &\leq\frac{C(r_j/r_\ell)}
        {w_{k,a}(B_a(x,r_\ell))+w_{k,a}(B_a(y,r_\ell))}
        \mathbf1_{\{\|\widehat x-\widehat y\|\leq Cr_\ell\}}.
    \end{align*}
    The two Schur integrals are bounded above by $Cr_j/r_\ell$. Transposition
    handles the case $r_\ell\leq r_j$.
     Hence
    \[
    \sup_j\sum_\ell\|D_jD_\ell\|_{2\to2}^{1/2}<\infty.
    \]
The Cotlar--Stein lemma yields bounds for every finite signed sum $\sum_j\epsilon_jD_j$
on $L^2_{k,a}(\mathbb{R}^n)$, uniformly for $\epsilon_j=\pm1$.
 For a finite set \(J\subset\mathbb Z\), averaging over all choices of signs gives
\[
\sum_{j\in J}\|D_jf\|_2^2
=
\frac{1}{2^{|J|}}
\sum_{\epsilon\in\{-1,1\}^{J}}
\left\|\sum_{j\in J}\epsilon_jD_jf\right\|_2^2
\leq C\|f\|_2^2,
\]
where \(C\) is independent of \(J\). Taking
\(J=\{-N,\ldots,N\}\) and letting \(N\to\infty\)
then proves the desired estimate.
\end{proof}

To estimate $D_jTD_\ell$, we need the following convolution and
molecule estimates. We state the convolution bound first and then
use cancellation to obtain the scale decay according to
\cite[Lemmas~2.9--2.11]{THHLL}.

\begin{lemma}\label{lem:volume-convolution}
    For $S\geq s>0$ and $\alpha,\beta>0$,
    \begin{equation}\label{eq:volume-convolution}
        \begin{aligned}
            &\int_{\mathbb R^n}\frac{(S/(S+d(x,u)))^\alpha}{V(x,u,S+d(x,u))}
            \frac{(s/(s+d(u,y)))^\beta}{V(u,y,s+d(u,y))}\,dw_{k,a}(u)\leq\frac C{V(x,y,S+d(x,y))}.
        \end{aligned}
    \end{equation}
\end{lemma}
\begin{proof}
    Write $D=d(x,y)$. If $D\leq S$, the first volume denominator
    is at least $w_{k,a}(B_a(x,S))$. Thus the left side is at most
    \[
    \frac1{w_{k,a}(B_a(x,S))}
    \int_{\mathbb R^n}\frac{(s/(s+d(u,y)))^\beta}{V(u,y,s+d(u,y))}\,dw_{k,a}(u)
    \leq\frac C{V(x,y,S+D)}.
    \]
    If $D>S$, every $u$ belongs to at least one of the sets
    $d(x,u)\geq D/2$ and $d(u,y)\geq D/2$. On the first set,
    \begin{align*}
        &\int_{d(x,u)\geq D/2}
        \frac{(S/(S+d(x,u)))^\alpha}{V(x,u,S+d(x,u))}
        \frac{(s/(s+d(u,y)))^\beta}{V(u,y,s+d(u,y))}\,dw_{k,a}(u)\\
        &\quad\leq\frac C{w_{k,a}(B_a(x,S+D))}
        \int_{\mathbb R^n}\frac{(s/(s+d(u,y)))^\beta}{V(u,y,s+d(u,y))}\,dw_{k,a}(u)\\
        &\quad\leq\frac C{V(x,y,S+D)}.
    \end{align*}
    On the second set,
    \begin{align*}
        &\int_{d(u,y)\geq D/2}
        \frac{(S/(S+d(x,u)))^\alpha}{V(x,u,S+d(x,u))}
        \frac{(s/(s+d(u,y)))^\beta}{V(u,y,s+d(u,y))}\,dw_{k,a}(u)\\
        &\quad\leq\frac C{w_{k,a}(B_a(y,D))}
        \int_{\mathbb R^n}\frac{(S/(S+d(x,u)))^\alpha}{V(x,u,S+d(x,u))}\,dw_{k,a}(u)\\
        &\quad\leq\frac C{V(x,y,S+D)}.
    \end{align*}
    The last comparisons use $D>S\geq s$ and reflection invariance.
    This proves \eqref{eq:volume-convolution}.
\end{proof}

\begin{lemma}\label{lem:molecule-AO}
    Suppose $S\geq s>0$,
    $F_S(x,\cdot)\in\widetilde{\mathcal M}(\varepsilon_0,\varepsilon_0,S,x)$
    and $G_s(\cdot,y)\in\mathcal M(\varepsilon_0,\varepsilon_0,s,y)$,
    with norms at most one. For $0<\varepsilon_1,\varepsilon_2<\varepsilon_0$,
    \begin{equation}\label{eq:weak-strong-AO}
        \left|\int_{\mathbb R^n} F_S(x,u)G_s(u,y)\,dw_{k,a}(u)\right|
        \leq C(s/S)^{\varepsilon_1}
        \frac{(S/(S+d(x,y)))^{\varepsilon_2}}{V(x,y,S+d(x,y))}.
    \end{equation}
    For the reverse order of scales, transpose the statement.
\end{lemma}
\begin{proof}
    It suffices to prove \eqref{eq:weak-strong-AO} with
    $\varepsilon=\max\{\varepsilon_1,\varepsilon_2\}<\varepsilon_0$.
   Since \(\int_{\mathbb R^n}G_s(u,y)\,dw_{k,a}(u)=0\), the constant
   \(F_S(x,y)\) may be subtracted inside the integral:
    \begin{align}
        \int_{\mathbb R^n} F_S(x,u)G_s(u,y)\,dw_{k,a}(u)
        &=\int_{\mathbb R^n}[F_S(x,u)-F_S(x,y)]G_s(u,y)\,dw_{k,a}(u),\notag\\
        \left|\int_{\mathbb R^n} F_SG_s\,dw_{k,a}\right|
        &\leq\int_{\|\widehat u-\widehat y\|\leq S}
        |F_S(x,u)-F_S(x,y)|\,|G_s(u,y)|\,dw_{k,a}(u)\notag\\
        &\quad+\int_{\|\widehat u-\widehat y\|>S}
        [|F_S(x,u)|+|F_S(x,y)|]|G_s(u,y)|\,dw_{k,a}(u)\notag\\
        &=: I+II.\label{eq:AO-two-regions}
    \end{align}
    For $h=\|\widehat u-\widehat y\|\leq S$,
    \begin{equation}\label{eq:AO-scale-extraction}
        (h/S)^{\varepsilon_0}(s/(s+h))^{\varepsilon_0}
        \leq(s/S)^\varepsilon(s/(s+h))^{\varepsilon_0-\varepsilon}.
    \end{equation}
    The difference condition for $F_S$ therefore gives
    \begin{align}
        I\leq{}&C(s/S)^\varepsilon
        \int_{h\leq S}\left\{
        \frac{(S/(S+d(x,u)))^{\varepsilon_0}}{V(x,u,S+d(x,u))}
        +\frac{(S/(S+d(x,y)))^{\varepsilon_0}}{V(x,y,S+d(x,y))}
        \right\}
        \frac{(s/(s+h))^{\varepsilon_0-\varepsilon}}{V(u,y,s+d(u,y))}
        \,dw_{k,a}(u)\notag\\
        =: {}&C(s/S)^\varepsilon(I_1+I_2).
        \label{eq:AO-I}
    \end{align}
    For $I_2$, replace $h$ by $d(u,y)$ in the numerator and use
    \eqref{eq:annuli} to obtain
    \[
    I_2\leq
    \frac{C(S/(S+d(x,y)))^{\varepsilon_0}}{V(x,y,S+d(x,y))}
    \leq\frac{C(S/(S+d(x,y)))^\varepsilon}{V(x,y,S+d(x,y))}.
    \]
    If $d(x,y)\leq2S$, Lemma~\ref{lem:volume-convolution} gives the
    same bound for $I_1$. If $d(x,y)>2S$, then $h\leq S$ implies
    $d(x,u)\geq d(x,y)/2$. Extracting an exponent $\varepsilon$ from
    the first numerator and using the convolution lemma with exponents
    $\varepsilon_0-\varepsilon$ gives
    \begin{align*}
        I_1
        &\leq C\left(\frac S{S+d(x,y)}\right)^\varepsilon
        \int_{\mathbb R^n}\frac{(S/(S+d(x,u)))^{\varepsilon_0-\varepsilon}}
        {V(x,u,S+d(x,u))}
        \frac{(s/(s+d(u,y)))^{\varepsilon_0-\varepsilon}}
        {V(u,y,s+d(u,y))}\,dw_{k,a}(u)\\
        &\leq\frac{C(S/(S+d(x,y)))^\varepsilon}{V(x,y,S+d(x,y))}.
    \end{align*}
    Thus
    \begin{equation}\label{eq:AO-I-bound}
        I\leq C(s/S)^\varepsilon
        \frac{(S/(S+d(x,y)))^\varepsilon}{V(x,y,S+d(x,y))}.
    \end{equation}

    For $II$, write $II=II_1+II_2$ according to the two terms in
    square brackets. The term containing $F_S(x,y)$ satisfies
    \begin{align*}
        II_2
        &\leq\frac{(S/(S+d(x,y)))^{\varepsilon_0}}{V(x,y,S+d(x,y))}
        \int_{h>S}\frac{(s/(s+h))^{\varepsilon_0}}
        {V(u,y,s+d(u,y))}\,dw_{k,a}(u)\\
        &\leq C(s/S)^\varepsilon
        \frac{(S/(S+d(x,y)))^\varepsilon}{V(x,y,S+d(x,y))}
        \int_{\mathbb R^n}\frac{(s/(s+d(u,y)))^{\varepsilon_0-\varepsilon}}
        {V(u,y,s+d(u,y))}\,dw_{k,a}(u)\\
        &\leq C(s/S)^\varepsilon
        \frac{(S/(S+d(x,y)))^\varepsilon}{V(x,y,S+d(x,y))}.
    \end{align*}
    For $II_1$, split into
    $d(x,u)\geq d(x,y)/2$ and $d(u,y)\geq d(x,y)/2$; their union
    covers the integration domain. On the first part extract
    $(S/(S+d(x,y)))^\varepsilon$ from $F_S$ and
    $(s/S)^\varepsilon$ from $G_s$, using $h>S$. The remaining
    integral is bounded by Lemma~\ref{lem:volume-convolution}, so
    \begin{align*}
        II_{11}
        &\leq C(s/S)^\varepsilon
        (S/(S+d(x,y)))^\varepsilon\\
        &\quad\times\int_{\mathbb R^n}
        \frac{(S/(S+d(x,u)))^{\varepsilon_0-\varepsilon}}
        {V(x,u,S+d(x,u))}
        \frac{(s/(s+d(u,y)))^{\varepsilon_0-\varepsilon}}
        {V(u,y,s+d(u,y))}\,dw_{k,a}(u)\\
        &\leq C(s/S)^\varepsilon
        \frac{(S/(S+d(x,y)))^\varepsilon}{V(x,y,S+d(x,y))}.
    \end{align*}
    On the second part, $h>S$ and $h\geq d(u,y)\geq d(x,y)/2$.
    Consequently,
    \begin{equation}\label{eq:AO-II12-factor}
        \left(\frac s{s+h}\right)^\varepsilon
        \leq C\left(\frac s{S+d(x,y)}\right)^\varepsilon
        =C(s/S)^\varepsilon
        \left(\frac S{S+d(x,y)}\right)^\varepsilon.
    \end{equation}
    Applying the same convolution estimate to the remaining factors,
    \begin{align*}
        II_{12}
        &\leq C(s/S)^\varepsilon
        (S/(S+d(x,y)))^\varepsilon\int_{\mathbb R^n}
        \frac{(S/(S+d(x,u)))^{\varepsilon_0}}{V(x,u,S+d(x,u))}
        \frac{(s/(s+d(u,y)))^{\varepsilon_0-\varepsilon}}
        {V(u,y,s+d(u,y))}\,dw_{k,a}(u)\\
        &\leq C(s/S)^\varepsilon
        \frac{(S/(S+d(x,y)))^\varepsilon}{V(x,y,S+d(x,y))}.
    \end{align*}
    This proves the assertion, since lowering the two exponents
    from $\varepsilon$ to $\varepsilon_1,\varepsilon_2$ increases its
    right side.

    Transposition handles the opposite scale order.
\end{proof}

\begin{lemma}\label{lem:T-AO}
    Suppose $T1=T^*1=0$ and \eqref{eq:testing-normalization} holds.
    For any $0<\theta,\delta<\varepsilon$, and
    $r=\max\{2^{-j},2^{-\ell}\}$,
    \begin{equation}\label{eq:T-AO-kernel}
        |D_jTD_\ell(x,y)|
        \leq\frac{C2^{-|j-\ell|\theta}}{V(x,y,r+d(x,y))}
        \left(\frac r{r+d(x,y)}\right)^\delta.
    \end{equation}
    Consequently,
    \begin{equation}\label{eq:T-AO-L2}
        \|D_jTD_\ell\|_{2\to2}\leq C2^{-|j-\ell|\theta}.
    \end{equation}
    These are estimates for blocks of the initially defined distribution;
    no $L^2_{k,a}(\mathbb{R}^n)$ bound for $T$ is used.
\end{lemma}
\begin{proof}
    Choose $\max\{\theta,\delta\}<\varepsilon_0<\delta_0<\varepsilon$.
    Suppose first $j\leq\ell$,  and put $S=2^{-j}$, $s=2^{-\ell}$.
    The support, derivative and cancellation estimates of
    Lemma~\ref{lem:approx} show that $D_j(x,\cdot)$ is, up to a
    uniform constant, a strong molecule in
    $\mathcal M(\varepsilon_0,\delta_0,S,x)$. By
    Theorem~\ref{thm:molecule},
    \[
    F_x(v)=\overline{T^*(D_j(x,\cdot))(v)}
    \in\widetilde{\mathcal M}(\varepsilon_0,\varepsilon_0,S,x)
    \]
    with uniform norm. Both $D_j$ and $D_\ell$ are real and symmetric,
    so the composition kernel is
    \begin{align*}
        D_jTD_\ell(x,y)
        &=\int_{\mathbb R^n} F_x(v)D_\ell(v,y)\,dw_{k,a}(v) \\  \end{align*}
     and \begin{align*}  |D_jTD_\ell(x,y)|
        &\leq C(s/S)^\theta
        \frac{(S/(S+d(x,y)))^\delta}{V(x,y,S+d(x,y))}
    \end{align*}
    by Lemma~\ref{lem:molecule-AO}.
    If $\ell\leq j$, apply the same argument to
    $(D_jTD_\ell)^*=D_\ell T^*D_j$. Its kernel is the conjugate
    transpose, and $V,d$ are symmetric. This gives
    \eqref{eq:T-AO-kernel} in the other order.

    Finally, by \eqref{eq:annuli} in both variables,
    \begin{align*}
        \sup_x\int_{\mathbb R^n}|D_jTD_\ell(x,y)|\,dw_{k,a}(y)
        &\leq C2^{-|j-\ell|\theta}\end{align*}
     and \begin{align*}    \sup_y\int_{\mathbb R^n}|D_jTD_\ell(x,y)|\,dw_{k,a}(x)
        &\leq C2^{-|j-\ell|\theta}.
    \end{align*}
    Schur's test proves \eqref{eq:T-AO-L2}.
\end{proof}

Fix $0<s<\min\{\varepsilon,1/2\}$ and recall the space $X_s$
defined in Section~\ref{sec:analytic-symbols}. We establish the
reproducing formula with convergence in both $L^2_{k,a}(\mathbb R^n)$
and $X_s$.

\begin{lemma}\label{lem:reproduction}
    For a sufficiently large integer $M$, put
    \begin{equation}\label{eq:UM}
        D_j^M=\sum_{|h|\leq M}D_{j+h} \qquad \mathrm{and} \qquad
        U_M=\sum_jD_jD_j^M.
    \end{equation}
    Then the following assertions hold:
    \begin{enumerate}[label=\textup{(\roman*)},ref=\thetheorem(\roman*),leftmargin=2.2em,labelsep=0.55em,itemsep=6pt,topsep=6pt,parsep=0pt]
        \item The sum defining $U_M$ converges strongly on
        $L^2_{k,a}(\mathbb R^n)$ and on $X_s$, and its finite interval
        sums are uniformly bounded on these spaces;

        \item $U_M$ is invertible on both spaces;

        \item\label{eq:reproduction}
        $\displaystyle f=\sum_jD_jD_j^MU_M^{-1}f$,
        with convergence in the corresponding space;

        \item\label{cond:reproduction-square}
        $\displaystyle \sum_j\|D_j^MU_M^{-1}f\|_2^2\leq C_M\|f\|_2^2$.
    \end{enumerate}
\end{lemma}
\begin{proof}
    We first verify the two endpoints of $\sum_jD_j$.
    For a compact smooth function $f$ and $r_j\leq1$,
    \begin{align*}
        \|S_jf-f\|_\infty&\leq Cr_j\|\nabla_{\widehat x}f\|_\infty,&
        [S_jf-f]_1&\leq C\|\nabla_{\widehat x}f\|_\infty,\\
       \end{align*} and \begin{align*}[S_jf-f]_s&\leq C r_j^{1-s}\|\nabla_{\widehat x}f\|_\infty.
    \end{align*}
    Its support lies in a fixed compact enlargement, so $L^2_{k,a}(\mathbb{R}^n)$ convergence
    also holds. At coarse scales, $S_jf(x)\ne0$ implies
    $\|\widehat x\|\leq Cr_j+C_f$. By homogeneity,
    \begin{align*}
        \|\partial_{\widehat x}^{\alpha}S_jf\|_\infty
        &\leq C_fr_j^{-n_{k,a}-|\alpha|}\|f\|_1,\qquad
        |\alpha|\leq1,
        \|S_jf\|_2^2&\leq\|S_jf\|_\infty\|S_jf\|_1
        \leq C_fr_j^{-n_{k,a}}\|f\|_1^2,\\ \end{align*}
     and \begin{align*}
        [S_jf]_s&\leq C_fr_j^{-n_{k,a}-s}\|f\|_1\longrightarrow0.
    \end{align*}
    Normalization and the kernel bounds give
    $\|S_jf\|_2+[S_jf]_s\leq C(\|f\|_2+[f]_s)$.
    Density therefore implies
    \begin{equation}\label{eq:S-limits}
        S_j\longrightarrow I\ (j\to+\infty),\qquad
        S_j\longrightarrow0\ (j\to-\infty),\qquad \mathrm{and} \qquad
        I=\sum_jD_j
    \end{equation}
    strongly on both spaces.

    On $L^2_{k,a}(\mathbb{R}^n)$, the remainder is
    \begin{equation}\label{eq:RM}
        I-U_M=R_M=\sum_{|j-\ell|>M}D_jD_\ell.
    \end{equation}
    For a fixed shift $h$, set $A_j=D_jD_{j+h}$. By
    Lemma~\ref{lem:approx}(viii),
    \begin{align*}
        \|A_j\|&\leq C2^{-|h|}
        \end{align*}
     and \begin{align*} \|A_j^*A_\ell\|+\|A_jA_\ell^*\|
        &\leq C\min\{2^{-2|h|},2^{-|j-\ell|}\}
        \leq C2^{-|h|}2^{-|j-\ell|/2}.
    \end{align*}
The Cotlar--Stein lemma, uniformly for finite sums, yields

    \begin{equation}\label{eq:RM-L2}
        \left\|\sum_jD_jD_{j+h}\right\|_{2\to2}\leq C2^{-|h|/2}
        \qquad \mathrm{and} \qquad \|R_M\|_{2\to2}\leq C2^{-M/2}.
    \end{equation}
  For each fixed shift $h$, the corresponding sums converge strongly
  by the square-function analysis and synthesis bounds.
  Decompose rectangular partial sums by $h=\ell-j$.
  The uniform bound $C2^{-|h|/2}$ is summable in $h$, so passage
  to the limit is justified and yields \eqref{eq:RM}.

    For the H\"older seminorm, telescope the inner sums:
    \begin{equation}\label{eq:RM-Holder}
        R_Mf=\sum_jD_j(I-S_{j+M})f+\sum_jD_jS_{j-M-1}f.
    \end{equation}
    For the first sum, cancellation and Lemma~\ref{lem:approx}(iii) imply
    \begin{equation}\label{eq:RM-first}
        \|\partial_{\widehat x}^{\alpha}[D_j(I-S_{j+M})f]\|_\infty
        \leq C2^{-Ms}2^{j(|\alpha|-s)}[f]_s,
        \qquad |\alpha|\leq1.
    \end{equation}
    For the second sum, its kernel has integral zero and support at
    coordinate distance at most $C2^{M-j}$. Subtracting the value at $x$
    in the large-scale factor gives
    \begin{align*}
        D_jS_{j-M-1}(x,y)
        &=\int_{\mathbb R^n} D_j(x,u)[S_{j-M-1}(u,y)-S_{j-M-1}(x,y)]\,dw_{k,a}(u),\\
        |D_jS_{j-M-1}(x,y)|
        &\leq\frac{C2^{-M}}{w_{k,a}(B_a(x,2^{M-j}))}
        \mathbf1_{\{\|\widehat x-\widehat y\|\leq C2^{M-j}\}},
        \end{align*}
     and \begin{align*} |\nabla_{\widehat x}(D_jS_{j-M-1})(x,y)|
        &\leq\frac{C2^j}{w_{k,a}(B_a(x,2^{M-j}))}
        \mathbf1_{\{\|\widehat x-\widehat y\|\leq C2^{M-j}\}}.
    \end{align*}
    Consequently,
    \begin{align}
        \|D_jS_{j-M-1}f\|_\infty
        &\leq C2^{-M(1-s)}2^{-js}[f]_s,\notag\\
        [D_jS_{j-M-1}f]_1
        &\leq C2^{Ms}2^{j(1-s)}[f]_s,\notag\\
        [D_jS_{j-M-1}f]_{1/2}
        &\leq C2^{-M(1/2-s)}2^{j(1/2-s)}[f]_s.
        \label{eq:RM-second}
    \end{align}
    The last line follows from
    $[u]_{1/2}\leq C\|u\|_\infty^{1/2}[u]_1^{1/2}$.
    The case $x=x'$ is immediate. For $x\ne x'$, put
    $h_0=\|\widehat x-\widehat{x'}\|>0$ and choose $j_0\in\mathbb Z$
    such that $2^{-j_0}\leq h_0<2^{-j_0+1}$. Use supremum bounds for $j\geq j_0$
    and first order or half-order differences for $j<j_0$. Equations
    \eqref{eq:RM-first}--\eqref{eq:RM-second} yield
    \begin{align}
        |R_Mf(x)-R_Mf(x')|
        &\leq C[f]_s\left[
        2^{-Ms}\left\{\sum_{j\geq j_0}2^{-js}
        +h_0\sum_{j<j_0}2^{j(1-s)}\right\}\right.\notag\\
        &\quad+2^{-M(1-s)}\sum_{j\geq j_0}2^{-js}\left.+2^{-M(1/2-s)}h_0^{1/2}
        \sum_{j<j_0}2^{j(1/2-s)}\right]\notag\\
        &\leq C\{2^{-Ms}+2^{-M(1/2-s)}\}h_0^s[f]_s.
        \label{eq:RM-smallness}
    \end{align}
    Both exponents of $M$ are positive.

    We also verify uniform boundedness and strong convergence of the
    finite sums. Put $E_j=D_jD_j^M$ within this proof. Cancellation,
    support at scale $C_M2^{-j}$, and the differentiated kernel bounds
    give
    \[
    \|E_jf\|_\infty\leq C_M2^{-js}[f]_s
    \qquad \mathrm{and} \qquad \|\nabla_{\widehat x}E_jf\|_\infty
    \leq C_M2^{j(1-s)}[f]_s.
    \]
    For a finite interval $I$, use the same $h_0>0$ and $j_0$ as above to obtain
    \begin{align*}
        \left|\sum_{j\in I}[E_jf(x)-E_jf(x')]\right|
        &\leq C_M[f]_s\left\{\sum_{j\geq j_0}2^{-js}
        +h_0\sum_{j<j_0}2^{j(1-s)}\right\}\\
        &\leq C_M h_0^s[f]_s
    \end{align*}
    The case $h_0=0$ is immediate.
    The $L^2_{k,a}(\mathbb{R}^n)$ bound follows from square-function synthesis:
    \[
    \left\|\sum_{j\in I}D_j(D_j^Mf)\right\|_2
    \leq C\left(\sum_{j\in I}\|D_j^Mf\|_2^2\right)^{1/2}
    \leq C_M\|f\|_2.
    \]
   For a compactly supported function $f$ that is smooth in the transformed
    coordinates, the fine-scale terms satisfy
    \[
    \|E_jf\|_\infty\leq C_M2^{-j}\|\nabla_{\widehat x}f\|_\infty,
    \quad [E_jf]_1\leq C_M\|\nabla_{\widehat x}f\|_\infty,
    \quad \mathrm{and} \quad [E_jf]_s\leq C_M2^{-j(1-s)}\|\nabla_{\widehat x}f\|_\infty.
    \]
    These seminorms are summable as $j\to+\infty$. At coarse scales,
    homogeneity and the support of $f$ give
    \[
    \|E_jf\|_\infty\leq C_{M,f}2^{jn_{k,a}}\|f\|_1,
    \quad [E_jf]_1\leq C_{M,f}2^{j(n_{k,a}+1)}\|f\|_1,
    \mathrm{and} \quad [E_jf]_s\leq C_{M,f}2^{j(n_{k,a}+s)}\|f\|_1.
    \]
    These are summable as $j\to-\infty$, since $n_{k,a}>0$.
    The $L^2_{k,a}(\mathbb{R}^n)$ tails tend to zero by square-function synthesis.
    The uniform bounds extend strong convergence from compact smooth
    functions to $X_s$.

    For sufficiently large $M$, \eqref{eq:RM-L2} and
    \eqref{eq:RM-smallness} give $\|R_M\|<1$ on both spaces. Hence
    \[
    U_M^{-1}=\sum_{r=0}^\infty R_M^r.
    \]
    The inverses agree on $X_s$. Finally,
    \[
    \sum_j\|D_j^MU_M^{-1}f\|_2^2
    \leq(2M+1)\sum_j\sum_{|h|\leq M}\|D_{j+h}U_M^{-1}f\|_2^2
    \leq C_M\|f\|_2^2.
    \]
    This proves (iii) and (iv), and completes the proof.
\end{proof}

\begin{proposition}\label{prop:cancellative}
    If $T1=T^*1=0$ and $T\in\WBP$, then
    \begin{equation}\label{eq:cancellative-bound}
        \|T\|_{2\to2}\leq C(C_K+\|T\|_{\WBP}).
    \end{equation}
\end{proposition}
\begin{proof}
    We work under \eqref{eq:testing-normalization}. The local extension
    established after Lemma~\ref{lem:Meyer}, using
    \eqref{eq:tail-extension}, shows that $T$ and $T^*$ act continuously
    from $X_s$ into distributions paired with compact output test functions.

    For compact test functions $f,g$, put $f_1=U_M^{-1}f$, $g_1=U_M^{-1}g$.
    Let $I,J\subset\Z$ be finite intervals. The block estimate gives
    \begin{align}
        &\left|\sum_{j\in I}\sum_{\ell\in J}
        \langle D_\ell TD_j(D_j^Mf_1),D_\ell^Mg_1\rangle\right|\notag\\
        &\quad\leq C\sum_{j\in I}\sum_{\ell\in J}
        2^{-|j-\ell|\theta}\|D_j^Mf_1\|_2\|D_\ell^Mg_1\|_2\notag\\
        &\quad\leq C\left(\sum_{h\in\Z}2^{-|h|\theta}\right)
        \left(\sum_j\|D_j^Mf_1\|_2^2\right)^{1/2}
        \left(\sum_\ell\|D_\ell^Mg_1\|_2^2\right)^{1/2}\notag\\
        &\quad\leq C\|f\|_2\|g\|_2.
        \label{eq:discrete-Young}
    \end{align}
    For noncompact $f_1,g_1$, these are the $L^2_{k,a}(\mathbb{R}^n)$ extensions of the blocks,
    not an assumed bounded extension of $T$.

    To identify the limit, choose compact test functions $u_N\to f_1$ and
    $v_N\to g_1$ in $X_s$, and finite windows $I_N,J_N$ so that
    \[
    f_N=\sum_{j\in I_N}D_jD_j^Mu_N\longrightarrow f,\qquad
    g_N=\sum_{\ell\in J_N}D_\ell D_\ell^Mv_N\longrightarrow g
    \quad\hbox{in }X_s.
    \]
    Each $f_N,g_N$ is compact and smooth. The finite identity and
    \eqref{eq:discrete-Young} imply, uniformly in $N,L$,
    \[
    |\langle Tf_N,g_L\rangle|\leq C\|u_N\|_2\|v_L\|_2.
    \]
    For fixed $L$, use the local continuity of $T$, then use the local
    continuity of $T^*$ with $f$ fixed:
    \begin{align*}
        |\langle Tf,g\rangle|
        &=\left|\lim_{L\to\infty}\lim_{N\to\infty}
        \langle Tf_N,g_L\rangle\right|\\
        &\leq C\|f_1\|_2\|g_1\|_2\leq C\|f\|_2\|g\|_2.
    \end{align*}
Density and the Riesz representation theorem give a bounded
    extension to $L^2_{k,a}(\mathbb R^n)$. The preceding limit argument
    shows that this extension agrees with the original operator on
    $\mathscr D_a$. Removing the normalization
    \eqref{eq:testing-normalization} proves
    \eqref{eq:cancellative-bound}.
\end{proof}

Finally, we consider the general case
$T1,T^*1\in\BMO(\mathbb R^n,dw_{k,a})$.
We use the discrete Calder\'on reproducing formula on the
coordinate metric measure space $(\mathbb R^n,\rho_a,dw_{k,a})$,
where $\rho_a(x,y)=\|\widehat x-\widehat y\|$.
The result below is an application of the general RD-space
theory in \cite[Theorems~4.12 and~4.13]{HMY}.
We use the form in which $D_j$ is the synthesis kernel;
see also the second reproducing formula in
\cite[Theorem~2.13]{THHLL}.

\begin{theorem}\label{thm:discrete-reproduction}
    Let $D_j=S_j-S_{j-1}$ and fix a sufficiently large integer $J_0$.
    Let $\mathcal Q_j$ consist of the pullbacks of coordinate dyadic
    cubes of side length $2^{-J_0}r_j$. There exists a family of
    dual analysis operators $\{\widetilde D_j\}$
    such that
    \begin{equation}\label{eq:discrete-reproduction}
        f(x)=\sum_j\sum_{Q\in\mathcal Q_j}w_{k,a}(Q)
        D_j(x,x_Q)\widetilde D_jf(x_Q).
    \end{equation}
 Here $Q$ is a pullback of a coordinate dyadic cube, and $x_Q\in Q$
 is a sample point. We choose the kernel representatives and sample
 points in the discrete construction so that the kernel estimates
 and the cancellation identity in the analysis variable hold at
 every selected point. In particular,
 \[
 \int_{\mathbb R^n}\widetilde D_j(x_Q,y)\,dw_{k,a}(y)=0
 \qquad\text{for every }j\in\mathbb Z,\quad Q\in\mathcal Q_j.
 \]
 For $f\in L^2_{k,a}(\mathbb R^n)$, the series converges in
 $L^2_{k,a}(\mathbb R^n)$. The formula also holds in the following
 dual sense: for $b\in\BMO(\mathbb R^n,dw_{k,a})$ and every compactly supported
 coordinate-smooth function $g$ with $\int g\,dw_{k,a}=0$,
 \[
 \langle b,g\rangle
 =\sum_{j,Q}w_{k,a}(Q)\widetilde D_jb(x_Q)
 \overline{D_jg(x_Q)}.
 \]
 Thus the identity for $b$ is understood modulo constants,
 after pairing with mean-zero test functions.
 Moreover, for some $0<\tau<1$, with
 $R=\|\widehat x-\widehat y\|$, the dual kernels satisfy
    \begin{enumerate}[label=\textup{(\roman*)},ref=\thetheorem(\roman*),leftmargin=2.2em,labelsep=0.55em,itemsep=6pt,topsep=6pt,parsep=0pt]
        \item\label{eq:dual-kernel-size}
        $\displaystyle
        |\widetilde D_j(x,y)|\leq
        \frac{C}{w_{k,a}(B_a(x,r_j+R))}
        \left(\frac{r_j}{r_j+R}\right)^\tau;
        $

        \item\label{eq:dual-kernel-difference}
        $\displaystyle
        \begin{aligned}[t]
            |\widetilde D_j(x,y)-\widetilde D_j(x,y')|
            &\leq\frac{C}{w_{k,a}(B_a(x,r_j+R))}
            \left(\frac{\|\widehat y-\widehat{y'}\|}{r_j+R}\right)^\tau
            \left(\frac{r_j}{r_j+R}\right)^\tau
        \end{aligned}
        $
        \par\smallskip
        for $\|\widehat y-\widehat{y'}\|\leq(r_j+R)/2$;

        \item\label{eq:dual-cancellation}
        $\displaystyle \int_{\mathbb R^n}\widetilde D_j(x,y)\,dw_{k,a}(x)=0$
        for $dw_{k,a}$-almost every $y\in\mathbb R^n$;

        \item\label{cond:dual-cancellation-y}
        $\displaystyle \int_{\mathbb R^n}\widetilde D_j(x,y)\,dw_{k,a}(y)=0$
        for $dw_{k,a}$-almost every $x\in\mathbb R^n$.
    \end{enumerate}
\end{theorem}

\begin{proof}
We apply the general reproducing theory to
$(\mathbb R^n,\rho_a,dw_{k,a})$, where
$\rho_a(x,y)=\|\widehat x-\widehat y\|$.
By Proposition~\ref{prop:geometry}, this is an unbounded RD-space
with positive lower growth exponent $\min\{n,n_{k,a}\}$.
The support, derivative, and mass estimates in
Lemma~\ref{lem:approx} show that $\{S_j\}$ is an approximation
to the identity with bounded support in the sense of \cite{HMY}.
We use ordinary dyadic cubes in the transformed coordinates,
with side length $2^{-J_0}r_j$ for a fixed sufficiently large
integer $J_0$.
Fix $0<\eta_0<1$ and $0<\beta_0,\delta_0<\eta_0$.
Then \cite[Theorem~4.12]{HMY}, in the last form of its equation~(4.69),
supplies the dual analysis kernels and the reproducing identity
on the homogeneous mean-zero test space
$\mathring{\mathcal G}_0^{\eta_0}(\beta_0,\delta_0)$ of that reference.
The kernel conclusions, through \cite[Theorem~3.11]{HMY}, give
the stated size, second-variable regularity, and the zero-mean property of $\widetilde D_j$ for an exponent $\tau\in(\eta_0,1)$; the volume
denominators are equivalent by doubling.

To justify the BMO version, let $B_*$ be the reference ball
in the definition of this test space. The size bound for a
mean-zero test function $g$ and the telescoping estimate for
BMO averages over concentric coordinate balls give
\[
\left|\int_{\mathbb R^n}(b-b_{B_*})\overline g\,dw_{k,a}\right|
\leq C\|b\|_{\BMO}\|g\|_{\mathcal G(\beta_0,\delta_0)}
\sum_{m=0}^{\infty}(m+1)2^{-m\delta_0}.
\]
Consequently, $b$ defines a continuous functional on
$\mathring{\mathcal G}_0^{\eta_0}(\beta_0,\delta_0)$, independently of
its additive constant. Applying \cite[Theorem~4.13]{HMY}
therefore gives the asserted dual identity. The factor
$\overline{D_jg(x_Q)}$ follows from the real symmetry of $D_j$.
These background reproducing results are independent of the
$T1$ theorem proved here.

To verify the $L^2$ assertion, write $d\mu=dw_{k,a}$ and apply
Lemma~\ref{lem:sample-square} below. Its proof uses only the
kernel estimates and cancellations already obtained, and does
not use the reproducing identity. By that lemma and the real
symmetry of $D_j$,
the analysis map $f\mapsto\{\widetilde D_jf(x_Q)\}_{j,Q}$ and
the synthesis map
$\{c_{j,Q}\}\mapsto\sum_{j,Q}\mu(Q)D_j(\cdot,x_Q)c_{j,Q}$
are bounded, respectively, from $L^2(\mu)$ to
$\ell^2(\{\mu(Q)\}_{j,Q})$ and back to $L^2(\mu)$.
Hence their composition is a bounded operator on $L^2(\mu)$,
and its defining series converges unconditionally in $L^2(\mu)$.
By the cited reproducing theorem this composition is the identity
on compactly supported coordinate-smooth mean-zero functions,
which are dense in $L^2(\mu)$ because $\mu(\mathbb R^n)=\infty$;
it therefore equals the identity on all of $L^2(\mu)$.
\end{proof}

We next establish the sampled estimate used above and in the
paraproduct construction.

\begin{lemma}
For the kernels and sample points in
Theorem~\ref{thm:discrete-reproduction}, every
$f\in L^2_{k,a}(\mathbb R^n)$ satisfies\label{lem:sample-square}
\begin{equation}\label{eq:sample-square}
    \begin{aligned}
        \sum_{j,Q}w_{k,a}(Q)|D_jf(x_Q)|^2+\sum_{j,Q}w_{k,a}(Q)|\widetilde D_jf(x_Q)|^2
        \leq C\|f\|_2^2.
    \end{aligned}
\end{equation}
\end{lemma}

\begin{proof}
    Write $d\mu=dw_{k,a}$ and
    $\rho(x,y)=\|\widehat x-\widehat y\|$. For each $j$, define
    \[
    A_jf(x)=\sum_{Q\in\mathcal Q_j}\mathbf1_Q(x)
    \widetilde D_jf(x_Q).
    \]
    We use representatives for which the stated kernel estimates and
    the cancellation in the second variable hold at every selected
    sample point. The cubes satisfy
    $\rho(x,x_Q)\leq c r_j$ for $x\in Q$, with $c$ independent of $j,Q$.
    The kernel size estimate and doubling therefore give uniform
    Schur bounds for $A_j$, so that $A_j$ and $A_j^*$ are bounded on
    $L^2(\mu)$.
    Fix $\sigma=\tau/4$. For each selected center $u=x_Q$, the row
    $\widetilde D_j(u,\cdot)$ has integral zero and satisfies the
    ordinary coordinate molecule estimates at scale $r_j$, with
    size and regularity exponents $\tau$. Indeed, the stated local
    regularity estimate implies the global molecule difference
    estimate when the two points are close; otherwise the size
    estimate suffices. Applying the proof of
    Lemma~\ref{lem:molecule-AO} with $d$ replaced by $\rho$ gives,
    for any selected centers $u$ at scale $r_j$ and $v$ at scale
    $r_\ell$,
    \[
    \begin{aligned}
        \left|\int_{\mathbb R^n}
        \widetilde D_j(u,y)\overline{\widetilde D_\ell(v,y)}\,d\mu(y)
        \right|\leq C2^{-\sigma|j-\ell|}
        \frac{1}{\mu(B_a(u,s+\rho(u,v)))}
        \left(\frac{s}{s+\rho(u,v)}\right)^\sigma.
    \end{aligned}
\]
Here $s=\max\{r_j,r_\ell\}$. The argument uses only the coordinate annular and
convolution estimates, which follow from doubling. For the
smaller-scale row, one subtracts from the larger-scale row its
value at the center of the smaller-scale row and uses the zero
integral of the smaller-scale row. Thus no regularity
in the first variable of $\widetilde D_j$ is needed.

If $x\in Q\in\mathcal Q_j$ and $z\in P\in\mathcal Q_\ell$,
the kernel of $A_jA_\ell^*$ is the integral above with
$u=x_Q$ and $v=x_P$. Since
$\rho(x,u)+\rho(z,v)\leq 2cs$, doubling yields
\[
|(A_jA_\ell^*)(x,z)|
\leq C2^{-\sigma|j-\ell|}
\frac{1}{\mu(B_a(x,s+\rho(x,z)))}
\left(\frac{s}{s+\rho(x,z)}\right)^\sigma.
\]
Coordinate annular integration in either variable and Schur's
test imply
\[
\|A_jA_\ell^*\|_{2\to2}
\leq C2^{-\sigma|j-\ell|}.
\]
Consequently, for every finitely supported sequence
$\{h_j\}\subset L^2(\mu)$,
\begin{align*}
    \left\|\sum_j A_j^*h_j\right\|_2^2
    &\leq C\sum_{j,\ell}2^{-\sigma|j-\ell|}
    \|h_j\|_2\|h_\ell\|_2\\
    &\leq C\sum_j\|h_j\|_2^2.
\end{align*}
By duality, followed by monotone convergence over finite index
sets,
\[
\sum_j\|A_jf\|_2^2
=\sum_{j,Q}\mu(Q)|\widetilde D_jf(x_Q)|^2
\leq C\|f\|_2^2.
\]
The same argument with $D_j$ in place of $\widetilde D_j$, using
Lemma~\ref{lem:approx}, proves the other term in
\eqref{eq:sample-square}.
\end{proof}

The paraproduct associated with the above reproducing formula has
the following properties.

\begin{lemma}\label{lem:paraproduct}
    For $b\in\BMO(\mathbb R^n,dw_{k,a})$, there is a classical
    coordinate Calder\'{o}n-Zygmund operator $\Pi_b$ such that
    \begin{enumerate}[label=\textup{(\roman*)},ref=\thetheorem(\roman*),leftmargin=2.2em,labelsep=0.55em,itemsep=6pt,topsep=6pt,parsep=0pt]
        \item\label{eq:Pi-properties} $\Pi_b1=b$ modulo constants;
        \item $\Pi_b^*1=0$;
        \item $\displaystyle \|\Pi_b\|_{2\to2}\leq C\|b\|_{\BMO}$.
    \end{enumerate}
    It also satisfies Definition~\ref{def:CZ} for any sufficiently
    small positive kernel exponent.
\end{lemma}
\begin{proof}
    Define, initially by finite sums,
    \begin{equation}\label{eq:Pi-definition}
        \Pi_bf(x):=\sum_j\sum_{Q\in\mathcal Q_j}w_{k,a}(Q)
        D_j(x,x_Q)\widetilde D_jb(x_Q)S_jf(x_Q).
    \end{equation}
Cancellation and BMO oscillation over coordinate annuli give
\begin{equation}\label{eq:b-coeff}
    |\widetilde D_jb(x_Q)|\leq C\|b\|_{\BMO}.
\end{equation}
For a compactly supported coordinate-smooth function $g$ of
integral zero, the support and derivative estimates for $D_j$ give
\[
\sum_{Q\in\mathcal Q_j}w_{k,a}(Q)|D_jg(x_Q)|
\leq C_g\min\{r_j,r_j^{-1}\}.
\]
The fine-scale estimate uses the cancellation of $D_j$,
while the coarse-scale estimate uses the cancellation of $g$.
Together with \eqref{eq:b-coeff}, this proves absolute convergence
of the dual pairing series above.
    We verify the Carleson estimate. For a coordinate ball $B$ of radius
    $r_B$, subtract $b_{2B}$ and write
    \[
    b-b_{2B}=(b-b_{2B})\mathbf1_{2B}
    +(b-b_{2B})\mathbf1_{(2B)^c}.
    \]
    By \eqref{eq:sample-square} and the square-integrable form of the John--Nirenberg inequality,
    \begin{align*}
        &\sum_{j:r_j\leq r_B}\sum_{Q\subset B}w_{k,a}(Q)
        |\widetilde D_j[(b-b_{2B})\mathbf1_{2B}](x_Q)|^2\\
        &\qquad\leq C\int_{2B}|b-b_{2B}|^2\,dw_{k,a}
        \leq Cw_{k,a}(B)\|b\|_{\BMO}^2.
    \end{align*}
    For the far part, the dual kernel decay gives some $\delta>0$ with
    \begin{align*}
        |\widetilde D_j[(b-b_{2B})\mathbf1_{(2B)^c}](x_Q)|
        &\leq C\|b\|_{\BMO}(r_j/r_B)^\delta
        \sum_{m=1}^\infty(m+1)2^{-m\delta},
        \quad x_Q\in B,\ r_j\leq r_B.
    \end{align*}
    At each scale $\sum_{Q\subset B}w_{k,a}(Q)\leq w_{k,a}(B)$; hence
    \begin{equation}\label{eq:Carleson}
        \sum_{j:r_j\leq r_B}\sum_{Q\subset B}w_{k,a}(Q)
        |\widetilde D_jb(x_Q)|^2
        \leq Cw_{k,a}(B)\|b\|_{\BMO}^2.
    \end{equation}
    The same bound holds for dyadic cubes. Indeed, for a coordinate
    dyadic cube $P$, choose a containing coordinate ball of radius
    $C2^{J_0}\ell(P)$. If $Q\in\mathcal Q_j$ and $Q\subset P$, then
    $r_j\leq2^{J_0}\ell(P)$. The ball estimate therefore includes
    all such terms, and doubling bounds the measure of the ball
    by $C_{J_0}w_{k,a}(P)$. Since $J_0$ is fixed, this gives
    the required dyadic Carleson bound.

    The support and size of $S_j$ imply
    $|S_jf(x_Q)|\leq C\inf_{x\in Q}\mathscr Mf(x)$, where $\mathscr M$ is the
      maximal operator defined in Section~\ref{sec:analytic-symbols}. Applying
    \eqref{eq:Carleson} to the maximal dyadic cubes in its level sets,
    \begin{align}
        \sum_{j,Q}w_{k,a}(Q)|\widetilde D_jb(x_Q)|^2|S_jf(x_Q)|^2
        &\leq C\|b\|_{\BMO}^2
        \int_0^\infty2\lambda\,w_{k,a}\{\mathscr Mf>\lambda\}\,d\lambda\notag\\
        &\leq C\|b\|_{\BMO}^2\|f\|_2^2.
        \label{eq:Carleson-embedding}
    \end{align}
    For finite sums, using the Cauchy--Schwarz inequality gives
    \begin{align*}
        |\langle\Pi_bf,g\rangle|
        &\leq\left(\sum_{j,Q}w_{k,a}(Q)|\widetilde D_jb(x_Q)|^2
        |S_jf(x_Q)|^2\right)^{1/2}
        \left(\sum_{j,Q}w_{k,a}(Q)|D_jg(x_Q)|^2\right)^{1/2}\\
        &\leq C\|b\|_{\BMO}\|f\|_2\|g\|_2.
    \end{align*}
    The two square-summable arrays give convergence of the pairings and
    the asserted bounded operator.
    Its kernel is
    \[
    K_{\Pi_b}(x,y)=\sum_{j,Q}w_{k,a}(Q)
    D_j(x,x_Q)\widetilde D_jb(x_Q)S_j(x_Q,y).
    \]
    At a fixed scale, \eqref{eq:b-coeff} and support give
    \[
    \sum_Qw_{k,a}(Q)|D_j(x,x_Q)\widetilde D_jb(x_Q)S_j(x_Q,y)|
    \leq\frac{C\|b\|_{\BMO}}{w_{k,a}(B_a(x,r_j))}
    \mathbf1_{\{\|\widehat x-\widehat y\|\leq Cr_j\}}.
    \]
    Writing $R=\|\widehat x-\widehat y\|$, lower growth yields
    \begin{equation}\label{eq:Pi-classical-size}
        |K_{\Pi_b}(x,y)|\leq
        C\|b\|_{\BMO}\sum_{r_j\geq R/C}\frac1{w_{k,a}(B_a(x,r_j))}
        \leq\frac{C\|b\|_{\BMO}}{w_{k,a}(B_a(x,R))}.
    \end{equation}
    For $h=\|\widehat x-\widehat{x'}\|\leq R/2$, differentiating the
    $D_j$ factor and integrating along the coordinate segment gives
    \begin{equation}\label{eq:Pi-classical-difference}
        |K_{\Pi_b}(x,y)-K_{\Pi_b}(x',y)|
        \leq C\|b\|_{\BMO}\sum_{r_j\geq R/C}
        \frac{h/r_j}{w_{k,a}(B_a(x,r_j))}
        \leq\frac{C\|b\|_{\BMO}}{w_{k,a}(B_a(x,R))}\frac hR.
    \end{equation}
   The corresponding difference estimate in the second variable
follows similarly from the derivative bounds for $S_j(x_Q,y)$ with
respect to the transformed coordinate $\widehat y$.

    For $d=d(x,y)>0$ and $0<\varepsilon_0<\min\{1,n,n_{k,a}\}$,
    \begin{equation}\label{eq:unreflected-in-mixed}
        \frac1{w_{k,a}(B_a(x,R))}
        \leq\frac C{V(x,y,d)}(d/R)^{\min(n,n_{k,a})}
        \leq\frac C{V(x,y,d)}(d/R)^{\varepsilon_0}.
    \end{equation}
    Also, if $h\leq d/2$,
    \[
    \frac1{w_{k,a}(B_a(x,R))}\frac hR
    \leq\frac C{V(x,y,d)}(h/R)^{\varepsilon_0}.
    \]
    Equations \eqref{eq:Pi-classical-size}--\eqref{eq:unreflected-in-mixed}
    prove the kernel conditions in Definition~\ref{def:CZ}.

We finish by verifying the constant tests. Let $\chi_R$ be a smooth
coordinate cutoff equal to one on $B_a(0,R)$ and supported in
$B_a(0,2R)$, with $0\leq\chi_R\leq1$.
The coordinate kernel difference estimates above show that, for
each fixed compactly supported mean-zero test, the tail corrections
in the definitions of $\Pi_b1$ and $\Pi_b^*1$ are bounded by
$C\|b\|_{\BMO}/R$, where the constant may depend on the test.
Thus these constant tests can be evaluated by the cutoff limits.

For $g\in\mathscr D_a$ with $\int g\,dw_{k,a}=0$, positivity and
normalization give $0\leq S_j\chi_R(x_Q)\leq1$ and
$S_j\chi_R(x_Q)\to1$ as $R\to\infty$.
The absolute pairing estimate proved above permits dominated
convergence in the series, and the dual reproducing identity yields
\[
\lim_{R\to\infty}\langle\Pi_b\chi_R,g\rangle
=\sum_{j,Q}w_{k,a}(Q)\widetilde D_jb(x_Q)
\overline{D_jg(x_Q)}
=\langle b,g\rangle.
\]
Consequently, $\Pi_b1=b$ modulo constants.

For the adjoint identity, take $f\in\mathscr D_a$ with
$\int f\,dw_{k,a}=0$.
For sufficiently large $R$, the support conditions and the
cancellation of $D_j$ make the terms with $r_j<cR$ vanish in
$\langle\Pi_bf,\chi_R\rangle$, where $c>0$ is fixed and small.
For $r_j\geq cR$, subtracting $S_j(x_Q,0)$ and using the
second-variable derivative bound for $S_j$ gives
\[
\sum_Qw_{k,a}(Q)|S_jf(x_Q)|\leq C_f/r_j.
\]
Here the contributing cubes lie in a coordinate ball of radius
$Cr_j$; doubling and the kernel size bound control their total
contribution. Since $\|D_j\chi_R\|_\infty\leq C$,
\eqref{eq:b-coeff} implies
\[
|\langle\Pi_bf,\chi_R\rangle|
\leq C_f\|b\|_{\BMO}\sum_{r_j\geq cR}r_j^{-1}
\leq C_f\|b\|_{\BMO}R^{-1}\longrightarrow0.
\]
Taking complex conjugates gives
$\langle\Pi_b^*\chi_R,f\rangle\to0$, and hence $\Pi_b^*1=0$
in the stated constant-test sense.
\end{proof}

\begin{proof}[Proof of Theorem~\ref{thm:T1-intro}]
    The necessary conditions follow from
    Proposition~\ref{prop:necessity}. For sufficiency, assume that
    $T1,T^*1\in\BMO(\mathbb{R}^n,dw_{k,a})$ and $T\in\WBP$.
    Choose $0<\varepsilon_0<\min\{\varepsilon,1,n,n_{k,a}\}$.
    The kernel conditions at exponent $\varepsilon$ imply those at
    $\varepsilon_0$. Define
    \begin{equation}\label{eq:T0}
        T_0=T-\Pi_{T1}-(\Pi_{T^*1})^*.
    \end{equation}
    By Lemma~\ref{lem:paraproduct}, $T_01=T_0^*1=0$.
    The paraproducts are $L^2_{k,a}(\mathbb{R}^n)$-bounded, so their WBP bounds follow from
    Proposition~\ref{prop:necessity}. Their kernel estimates give
    \begin{align*}
        C_{K_{T_0}}+\|T_0\|_{\WBP}
        &\leq C\bigl(C_K+\|T\|_{\WBP}
        +\|T1\|_{\BMO}+\|T^*1\|_{\BMO}\bigr),\\
        \|T_0\|_{2\to2}
        &\leq C(C_{K_{T_0}}+\|T_0\|_{\WBP}),\\
        \|T\|_{2\to2}
        &\leq\|T_0\|_{2\to2}+\|\Pi_{T1}\|_{2\to2}
        +\|\Pi_{T^*1}\|_{2\to2}\\
        &\leq C\bigl(C_K+\|T\|_{\WBP}
        +\|T1\|_{\BMO}+\|T^*1\|_{\BMO}\bigr).
    \end{align*}
    The second line is Proposition~\ref{prop:cancellative}. This proves
    both sufficiency and \eqref{eq:T1-bound}.
\end{proof}

\subsection*{Acknowledgements}  Y.
Liu was supported by the National Natural Science Foundation of
China (No.12471089, No.12271042), the Beijing Natural Science
Foundation of China (No.1232023)  and the Fundamental Research Funds
for the Central Universities (No.FRF-BRA-26-007).

\vspace{0.3cm} \hspace{-0.5cm}{\bf Competing interests.}
The authors declare that they have no competing interests.

\vspace{0.3cm} \hspace{-0.5cm}{\bf AI declaration.} The authors used OpenAI's Codex for English-language editing
and assistance with checking and refining the exposition and
selected proof details. The authors take full responsibility
for the mathematical arguments and the final text.

\end{document}